\documentclass[11pt]{article}
\usepackage{mlmodern}

\usepackage[letterpaper,margin=1in]{geometry}
\newcommand\blfootnote[1]{%
\begingroup
  \renewcommand\thefootnote{}\footnote{#1}%
  \addtocounter{footnote}{-1}%
  \endgroup
  }

\usepackage[dvipsnames,table,xcdraw]{xcolor}
\usepackage[utf8]{inputenc}

\def\colorful{1}

\ifnum\colorful=1
\usepackage[normalem]{ulem}

\newcommand{\light}[1]{{\grey #1}}
  \newcommand{\anote}[1]{\footnote{{\bf [Ankit: {#1}\bf ]}}}
  \newcommand{\pnote}[1]{\footnote{{\bf [Po-Ling: {#1}\bf ]}}}
  \newcommand{\vnote}[1]{\footnote{{\bf [Varun: {#1}\bf ]}}}
\newcommand{\todo}[1]{{\textbf{ [{\red{Todo}}: {#1}]}}}
\else

\newcommand{\light}[1]{{}}
\newcommand{\anote}[1]{}
\newcommand{\pnote}[1]{}
\newcommand{\vnote}[1]{}
\newcommand{\todo}[1]{}
\fi

\usepackage[T1]{fontenc}    %
\usepackage{url}      %
\usepackage{booktabs}     %
\usepackage{nicefrac}     %

\usepackage{bbm}
\usepackage{caption}
\usepackage{amsmath,amsfonts,amssymb,amsthm}
\usepackage{algorithmic}
\usepackage{algorithm}
\usepackage{array}
\usepackage[caption=false,font=normalsize,labelfont=sf,textfont=sf]{subfig}
\usepackage{textcomp}
\usepackage{stfloats}
\usepackage{url}
\usepackage{verbatim}
\usepackage{graphicx}
\usepackage{enumitem}
\usepackage[
backend=biber,
style=alphabetic,
maxbibnames=15,
maxalphanames=10,
minalphanames=6,
doi=false,
isbn=false,
url=false,
eprint=false,
backref=true,
]{biblatex}
\definecolor{lb}{RGB}{0, 100, 200}
\definecolor{green2}{RGB}{60, 120, 0}

\usepackage[colorlinks,citecolor=green2,linkcolor=lb,bookmarks=true]{hyperref}

\usepackage[nameinlink,capitalise]{cleveref}
\usepackage{thm-restate}

\def\E{\mathbb E}
\def\P{\mathbb P}
\def\R{\mathbb R}

\def\N{\mathbb N}

\newcommand\numberthis{\addtocounter{equation}{1}\tag{\theequation}}

\newcommand{\ba}{\overline{\alpha}}
\newcommand{\bl}{\overline{\lambda}}
\newcommand{\bg}{\overline{\gamma}}
\newcommand{\bp}{\overline{p}}
\newcommand{\bq}{\overline{q}}

\newcommand{\eps}{\epsilon}

\newcommand{\poly}{\mathrm{poly}}

\newcommand{\cN}{\mathcal{N}}

\newcommand{\cE}{\mathcal{E}}

\newcommand{\cF}{\mathcal{F}}
\newcommand{\cP}{\mathcal{P}}
\newcommand{\cY}{\mathcal{Y}}

\newcommand{\cH}{\mathcal{H}}

\newcommand{\cR}{\mathcal{R}}

\newcommand{\cA}{\mathcal{A}}
\newcommand{\cC}{\mathcal{C}}
\newcommand{\cX}{\mathcal{X}}
\newcommand{\cM}{\mathcal{M}}

\newcommand{\cB}{\mathcal{B}}

\newcommand{\dtv}{\mathrm{TV}}
\newcommand{\hel}{\mathrm{h}}
\newcommand{\kl}{\mathrm{KL}}
\newcommand{\js}[1]{\mathrm{JS}_{#1}}
\newcommand{\hbin}{h_{\mathrm{bin}}}
\newcommand{\hent}{H_{\mathrm{ent}}}
\newcommand{\helconst}{0.125}

\newcommand{\prior}{\pi}

\newcommand{\tagm}[1]{\tag*{$\left(#1\right)$}}

\newcommand{\nstar}{{n}^*}
\newcommand{\bay}{\mathrm{B}}
\newcommand{\baye}{\mathrm{B\text{-}E}}
\newcommand{\tot}{\mathrm{total}}
\newcommand{\und}{\mathrm{undetected}}
\newcommand{\era}{\mathrm{erasure}}

\newcommand{\pfht}{\mathrm{PF}}
\newcommand{\pfe}{\mathrm{PF\text{-}E}}
\newcommand{\seq}{\mathrm{SEQ}}
\newcommand{\comm}{\mathrm{comm}}
\newcommand{\priv}{\mathrm{priv}}

\newcommand{\real}{\ensuremath{\mathbb{R}}}
\newcommand{\Ber}{\mathrm{Ber}}

\crefformat{equation}{(#2#1#3)}
\crefname{equation}{Equation}{Equations}
\crefname{lemma}{Lemma}{Lemmata}
\crefname{claim}{Claim}{Claims}
\crefname{fact}{Fact}{Facts}
\crefname{theorem}{Theorem}{Theorems}
\crefname{proposition}{Proposition}{Propositions}
\crefname{corollary}{Corollary}{Corollaries}
\crefname{remark}{Remark}{Remarks}
\crefname{definition}{Definition}{Definitions}
\crefname{question}{Question}{Questions}
\crefname{condition}{Condition}{Conditions}
\crefname{figure}{Figure}{Figures}

\newtheorem{theorem}{Theorem}[section]
\newtheorem{lemma}[theorem]{Lemma}

\newtheorem{claim}[theorem]{Claim}
\newtheorem{proposition}[theorem]{Proposition}
\newtheorem{corollary}[theorem]{Corollary}

\theoremstyle{definition}
\newtheorem{fact}[theorem]{Fact}
\newtheorem{definition}[theorem]{Definition}
\newtheorem{example}[theorem]{Example}
\newtheorem{question}[theorem]{Question}

\newtheorem{remark}[theorem]{Remark}

\begin{document}

\title{The Sample Complexity of Simple Binary Hypothesis Testing\blfootnote{This paper was presented in part at COLT 2024 as an extended abstract~\cite{PenJL24-colt}.}}

\author{
Ankit Pensia
\\
Simons Institute for the Theory of Computing\\
{\tt ankitp@berkeley.edu}\\
\and
Varun Jog\\
University of Cambridge\\
{\tt vj270@cam.ac.uk}
\and
Po-Ling Loh\\ University of Cambridge\\
{\tt pll28@cam.ac.uk}
}

\maketitle

\begin{abstract}%
The sample complexity of simple binary hypothesis testing is the smallest number of i.i.d.\ samples required to distinguish between two distributions $p$ and $q$ in either: (i) the prior-free setting, with type-I error at most $\alpha$ and type-II error at most $\beta$; or (ii) the Bayesian setting, with Bayes error at most $\delta$ and prior distribution $(\prior, 1-\prior)$. This problem has only been studied when $\alpha = \beta$ (prior-free) or $\prior = 1/2$ (Bayesian), and the sample complexity is known to be characterized by the Hellinger divergence between $p$ and $q$, up to multiplicative constants. In this paper, we derive a formula that characterizes the sample complexity (up to multiplicative constants that are independent of $p$, $q$, and all error parameters) for: (i) all $0 \le \alpha, \beta \le 1/8$ in the prior-free setting; and (ii) all $\delta \le \prior/4$ in the Bayesian setting. In particular, the formula admits equivalent expressions in terms of certain divergences from the Jensen--Shannon and Hellinger families. The main technical result concerns an $f$-divergence inequality between members of the Jensen--Shannon and Hellinger families, which is proved by a combination of information-theoretic tools and case-by-case analyses. 
We explore applications of our results to (i) robust hypothesis testing, (ii) distributed (locally-private and communication-constrained) hypothesis testing, (iii) sequential hypothesis testing, and (iv) hypothesis testing with erasures. 
\end{abstract}

\thispagestyle{empty}

\newpage
\thispagestyle{empty}
\tableofcontents
\thispagestyle{empty}

\thispagestyle{empty}
\newpage
\setcounter{page}{1}

\section{Introduction} %
\label{sec:introduction}

Hypothesis testing is a fundamental problem in statistical
inference that seeks the most likely hypothesis corresponding to a given set of observations. 
The simplest formulation of hypothesis testing, which is also the focus of this paper, is %
\emph{simple
binary hypothesis testing}. Here, the two hypotheses correspond to distributions $p$ and $q$ over a domain $\cX$, and the set of observations comprises $n$ i.i.d.\ samples $X_1,\dots,X_n$ from $\theta \in \{p,q\}$. The goal is to identify which distribution generated the samples; i.e., to produce $\widehat{\theta} := \widehat{\theta}(X_1,\dots,X_n)$ so that $\widehat{\theta}= \theta$ with high probability.

Simple binary hypothesis testing is a crucial building block in statistical inference procedures. 
Consequently, much research has analyzed the best algorithms and their performance~\cite{LehRC86}. The famous Neyman--Pearson lemma~\cite{NeyPea33} provides the optimal procedure for $\widehat \theta$, and subsequent works have completely characterized its error probability in two regimes: the single-sample setting ($n=1$) and the infinite-sample setting ($n\to \infty$). While the single-sample regime is relatively straightforward, much historical work in statistics and information theory has focused on the infinite-sample (or asymptotic) regime, where the asymptotic error admits particularly neat expressions in terms of information-theoretic divergences between $p$ and $q$.

Although asymptotic results provide crucial insight into the problem structure, they offer no concrete guarantees for finite samples. In particular, asymptotic bounds cannot satisfactorily answer the question of how many samples are needed (i.e., what is the sample complexity) to solve hypothesis testing with a desired level of accuracy. This limits their applicability in practice, as well as in learning theory research, where sample complexity bounds are paramount.  

While some prior work has established non-asymptotic bounds on the error probability in simple binary hypothesis testing, we observe that the sample complexity perspective has remained largely unaddressed. (See \Cref{sec:related_work} for details.) 
Our main contribution is to fill this gap by developing tight results for the sample complexity of simple binary hypothesis testing. In what follows, we precisely define the problem and introduce terminology necessary to describe our results. 

In the context of simple binary hypothesis testing, an algorithm can make two types of
errors: (i) $\P(\widehat{\theta} = q| \theta =
p)$, termed type-I error, and (ii)
$\P(\widehat{\theta} = p|\theta = q)$, termed
type-II error. 
These two types of errors may have different operational consequences  
in different situations, e.g., false positives vs.\ false negatives of a lethal illness, where the cost incurred by the former is significantly smaller than that of the latter. 
A natural way to combine these metrics is by considering a weighted sum of
the errors, which is equivalent to the
well-studied Bayesian formulation. We define the Bayesian version of simple binary hypothesis testing and its corresponding sample complexity as follows:
\begin{restatable}[Bayesian  simple binary hypothesis testing]{definition}{DefBaySimpHypTesting}
\label{prob:bay-intro}
Let $p$ and $q$ be two distributions over a finite domain $\cX$, and let $\prior, \delta \in (0,1)$.
We say that a 
test $\phi: \cup_{n=0}^\infty \cX^n \to \{p,q\}$
\textit{solves the Bayesian simple binary hypothesis testing
problem} of distinguishing $p$ versus $q$, under the prior
$(\prior, 1- \prior)$, with sample complexity $n$ and probability
of error $\delta$, if for $X := (X_1,\dots,X_n)$, we have
\begin{align}
\label{eq:bayesian-error-prob-intro}
\prior \cdot \P_{X \sim p^{\otimes n}} \left\{   \phi(X) \neq p\right\} + (1 - \prior) \cdot \P_{X \sim q^{\otimes n}} \left\{   \phi(X) \neq q\right\}  \leq \delta\,.
\end{align}
We use $\cB_\bay(p,q,\prior,\delta)$ to denote this problem, and define its sample complexity $\nstar_\bay(p,q,\prior,\delta)$
to be the smallest $n$ such that there exists a test $\phi$ which
solves $\cB_\bay(p,q,\prior,\delta)$.
\end{restatable}
Note that the left-hand side of \Cref{eq:bayesian-error-prob-intro} may also be written as
$\P( \widehat{\Theta} \neq \Theta)$, where $\Theta$ is supported on
$\{p,q\}$ with $\P(\Theta=p)=\prior$, $\widehat{\Theta} = \phi(X)$, and for any $n\in \N$, the
conditional distribution is $X|{\Theta = \theta} \sim \theta^{\otimes n}$.

We also consider a prior-free version, where 
these two errors are analyzed separately.
\begin{restatable}[Prior-free simple binary hypothesis testing]{definition}{DefPriorFreeSimpHypTesting}
\label{prob:pfht-intro}
Let $p$, $q$, and $\phi$ be as in \Cref{prob:bay-intro}. We say that $\phi$ \emph{solves the simple binary
hypothesis testing problem} of $p$ versus $q$ with type-I error
$\alpha$, type-II error $\beta$, and sample complexity $n$, if
\begin{align}
\P_{X \sim p^{\otimes n}} \left\{   \phi(X) \neq p\right\} \leq \alpha, \qquad \text{and} \qquad  \P_{X \sim q^{\otimes n}} \left\{   \phi(X) \neq q\right\}  \leq \beta.
\end{align}
We use $\cB_\pfht(p,q,\alpha,\beta)$ to denote this
problem, and define its sample complexity $\nstar_\pfht(p,q,\alpha,\beta)$ to be the smallest $n$ such that
there exists a test $\phi$ which solves
$\cB_\pfht(p,q,\alpha,\beta)$.%
\end{restatable}

In the remainder of this section, we restrict our attention to the Bayesian setting (\Cref{prob:bay-intro}) for  ease of presentation, but analogous claims hold for the prior-free setting, as well.

\looseness=-1
It is well-known that the likelihood ratio test 
(with a suitable threshold depending on $\alpha$) minimizes the average error probability~\cite{NeyPea33}.
A natural question is to understand the relationship between the Bayes error (the left-hand side of \Cref{eq:bayesian-error-prob-intro} when $\phi$ is the optimal estimator) and the sample size $n$ for the optimal estimator.
The
Chernoff--Stein lemma gives the precise \emph{asymptotic} dependence of the Bayes error on $n$, stating that
the Bayes error equals $e^{-n \left( \mathrm{CI}(p,q) + o(1) \right)}$, where $\mathrm{CI}(p,q)$
is the Chernoff information
between $p$ and $q$~\cite{CovTho06}.
Observe that the asymptotic error rate does not depend on the prior
$\alpha$ at all, summarized as follows in \cite[Chapter
11]{CovTho06}: ``\emph{Essentially, the effect of the prior is washed
out for large sample sizes.}''
Consequently, the error rate is symmetric in $p$ and $q$.
Furthermore, the asymptotic rate suggests a guess for the sample
complexity to obtain error $\delta$:
\begin{align}
\label{eq:guess-asymptotic}
\nstar_\bay(p,q, \prior, \delta) \stackrel{?}{\approx} \frac{\log(1/\delta)}{\mathrm{CI}(p,q)} \asymp \frac{\log(1/\delta)}{\hel^2(p,q)},
\end{align}
where $\hel^2(p,q)$ is the Hellinger divergence between
$p$ and $q$, which satisfies $\hel^2(p,q) \asymp \mathrm{CI}(p,q)$~\cite[Chapter 16]{PolWu23}.

Interestingly, the sample complexity suggested by \Cref{eq:guess-asymptotic} turns out to be accurate for the special case of the uniform prior ($\prior=1/2$). More precisely, \cite{Ziv02} states the sample complexity satisfies
$\nstar_\bay(p,q,1/2, \delta) \asymp
\frac{\log(1/\delta)}{\hel^2(p,q)}$ under mild technical assumptions such as $\hel^2(p,q) \le 0.5$ and $\delta \leq 1/4$. %
We consider this to be a success story for asymptotic results:
precise asymptotics yielding accurate predictions for the
finite-sample regime.

However, this success is rather brittle, and gaps appear between the
asymptotic and non-asymptotic regimes as soon as $\prior \rightarrow 0$.
Consider the simple example where both $p$ and $q$ are Bernoulli
distributions, with biases $0$ and $\epsilon \in (0,1/4)$, respectively. 
Suppose that for a given prior $\prior \in (0,1/2)$, we ask for the
error probability $\delta:= \prior/10$, i.e., a non-trivial
prediction. Direct calculations show that $\nstar_\bay(p,q,\prior, \prior/10)
\asymp \frac{\log(1/\prior)}{\epsilon}$, whereas 
$\nstar_\bay(q,p,\prior, \prior/10) \asymp \frac{1}{\epsilon}$. 
In contrast, the sample complexity suggested by \Cref{eq:guess-asymptotic} for both problems is $\asymp \frac{\log(1/\prior)}{\epsilon}$, which does not correctly capture the sample complexity for the second problem.
This example highlights that the sample complexity may be (i)
\emph{asymmetric} in $p$ and $q$, (ii) \emph{prior-dependent}, and
(iii) \emph{much smaller} than predicted, particularly, {without the
$\log(1/\delta)$ factor in \Cref{eq:guess-asymptotic}; 
all of these contradict the predictions from the asymptotic regime (cf.\ \Cref{eq:guess-asymptotic}).

For general $\prior\in(0,1/2]$, \cite{Ziv02} implies the following non-asymptotic bounds:
\begin{align}
\label{eq:bound-literature}
\frac{\log(\prior/\delta)}{\hel^2(p,q)} \lesssim
\nstar_\bay(p,q,\prior,\delta) \lesssim \frac{\log(1/\delta)}{\hel^2(p,q)}\,.
\end{align}
The lower and upper bounds are not tight up to universal constants, and are clearly seen to diverge when $\prior \to 0$ and
$\delta/\prior$ is relatively large.\footnote{
For example, when $\delta = \prior/10$ and $\prior \to 0$:
As a sanity check, for any fixed $\prior$, as 
$\delta \to 0$, the non-asymptotic bounds in 
\Cref{eq:bound-literature}	converge to the asymptotic prediction
of \Cref{eq:guess-asymptotic}.
In fact, the asymptotic regime kicks in as soon as 
$\delta \leq \prior^2$. \label{footnote: delta}
}
The previous example involving the Bernoulli distributions
demonstrates that \Cref{eq:bound-literature} 
cannot be improved in general; i.e., it is the best possible bound in terms of the Hellinger divergence.
The small-$\prior$ regime corresponds to the high-stakes
situations mentioned earlier.
Thus, a testing procedure  informed by \Cref{eq:bound-literature} 
might need to collect more samples than needed, leading to the following question:
\begin{question}
\centering
\label{ques:sample-complexity}
What is $\nstar(p,q,\prior,\delta)$ as a function of $p$, $q$, $\prior$, and $\delta$? 
\end{question}
We remark that we are looking for a 
``simple'' expression for 
$\nstar_\bay(p,q,\prior,\delta)$ in terms of $p$, $q$, $\prior$, and $\delta$ (tight up to multiplicative constants).
It is not hard to see that the sample complexity is dictated by how fast a certain $f$-divergence between the product
distributions, $p^{\otimes n}$ and $q^{\otimes n}$, grows with $n$.\footnote{For the uniform prior, the ``divergence'' is measured in terms of the total variation distance, which satisfies a two-sided polynomial relationship with the Hellinger divergence (and the Hellinger divergence tensorizes).
For a non-uniform prior, the ``divergence'' is measured 
by the $E_\gamma$-divergence (\Cref{def:E-gamma}),
which does not satisfy any two-sided inequalities with any standard tensorizing $f$-divergence. This is because $E_\gamma$ can be zero for distinct distributions.
}
The divergence between $p^{\otimes n}$ and $q^{\otimes n}$ need not admit a simple expression in terms of $p$ and $q$. Somewhat surprisingly, the Hellinger divergence neatly
captures the sample complexity under a uniform prior. Unfortunately, the proof technique for the uniform prior breaks down for a general prior. We observe that the kind of simple expressions sought in \Cref{ques:sample-complexity} can be extremely beneficial in algorithm design, allowing the statistician to quickly identify the most sample-efficient test among $\cB(p,q,\prior,\delta)$ and $\cB(p',q',\prior',\delta')$ by comparing  $\nstar_\bay(p,q,\prior,\delta)$ with $\nstar_\bay(p',q',\prior',\delta')$.
Existing bounds \Cref{eq:bound-literature} do not even allow for comparisons between $\nstar_\bay(p,q,\prior,\delta)$ and
$\nstar_\bay(q,p,\prior,\delta)$, which can differ by as much as $\log(1/\delta)$.
Moreover, these expressions may highlight additional properties of the sample complexity, which otherwise may be difficult to establish by its definition.
For example, reexamining the uniform prior setting, the Hellinger divergence satisfies additional properties such as joint convexity in $p$ and
$q$ (being an $f$-divergence).
These properties have been crucial in developing efficient
algorithms for hypothesis testing under the
additional constraints of communication and privacy (with a uniform prior)~\cite{PenJL24, PenAJL23}. This leads to our second question:
\begin{question}
\label{ques:structural-properties}
\centering
Do the techniques developed for the uniform prior setting\\ transfer to the general prior setting?
\end{question}
We provide affirmative answers to both questions in \Cref{sec:our_results}.

\section{Our results}
\label{sec:our_results}
\looseness=-1
In this section, we present our main results. \Cref{sub:sample-complexity-of-bayesian,sub:sample-complexity-of-prior} provide tight sample complexity results for the Bayesian and prior-free simple binary hypothesis testing problems (\Cref{prob:bay-intro,prob:pfht-intro}). \Cref{sub:distributed-simple-binary,sub:results-robust-binary} focus on applications to distributed and robust hypothesis testing. \Cref{sub:results-large-failure} considers the weak detection regime.

\subsection{Sample complexity of Bayesian simple hypothesis testing}
\label{sub:sample-complexity-of-bayesian}

We begin with some notation. For a parameter $\lambda\in[0,1]$, define the $f$-divergence $\cH_{\lambda}(p,q) := 1 - \sum_{i}p_i^\lambda 
q_i^{1 - \lambda}$, which is equal to the
Hellinger divergence when $\lambda=1/2$.
We use $I(X;Y)$ to denote the
mutual information between random variables $X$ and $Y$.
A particularly important quantity will be $I(\Theta;X_1)$, 
with $\Theta$ and $X_1$ as in \Cref{prob:bay-intro}.
Observe that $I(\Theta;X_1)$ depends on $\prior$, $p$, and $q$.

Our main result is \Cref{thm:main-result-intro-bay}, which 
characterizes $\nstar_\bay(p,q,\prior,\delta)$ in terms of 
$I(\Theta;X_1)$ and $\cH_{\lambda}(p,q)$.
As mentioned earlier (in~\Cref{footnote: delta}), the sample complexity crucially depends 
on how small $\delta$ is compared to $\prior$.
The following result handles these regimes separately:

\begin{restatable}[Bayesian simple hypothesis testing]{theorem}{ThmMainResultIntroBay}
\label{thm:main-result-intro-bay}
Let $p$ and $q$ be two arbitrary discrete distributions satisfying
$\hel^2(p,q) \leq \helconst$.	%
Let $\prior \in (0,1/2]$ be the prior parameter and $\delta \in
(0,\prior/4)$ be the desired average error probability.
\begin{enumerate}[noitemsep]
\item (Linear error probability) If $\delta \in [\prior/100, \prior/4]$, then for $\lambda = \frac{(0.5 \log 2)}{\log(1/\prior)}$, we have
\begin{align}
\label{eq:main-thm-constant-failure}
\nstar_\bay(p,q,\prior,\delta) 
\asymp \frac{\prior \log(1/\prior)}{I(\Theta;X_1)}
\asymp \frac{1}{\cH_{1 - \lambda} (p,q)}\,.
\end{align}
\item (Sublinear error probability) If $\delta \in \left[\prior^2, \frac{\prior}{100}\right]$, then for $T = \big\lfloor \frac{\log\left( \prior/\delta\right)}{ \log 8} \big\rfloor$ and $\prior' := \prior^{\frac{1}{T}}$, we have
\begin{align}
\nstar_\bay(p,q,\prior,\delta) &\asymp \log\left( \prior/\delta\right) \cdot \nstar_\bay(p,q,\prior',\prior'/4) \,,
\end{align}
where the last expression can be evaluated using \Cref{eq:main-thm-constant-failure}.

\item (Polynomial error probability) If $\delta \leq \prior^2$, then $\nstar_\bay(p,q,\prior,\delta) \asymp \frac{\log(1/\delta)}{h^2(p,q)}$.
\end{enumerate}
\end{restatable}
\Cref{thm:main-result-intro-bay} gives the sample complexity for all $p$, $q$, $\prior$, and $\delta$ up to constant multiplicative factors, answering \Cref{ques:sample-complexity}. We explain the intuition for each regime below:

\begin{enumerate}[itemsep=1pt,topsep=4pt]
\item First suppose $\delta \in [\prior/100,\prior/4]$. If $\delta\geq \prior$, the problem is trivial (we
can simply output ``$q$'' without observing any samples); thus, this regime
requires a probability of error that is a constant
fraction better than the trivial test.
We consider this to be the most insightful regime.
The result \Cref{eq:main-thm-constant-failure} states that the
sample complexity is $\asymp \frac{\hbin(\prior)}{I(\Theta;X_1)}$,
where $\hbin(\cdot)$ is the binary entropy function.
A loose interpretation is as follows: initially,
the uncertainty of $\Theta$ is $\hbin(\prior)$ and each sample
provides roughly $I(\Theta;X_1)$ amount of information.
This argument, when made rigorous, is simply Fano's inequality, leading to a lower bound on the sample complexity in terms of this expression. Proving a matching upper bound is rather
sophisticated and constitutes a key technical result of our work.
\item \looseness=-1 In the sublinear regime of $\delta$, we follow a reduction-based 
approach to show that solving $\cB_B(p,q,\prior,\delta)$ is roughly equivalent 
to solving $T$ independent instances of $\cB_B(p,q,\prior',\delta')$
for $\prior' = \prior^{1/T}$ and $\delta'= \delta^{1/T}$.
The sample complexity upper bound follows from the popular \emph{success-boosting} procedure,
i.e., going from a larger error probability of $\delta'$ to a
smaller error probability $\delta$, using the median output of the $T$ tests.
Somewhat surprisingly, the lower bound by a straightforward application of Fano's inequality turns out to be loose.
The tight lower bound is established using a new \emph{failure-boosting} argument, which shows that the median is near-optimal.
We can choose $T$ so the modified
instance $\cB_B(p,q,\prior', \delta')$ lies in the regime of
\Cref{eq:main-thm-constant-failure}.
\item The final regime of $\delta$ corresponds to the
asymptotic regime, and the result follows by
\Cref{eq:bound-literature}.
\end{enumerate}

\subsection{Sample complexity of prior-free simple hypothesis testing}
\label{sub:sample-complexity-of-prior}

\Cref{thm:main-result-intro-bay} can be used to find the sample complexity of the prior-free formulation of hypothesis testing from \Cref{prob:pfht-intro}, as shown in the following theorem: 

\begin{restatable}[Prior-free simple hypothesis testing]{theorem}{ThmPriorFreeSBHT}
\label{thm:prior-free-simple-hyp-testing-intro}
Let $p$ and $q$ be two arbitrary distributions satisfying 
$\hel^2(p,q) \leq \helconst$.
Let $\alpha \in (0,1/8]$ be the type-I error and $\beta \in (0,1/8]$ be the type-II error.
Without loss of generality, let $\alpha \le \beta$.
Then
$\nstar_\pfht\left( p,q,\alpha,\beta \right) \asymp \nstar_\bay\Big( p, q, \frac{\alpha}{\alpha+
\beta}, \frac{\alpha \beta}{\alpha + \beta} \Big)$,
where the last expression can be calculated using \Cref{thm:main-result-intro-bay}. 
\end{restatable}

}

Thus, \Cref{thm:main-result-intro-bay} also implies tight sample complexity results for the prior-free setting. 

\subsection{Distributed simple binary hypothesis testing}
\label{sub:distributed-simple-binary}
The sample complexity formula established in \Cref{thm:main-result-intro-bay} (cf.\ \Cref{ques:structural-properties}) has numerous consequences. Recall that \Cref{thm:main-result-intro-bay} characterises the sample complexity in terms of either the mutual information
or the Hellinger divergence. Both structures have favorable properties and lead to algorithmic and statistical results in the distributed setting, as we show
below.

Distributed  statistical inference is typically studied in the following framework: Samples are collected by $n$  observers, who then communicate information about their samples (such as summary statistics) to a central server, which tries to solve the inference problem. The transmission from  an observer to the central server may be subject to different constraints, such as a bandwidth constraint that limits the number of bits the observer might transmit, or a privacy constraint that requires all transmissions to protect the privacy of the observer's data. 
We define the generic problem below\footnote{For simplicity, we
consider the special case where each user employs the same modifier $f$.
}:
\begin{definition}[Distributed simple binary hypothesis setting]
\label{def:hyp-testing-generic}
Let $p$ and $q$ be two distributions over $\cX$ and define the prior
parameter $\prior \in (0,1/2)$ and the error probability $\delta
\in (0,1)$.
Let $\cY$ be a fixed domain and let $\cC$ be a fixed set of
stochastic maps from $\cX \to \cY$, representing the constraints.
We say that a test-modifier pair $(\phi,f)$ with test $ \phi:
\cup_{n=0}^\infty \cX^n \to \{p,q\}$ and $ f\in \cC$
\textit{solves the Bayesian simple binary hypothesis testing
problem} of distinguishing $p$ versus $q$, under the prior
$(\prior, 1- \prior)$, with sample complexity $n$, probability of
error $\delta$, under constraints $\cC$, if for $Y_1,\dots,Y_n$,
generated independently (conditioned on the $X_i$'s) as $Y_i =
f(X_i)$, the following holds:
\begin{align*}
\prior \cdot \P_{X \sim p^{\otimes n}, Y_1,\dots,Y_n} \left\{\phi\left( Y_1,\dots, Y_n \right) \neq p\right\} + (1 - \prior) \cdot \P_{X \sim q^{\otimes n},Y_1,\dots,Y_n} \left\{   \phi\left(Y_1,\dots,Y_n)\right) \neq q\right\}  \leq \delta\,.
\end{align*}
We use $\cB_{\bay,\mathrm{constrained}}
(p,q,\prior,\delta,\cC)$ to denote this problem, and define its
sample complexity to be the smallest $n$ such that there exists a
test-modifier pair $(\phi,f)$ which solves
$\cB_{\bay,\mathrm{constrained}}(p,q,\prior,\delta,\cC)$; we
denote the sample complexity by
$\nstar_{\bay,\mathrm{constrained}}(p,q,\prior,\delta,\cC)$.
\end{definition}
We are interested in understanding
the ``costs'' of constraints $\cC$, measured in terms of statistical
complexity and runtime.
Starting with the statistical cost, we wish to compare
$\nstar_{\bay,\mathrm{constrained}}(p,q,\prior,\delta,\cC)$ and $\nstar_{\bay}(p,q,\prior,\delta)$.
The computational cost refers to the runtime of finding a modifier $f
\in \cC$ such that the resulting sample complexity with the choice of
modifier $f$ is comparable to the optimal choice,
$\nstar_{\bay,\mathrm{constrained}}(p,q,\prior,\delta,\cC)$.\footnote{Once the modifier $f$ is fixed, the optimal test $\phi$ is the likelihood ratio test between the push-forward distribution of $f$ under $p$ and $q$.}
The uniform prior setting was studied extensively in
\cite{PenJL24} for communication
constraints (see also \cite{BhaNOP21}) and \cite{PenAJL23} for local privacy constraints.
At a high level, \cite{BhaNOP21,PenJL24,PenAJL23} characterized the contraction of the
Hellinger divergence under these constraints and the time complexity
of finding a good $f \in \cC$.
Due to the reliance on the Hellinger divergence, a naive application of
these results to an arbitrary $\prior$ and $\delta$ would incur a 
superfluous multiplicative $\log(1/\delta)$ factor in the resulting 
sample complexity, akin to \Cref{eq:bound-literature}.
In this section, we show much improved bounds by using the structural
properties of the sample complexity established in \Cref{thm:main-result-intro-bay} 
above.

We begin with communication-efficient simple
binary hypothesis testing, also known as ``decentralized detection''~\cite{Tsitsiklis93}.
The following result is obtained by using the characterization of the sample complexity in terms of mutual information in \Cref{thm:main-result-intro-bay}, with the result of \cite{BhaNOP21} on the contraction of mutual information under quantization:
\begin{restatable}[Statistical and computational costs of communication for hypothesis testing]{theorem}{ThmCommHypTesting}
\label{thm:comm-hyp-testing}
Let $p$ and $q$ be two arbitrary distributions over $\cX$,  and let $\prior \in (0,1/2]$ and $\delta \in (0,\prior/4)$.
Let $D \geq 2$ be an integer denoting the communication budget.
Consider the setting in \Cref{def:hyp-testing-generic} with $\cY
= [D]$ and $\cC$ being the set of all channels from $\cX \to \cY$, and denote the resulting
sample complexity by
$\nstar_{\bay,\comm}(p,q,\prior,\delta,D)$.
Let $\nstar := \nstar_{\bay}(p,q,\prior,\delta)$.
For any integer $D\geq 2$, we have the following:
\begin{enumerate}
\item (Statistical cost) $\nstar \lesssim \nstar_{\bay,\comm}(p,q,\prior,\delta,D) \lesssim \nstar \max\left(1 , \frac{\log(\nstar/ \prior)}{D}  \right)$. 

\item (Computational cost)
There exists an algorithm that takes as input $p, q,
\prior, \delta$, and $D$ and runs in time
$\poly(|\cX|,D,\log(1/\prior),\log(1/\delta))$, and
outputs a quantizer $\widehat{f} \in \cC$ such that the
resulting sample complexity is at most $\nstar
\max\left(1 , \frac{\log(\nstar/ \prior)}{D} \right)$.
\end{enumerate}
\end{restatable}
Observe that for large enough $D$, the sample complexity is comparable to the unconstrained sample complexity.
A naive application of \cite{BhaNOP21,PenJL24} with
\Cref{eq:bound-literature} would instead yield a sample complexity of  $ n' \cdot \max\left(1, \frac{\log(n')}{D}\right)\cdot\log(1/\delta)$, for $ n':= \nstar_{B}(p,q,1/2, 1/8) \asymp 1/\hel^2(p,q)$, which can be larger than $\nstar$ by a factor of $\log(1/\delta)$ for all values of $D$.

It is likely that a more whitebox analysis of \cite{BhaNOP21} can replace $\log(\nstar/\prior)$ from \Cref{thm:comm-hyp-testing} with $\log(\nstar)$.
 The proof of \Cref{thm:comm-hyp-testing} is given in \Cref{app:hypothesis-testing-under-communication-constraints} and crucially relies on the mutual information characterization of sample complexity in \Cref{thm:main-result-intro-bay}.

We now turn our attention to privacy constraints, starting with a definition:
\begin{definition}[Local differential privacy]
For two finite domains $\cX$ and $\cY$, we say that a randomized mechanism $\cM: \cX \to \cP(\cY)$ is $\epsilon$-LDP for an $\epsilon > 0$ if for all $x,x' \in \cX$ and $y \in \cY$, we have $\P(\cM(x) = y) \leq e^{\epsilon} \cdot \P(\cM(x') = y)$.
\end{definition}

Similar to communication
constraints, the two main goals
are understanding the
statistical and computational costs
of privacy.
Even under a uniform prior, the statistical cost of privacy
(i.e., contraction of the Hellinger divergence under LDP constraints) can be
rather delicate, as shown in \cite{PenAJL23}.
Though one can similarly study the contraction of
$\cH_\lambda(\cdot,\cdot)$ under local privacy, we instead focus on
understanding the computational costs of privacy,
which effectively asks us to find the (approximate) minimizer of
$\nstar(\cM(p), \cM(q), \prior ,\delta)$ over private channels $\cM$.\footnote{Here, $\cM(p)$ denotes the push-forward measure of $p$ under the mechanism $\cM$.}
While it is unclear how to solve such a generic optimization problem in a computationally-efficient manner,
\cite{PenAJL23} gave an efficient algorithm for minimizing $g(\cM(p), \cM(q))$ over quantized private channels when $g(\cdot,\cdot)$ is a jointly quasi-concave function. 
Fortunately, \Cref{thm:main-result-intro-bay} states that $\nstar_\bay$ is also characterized by an $f$-divergence, so \cite{PenAJL23} is applicable. 

\begin{restatable}[Computational costs of privacy for hypothesis testing]{theorem}{ThmCompCostPrivacy}
\label{thm:computational-cost-privacy}
Let $p$ and $q$ be two arbitrary distributions over $\cX$.
Let $\prior \in (0,1/2]$ and $\delta \in (0,\prior/4)$.
Let $\epsilon>0$ be the privacy budget.
Consider the setting in \Cref{def:hyp-testing-generic}, with $\cY
= \cX$ and $\cC$ equal to the set of all $\epsilon$-LDP channels
from $\cX \to \cY$, and denote the resulting
sample complexity by
$\nstar_{\bay,\priv}(p,q,\prior,\delta,\epsilon)$.%

Let $\nstar_{\priv} :=
\nstar_{\bay,\priv}(p,q,\prior,\delta, \epsilon)$.
There exists an algorithm $\cA$ that takes as inputs $p, q,
\prior, \delta$, $\epsilon$, and communication budget parameter
$D$, runs in time
$\poly(|\cX|^{D^2},2^{\poly(D)},\log(1/\prior),\log(1/\delta))$, and
outputs an $\epsilon$-LDP mechanism $\cM:\cX \to [D]$, such that
the resulting sample complexity with $\cM$ is at most
$\nstar_{\priv} \cdot \max \left(1,
\frac{\log (\nstar_{\priv}/\prior)}{D} \right)$.
\end{restatable}
We remind the reader that the resulting sample complexity has
no explicit multiplicative factor of $\log(1/\delta)$ for $D$ moderately large.\footnote{A naive application of \cite{PenAJL23} would yield an $\epsilon$-LDP mechanism with 
sample complexity $\log(1/\delta) \cdot n' \cdot \max\left (1, \frac{\log(n')}{D} \right)$, for $ n':= \nstar_{\bay,\priv}(p,q,1/2, 1/8)$.}
We defer the proof of \Cref{thm:computational-cost-privacy} to \Cref{app:hypothesis-testing-under-local-differential-privacy}.

\subsection{Robust simple binary hypothesis testing}
\label{sub:results-robust-binary}

In some scenarios, it may not be possible to know the
distributions $p$ and $q$ exactly, but only up to a small error.
In this setting, we want to guess ``$p$'' (similarly, ``$q$'') even
if the samples are generated by $p' \in D_1$ (respectively, $q' \in
D_2$), where $D_1$ and $D_2$ are the sets of all distributions within
an $\epsilon$ total variation distance\footnote{Other distance
metrics over distributions could be considered, e.g., using the
Hellinger divergence.
Our results should continue to hold as long as they fit into the
least favorable distributions framework of \cite{Huber65,HubStr73}.
}, respectively.
This setting is termed robust hypothesis testing, and has a rich
history in statistics, with seminal contributions by~\cite{Huber65}.

\begin{definition}[Robust Bayesian simple binary hypothesis testing]
\label{prob:bay-rob}
Consider the setting in \Cref{prob:bay-intro}.
Let $\epsilon < \dtv(p,q)/2$, and define $D_1 := \{p': \dtv(p,p')
< \epsilon\}$ and $D_2 := \{q': \dtv(q,q') < \epsilon\}$.
We say that a (randomized) test $\phi: \cup_{n=0}^\infty \cX^n \to \{p,q\}$
\textit{solves the Bayesian simple binary hypothesis testing
problem} of distinguishing $p$ versus $q$, under the prior
$(\prior, 1- \prior)$, with sample complexity $n$, probability of
error $\delta$, and contamination level $\epsilon$, if
\begin{equation*}
\prior \cdot \sup_{p' \in D_1} \P_{X \sim p'^{\otimes
n}} \left\{ \phi(X) \neq p\right\} + (1 - \prior) \cdot \sup_{q'
\in D_2} \P_{X \sim q'^{\otimes n}} \left\{ \phi(X) \neq
q\right\} \leq \delta.
\end{equation*}
We denote the problem by $\cB_{\bay,
\mathrm{rob}}(p,q,\prior,\delta, \epsilon)$ and its sample
complexity by $\nstar_{\bay,\mathrm{rob}}(p,q,\prior,\delta,
\epsilon)$.
\end{definition}
\cite{Huber65} established that
\Cref{prob:bay-rob} reduces to solving
$\cB_\bay(p',q',\prior,\delta)$ for the least favorable
distribution (LFD) pair $p' \in D_1$ and $q' \in D_2$.
Moreover, the LFD pair can be calculated in a computationally-efficient manner by an algorithm
given $p$, $q$, and $\epsilon$.
Applying our result to the LFD pair yields the following result:

\begin{corollary}[Corollary of \Cref{thm:main-result-intro-bay} and \cite{Huber65}]
Let $p$ and $q$ be two arbitrary distributions, and suppose $\prior \in
(0,1/2]$, $\delta \in (0, \prior/4)$, and $\epsilon \in (0, \dtv(p,q)/2)$.
There is a computationally efficient algorithm that takes $p$, $q$, and
$\epsilon$ as inputs and computes distributions $p'$ and $q'$ such that
$\nstar_{\bay,\mathrm{rob}}\left( p, q, \prior, \delta,
\epsilon \right) \asymp\nstar_\bay\left( p',q' ,\prior, \delta
\right)$,
where the last expression can be evaluated using
\Cref{thm:main-result-intro-bay}.%
\end{corollary}

\subsection{Large error probability regime: Weak detection}
\label{sub:results-large-failure}

Finally, we turn our attention to the large error probability
regime, often called the ``weak detection regime.''
In this regime, we allow the estimators to make an error with probability $
\delta = \prior(1 - \gamma)$, for $\gamma$ close to 0 in $(0,1)$.
Recall that $\delta = \prior$ is trivial, so the regime
$\delta= \prior(1 - \gamma)$ asks the estimator to be
``very slightly'' better than trivial.
Even for a uniform prior, it is rather surprising that the exact sample
complexity is not known, with the tightest known results being
\begin{align}
\frac{\gamma^2}{\hel^2(p,q)}\lesssim 
\nstar_{\bay}\left( p, q, 1/2, \frac{1 -\gamma}{2} \right) 
\lesssim \frac{\gamma}{\hel^2(p,q)}
\end{align}
(see, e.g., \cite[Eq. (14.19)]{PolWu23}). In fact, as we show, these bounds happen to be tight in the worst case,
i.e., there exist $(p',q')$ and $(p'',q'')$ with $\hel^2(p',q') =
\hel^2(p'',q'')$ such that the sample complexity of $(p',q')$ is
given by the lower bound, while that of $(p'',q'')$ is given by the
upper bound (cf.\ \Cref{example:uniform-prior}).
Thus, surprisingly, even for a uniform prior, the Hellinger divergence does
not exactly characterize the sample complexity in the weak detection
regime.
Even more surprisingly, for $\prior \neq 1/2$, the sample complexity
in the weak detection regime may not even depend on $\gamma$ for all small values of $\gamma$. %
We establish the following result:
\begin{theorem}[Surprises in weak detection; Informal version] Consider the setting in \Cref{prob:bay-intro}.
Let $\delta = \prior(1 - \gamma)$, for some $\gamma \in (0,1/2)$.
Then 			
\begin{align*} 
\gamma \max\left\{\gamma,|0.5 - \prior|\right\} \lesssim \frac{\nstar_\bay\left( p, q, \prior, \prior(1 - \gamma) \right)}{\nstar_\bay\left( p, q, \prior, \prior/4 \right)}  
\lesssim \max\left\{\gamma, |0.5 - \prior|\right\}.
\end{align*}
In particular, for $\prior \leq 1/4$, we have $\gamma  \lesssim \frac{\nstar_\bay\left( p, q, \prior, \prior(1 - \gamma) \right)}{\nstar_\bay\left( p, q, \prior, \prior/4 \right)}  
\lesssim 1$, and the upper bounds are tight in the worst case (cf.\ \Cref{ex:upper-bounds-tight-weak-detection}).
\end{theorem}
The formal version of this result is presented in \Cref{thm:bayes-all-regimes}.
\looseness=-1That is, the upper bounds and the lower bounds differ by a 
factor of $\gamma$, with the upper bound not decreasing with $\gamma$
for any non-uniform prior.
Perhaps surprisingly, there exist distribution pairs whose sample
complexity with error probability $\prior/4$ is almost as large as the
sample complexity with error probability $\prior(1- \gamma)$, for
any tiny $\gamma$. To be precise, for any $n_0 \in \N$ and  non-uniform prior, there exist distributions, given in \Cref{ex:upper-bounds-tight-weak-detection},
such that the error probability with $n$ samples for all $n\in\{0,\dots,n_0\}$ is $\prior$---the trivial error, whereas $\nstar_\bay(p,q,\prior, \prior/4) \lesssim n_0 $.
Hence, weak detection might be as hard as strong detection (up to constant factors).

\section{Related work}
\label{sec:related_work}

Simple binary hypothesis testing is a classical problem that has been studied for over a century. The famous Neyman--Pearson lemma~\cite{NeyPea33} states that the optimal test is obtained by thresholding the likelihood ratio. The lemma allows calculation of the decision regions, but does not provide the resulting errors of the optimal test. It is of great practical and theoretical interest to determine how many i.i.d.\ samples are necessary to obtain errors that are below some desired threshold, and this problem has accordingly been studied extensively within the statistics and information theory communities. Almost all of this work focuses on asymptotic properties of the errors for large sample sizes; i.e., identifying exponents $E$ such that the error(s) decay as $e^{-n(E+o(1))}$. These results are included in textbooks such  as~\cite{CovTho06, PolWu23}, and we briefly survey them: the \emph{Stein regime} considers the exponential rate of decay of the type-II error when the type-I error is fixed; the \emph{Chernoff regime} considers all possible pairs $(E_0, E_1)$ such that the type-I and type-II errors decay exponentially fast with exponents $E_0$ and $E_1$; and the \emph{Bayesian setting} considers the exponential rate of decay of the Bayes error. In all of these problems, the exponents are neatly characterized in terms information-theoretic quantities such as the KL divergence and the Chernoff information between $p$ and $q$. Non-asymptotic bounds on the errors have appeared in prior works such as \cite{Stra62, PolPooVer10}, using Berry--Esseen type bounds to analyse the concentration of the likelihood ratio. These works derive upper and lower bounds on the type-II error in terms of the type-I error and the number of samples $n$. Deriving analytical sample complexity from these results appears difficult. We may derive sample complexity bounds computationally, but these are unlikely to be tight within constants since they rely on moments of the likelihood ratio and not on the divergences identified in this paper. The problem of deriving sample complexity bounds that are tight up to multiplicative constants gained traction in the early 2000s with the emergence of property testing~\cite{Gold17}, and in particular, distribution testing~\cite{Can22}, within the theoretical computer science community. To our knowledge, the sample complexity of simple binary hypothesis testing was first explicitly stated in~\cite{Ziv02}, although it was likely well-known before that. Here, it was shown that the sample complexity is $\asymp \log(1/\delta)/\hel^2(p,q)$, where $\delta$ is the desired type-I and type-II error or the Bayes error. Distributed hypothesis testing under information constraints dates back to the work of~\cite{Tsitsiklis93} which studied an asymptotic version of hypothesis testing under communication constraints. A relatively recent resurgence in interest has led to numerous contributions towards obtaining sample complexity bounds for distribution estimation under communication and privacy constraints~\cite{AchCT20-I, AchCT20-II,CheKO23}. The problem of simple hypothesis testing (binary and $M$-ary) under constraints has been studied in~\cite{CanKMSU19,BunKSW19,GopKKNWZ20,PenJL24,PenAJL23,AliBS23}.

\section{Preliminaries} %
\label{sec:preliminaries}

\paragraph{Notation:} Throughout this paper, we assume that the domain $\cX$ is discrete, taken
to be $\mathbb N$ without loss of generality.
For distributions $p_1,\dots,p_n$ over $\cX$, we use $\prod_{i=1}^n
p_i$ to denote the product distribution.
When $p_1,\dots,p_n$ are identically equal to a distribution $p$, we
denote the joint distribution by $p^{\otimes n}$.
When we have a setting with just two distributions $p$ and $q$, we
use $p_i$ and $q_i$ to denote $p(i)$ and $q(i)$, respectively.
We denote the measure of a set $A \subseteq \cX$ under $p$ by
$\mathbb P_p(A)$.
For a bias $ a \in [0,1]$, we denote the Bernoulli distribution with bias $a$ by
$\Ber(a)$.

For two distributions $p$ and $q$, we use
$\dtv(p,q) = 0.5\E_q\left[\left|\frac{dp}{dq} -
1\right|\right]$ to denote the total variation
distance; $\hel^2(p,q) = 0.5
\E_q\left[\left|\sqrt{\frac{dp}{dq}} -
1\right|^2\right]$ to denote the Hellinger
divergence;
and $\kl(p,q) = \E_q\left[\frac{dp}{dq} \log\left(
\frac{dp}{dq} \right)\right]$, whenever finite, to
denote the Kullback-Leibler divergence, where
$\frac{dp}{dq}$ denotes the Radon--Nikodym
derivative.
We define additional $f$-divergences in
\Cref{sec:probability_divergences}.
For two random variables $X$ and $Y$, we use $I(X;Y)$ to denote the
mutual information between $X$ and $Y$.

For a scalar $x \in [0,1]$, shall use
$\overline{x}$ as shorthand to denote $1-x$.
All logarithms are natural logarithms. We use $\poly(\cdot)$ to denote an expression that is polynomial in this arguments.
We also use the standard inequalities $\lesssim, \gtrsim, \asymp$ to hide absolute constant factors in the relationships.

Let $P_e^{(n)}$ denote the Bayesian error probability with $n$ i.i.d.\ samples, i.e., the minimum error probability over all tests $\phi$ in \Cref{eq:bayesian-error-prob-intro}.
When $n = 1$, we also use the notation $P_e$.

Finally, we will focus on pairs of distributions $p$ and $q$ that are sufficiently far apart, in the sense of $\nstar_\pfht\left(p,q,\frac{1}{4}, \frac{1}{4}\right) \geq 2$. 
This condition implies that one needs to observe at least two samples to ensure that both type-I and type-II errors are at most $1/4$. We consider this to be a fairly mild regularity condition.
If such a condition does not hold, then the sample complexity might display some degenerate behavior, contradicting even \Cref{eq:bound-literature}.
For example, if $p$ and $q$ have non-overlapping supports, then $\nstar_\bay(p,q,\prior, \delta) \leq 1$ and $\nstar_\pfht(p,q,\alpha, \beta) \leq 1$ for all parameter choices $\alpha, \delta, \beta, \prior$.
We note that a non-trivial upper bound on $\hel^2(p,q)$ is also assumed in \cite[Theorem 4.7]{Ziv02}.
The condition $\nstar_\pfht(p,q,1/4,1/4) \geq 2$ is guaranteed if $\hel^2(p,q) \leq \helconst$ because $\nstar_\pfht(p,q,1/4,1/4) \geq \nstar_\bay(p,q,0.5,1/4)$ by \Cref{claim:reln-bayesian-prior-free}, and $\nstar_\bay(p,q,0.5,1/4) \geq 2$ if $0.5(1 - \dtv(p,q)) < 1/4$ by \Cref{prop:E-gamma}, which is equivalent to $\dtv(p,q) > 1/2$. 
Finally, $\dtv(p,q) < \sqrt{2} \hel(p,q)$ for distinct $p$ and $q$ (see, for example, \cite[Proposition 2.38]{Ziv02}) implies that $\hel(p,q) \leq 1/2 \sqrt{2}$ suffices for $\dtv(p,q) < 1/2$.

\subsection{Problem definitions}

\subsubsection{Bayesian hypothesis testing} %
\label{par:bayesian_hypothesis_testing}

Recall the Bayesian hypothesis
testing problem:

\DefBaySimpHypTesting*

Observe that the problem is trivial when $\delta
\geq \min(\prior,1- \prior)$.
The case $\prior= 1/2$ corresponds to the typical setting of a
uniform prior, which is well-understood;
in particular, the sample complexity with a uniform prior is
characterized by the Hellinger divergence.
We record both of these results below:

\begin{proposition}[Standard bounds on $\nstar_\bay$ (Folklore)]
We have the following bounds:
\label{prop:trivial-bounds-bay}
\begin{enumerate}
\item (Vacuous Regime) $\nstar_\bay(p,q,\prior,\delta) = 0$ if $\delta \geq \min(\prior,1-\prior)$.
\item Let $p$ and $q$ satisfy $\hel^2(p,q) \leq \helconst$.
 For all $\delta \leq \prior/4$, we have
$\frac{\log(\prior/\delta)}{\hel^2(p,q)} \lesssim  
\nstar_\bay(p,q,\prior,\delta) \lesssim
\frac{\log(1/\delta)}{\hel^2(p,q)}$.
Moreover, these bounds are tight in the worst case, as shown in \Cref{ex:bern-symmetry}.

\item In particular, for $\prior = 1/2$, any $\delta \leq
1/8$, and $\hel^2(p,q) \leq \helconst$, we have
$\nstar_\bay(p,q,0.5,\delta) \asymp
\frac{\log(\frac{1}{\delta})}{\hel^2(p,q)}$~\cite{Ziv02,PolWu23}.
\end{enumerate}
\end{proposition}
The first claim in the proposition follows by always picking the output
hypothesis (without even looking at the samples) with prior $1 -
\prior$ (respectively, $\prior$), leading to a test with average
error $\prior$ (respectively, $1- \prior$).
The second claim follows by \Cref{prop:trivial-bounds-fre} and
\Cref{cor:equivalence-bayesian-pf}, proved later.

\Cref{prop:trivial-bounds-bay} implies that the ``interesting'' regime 
of the parameter $\delta$ is when $\delta$ is not too small: 
In particular, if $\delta \leq \prior^{1+c}$ for a constant $c > 0$, the sample complexity is $\log(1/\delta)/h^2(p,q)$ (up to a multiplicative factor of $1/c$),
i.e., it matches the uniform prior rates.
This is in line with the asymptotic result that the prior does not
play a role in the asymptotic regime~\cite{CovTho06}.
Our main contribution will be to provide tight bounds when $\delta \geq
\prior^{1+c}$.
A particularly important case will be the linear regime, i.e., the
case when $\delta \in (\prior/100, \prior/4)$.
As we shall see, a generic instance $\cB_\bay(p,q,\prior,\delta)$
with error probability $\delta \in (\prior^2,\prior/4)$
can be related to that of 
$\cB_\bay(p,q,\prior',\delta')$, where $\delta' \in (\prior'/16,
\prior'/4)$ (cf.%
\ \Cref{prop:reduction}).

Moreover, the sample complexity is
not symmetric in the first two
arguments when $\prior
\neq 0.5$, i.e., in general,
$\nstar_\bay(p,q,\prior, \delta)
\neq \nstar_\bay(q,p,\prior,
\delta)$.
For example, for every $\epsilon \in (0, 1)$,
there exist $p$ and $q$ with $\hel^2(p,q) = \epsilon$ such that for
all $\prior \leq 1/2$, we have $\nstar_\bay\left(p,q,\prior,\frac{\prior}{10}\right) \asymp \frac{\log(1/\prior)}{\epsilon}$  and   $\nstar_\bay\left(q,p,\prior, \frac{\prior}{10}\right) \asymp \frac{1}{\epsilon}$.
Formally, we have the following:
\begin{example}
\label{ex:bern-symmetry}
Let $\epsilon \in (0,\eps_0)$ for a small absolute constant $\eps_0$. and $\prior \in (0,1/2)$. 
Let $\delta \in (0,\prior/4)$.
Let $p = \Ber(0)$ and $q = \Ber(\eps)$ be the Bernoulli distributions with bias $0$ and $\epsilon$, respectively.
 Then $\hel^2(p,q) \asymp \epsilon$ and
  \begin{align*}
  \nstar_\bay(p,q,\prior,\delta) \asymp \frac{\log(1/\prior)}{\hel^2(p,q)} \,\,\, \text{and}\,\,  \nstar_\bay(q,p,\prior,\delta) \asymp \frac{\log(\prior/\delta)}{\hel^2(p,q)}.
  \end{align*}
\end{example}

\begin{proof}
Let the support of $p$ and $q$ be $\{A,B\}$, such that $\P(A) = 0$ and $\P(B) = 1$ under $p$.

We first show that $\nstar_\bay(q,p,\prior) \asymp \frac{\log(\prior/\delta)}{\epsilon}$. The lower bound follows from \Cref{prop:trivial-bounds-bay}, so we focus on establishing the upper bound.
Consider the test $\phi$ that picks $q$ if the data contains a nonzero number of $A$'s, otherwise picks $B$.
Then the probability of error under $p$ is $0$, whereas the probability of error under $q$ is $(1 -\epsilon)^n$.
Thus, the expected error is $\prior(1 - \epsilon)^n$, which is less than $\delta$ if $n \lesssim \log(\prior/\delta)/\epsilon$.

We now consider $\nstar_\bay(p,q,\prior,\delta)$, where we establish the lower bound; the upper bound follows from \Cref{prop:trivial-bounds-bay}.
Since the probability of observing ``A'' is zero, the minimum probability of error has a simple expression. \Cref{prop:E-gamma}, proved later, implies the following expression for the minimum expected error:
\begin{align*}
P_e^{(n)} &= \sum_{x_1,\dots,x_n \in \{A,B\}} \min\left( \prior p^{\otimes n} (x_1,\dots,x_n), (1- \prior) q^{\otimes n}(x_1,\dots,x_n) \right) \\ 
&= \min\left( \prior, (1- \prior) (1 -\epsilon)^n \right),
\end{align*}
which is less than $\delta$ if and only if $(1- \prior) (1 -\epsilon)^n \leq \delta$, or equivalently, $n \geq \frac{\log((1 -\prior)/\delta)}{-\log(1 - \epsilon)}$.
Since $\prior \leq 1/2$ and $\epsilon\leq 1/2$, the aforementioned condition implies that  $n \gtrsim  \frac{\log(1/\delta)}{\epsilon}$.
\end{proof}

In terms of the optimal test for
$\cB_\bay(p,q,\prior,\delta)$, the classical
result of Neyman and Pearson implies that the
likelihood ratio test (with a particular choice of
threshold) is optimal.
\begin{definition}[Likelihood ratio test]
Given a threshold $T \in [0, \infty]$ and two distributions $p$ and $q$,
the likelihood ratio test takes in $n$ samples $x_1,\dots,x_n$ and outputs $p$ if the likelihood ratio, $\prod_{i=1}^n \frac{p(x_i)}{q(x_i)}$,
is at least $T$ and $q$ otherwise.
\end{definition}

\begin{fact}[Optimal test \cite{NeyPea33}]
\label{fact:optimal-neyman-pearson}
For any $n$, the likelihood ratio test with the threshold
$\frac{1-\prior}{\prior}$ minimizes the probability of error  
in \Cref{eq:bayesian-error-prob-intro}.
\end{fact}

\subsubsection{Prior-free hypothesis testing}
We now turn our attention to the prior-free setting, recalled below.

\DefPriorFreeSimpHypTesting*
If both $\alpha$ and $\beta$ are larger than $1/2$, the problem is
trivial.
The case $\alpha = \beta$ is roughly equivalent to the uniform
prior setting with error probability $\alpha$ (which is characterized
by the Hellinger divergence).

\begin{proposition}[Standard Bounds on $\nstar_\pfht$]
We have the following bounds:
\label{prop:trivial-bounds-fre}
\begin{enumerate}
\item (Vacuous Regime) $\nstar_\pfht(p,q,\alpha,\beta) =0$ if $\alpha+ \beta \geq 1$.
\item   \label{cl:trivial-bounds-2b}
For all $\alpha, \beta,p$, and $q$ such that $\max(\alpha,\beta)< 1/4$ and $\hel^2(p,q) \leq \helconst$, we have
\begin{align}
\label{eq:trivial-bounds-2b}
\frac{\log(1/\max(\alpha,\beta))}{h^2(p,q)}  \lesssim \nstar_\pfht(p,q,\alpha,\beta) \lesssim \frac{\log(1/\min(\alpha,\beta))}{h^2(p,q)}\,.
\end{align}
\end{enumerate}
\end{proposition}
\begin{proof}
The test that chooses the distribution $p$ with probability $\overline
\alpha$ and $q$ with probability $\alpha$ has type-1 error
$\alpha$ and type-2 error $1 -\alpha$, which is at most $\beta$
if $\alpha + \beta \geq 1$.
This test does not use any samples.
The second claim in \Cref{eq:trivial-bounds-2b} follows by the
folklore result \cite{Ziv02,PolWu23} that $\nstar_\pfht(p,q,\delta, \delta) \asymp
\log(1/\delta)/h^2(p,q)$ for any $\delta \leq 1/4$ and $\hel^2(p,q) \leq \helconst$, combined with the following monotonocity property
\begin{align*}
 \nstar_\pfht(p,q,
\max\left( \alpha, \beta \right), \max\left( \alpha, \beta
\right))\leq \nstar_\pfht(p,q,\alpha,\beta) \leq
\nstar_\pfht(p,q, \min\left( \alpha, \beta \right), \min\left(
\alpha, \beta \right))\,. 
 \end{align*}
\end{proof}

\Cref{prop:trivial-bounds-fre} implies that when $\max(\alpha,\beta)$ 
is polynomially related to $\min(\alpha,\beta)$, the sample complexity 
is characterized by $\log(1/\alpha)/\hel^2(p,q)$.
Thus, the ``interesting'' regime of the parameters is when $\alpha$
and $\beta$ diverge.
Again, the sample complexity $\nstar_\pfht(p,q,\alpha,\beta)$ is not
symmetric in the first two arguments.

\subsubsection{Relation between these two problems} %
It is easy to see that the sample complexities of the Bayesian and prior-free hypothesis testing problems are closely related to each other.
\begin{claim}[Relation between Bayesian and prior-free sample complexities]
\label{claim:reln-bayesian-prior-free}
For any $\alpha, \beta \in (0,1)$, we have
\begin{align}
\label{eq:reln-bayesian-prior-free-1}
\nstar_\bay\left( p,q, \frac{\beta}{\alpha+\beta}, \frac{2\alpha \beta}{\alpha+\beta} \right) \le  \nstar_\pfht(p,q,\alpha,\beta) \le \nstar_\bay \left(p,q, \frac{\beta}{\alpha + \beta}, \frac{\alpha \beta}{\alpha + \beta}  \right),
\end{align}
and for all $\prior, \delta \in (0,1)$, we have
\begin{align}
\label{eq:reln-bayesian-prior-free-2}
\nstar_\pfht\left(p,q, \frac{\delta}{\prior}, \frac{\delta}{1 -\prior}\right) \le \nstar_\bay\left( p,q, \prior, \delta \right) \le \nstar_\pfht\left(p,q, \frac{\delta}{2 \prior}, \frac{\delta}{2(1 - \prior)}  \right).
\end{align}

\end{claim}
\begin{proof}
We begin with \Cref{eq:reln-bayesian-prior-free-1}.
\paragraph{Proof of \Cref{eq:reln-bayesian-prior-free-1}:}
The first inequality in \Cref{eq:reln-bayesian-prior-free-1} can
be checked by noting that the average error (with $p$ chosen with
probability $\beta/(\alpha + \beta)$) of the test achieving the
guarantee for $ \nstar_\pfht(p,q,\alpha,\beta)$ is at most $2
\alpha \beta/(\alpha+\beta)$, thus solving $\cB_\bay\left( p,q,
\frac{\beta}{\alpha+ \beta}, \frac{2 \alpha \beta }{\alpha +
\beta} \right)$.

For the second inequality, consider a test that solves
$\nstar_\bay \left(p,q, \frac{\beta}{\alpha + \beta},
\frac{\alpha \beta}{\alpha + \beta} \right)$.
Observe that the type-I error cannot be larger than $ \alpha$ and
the type-II error cannot be larger than $\beta$, or else the
weighted average error would exceed $ \alpha \beta/(\alpha +
\beta)$.

\paragraph{Proof of \Cref{eq:reln-bayesian-prior-free-2}:}

For the first inequality in \Cref{eq:reln-bayesian-prior-free-2},
consider a test $\phi$ that solves $\cB_\bay\left( p,q, \prior,
\delta \right)$.
Then the type-I error of the test $\phi$ cannot be more than
$\delta/\prior$; otherwise, the average error would be larger
than $\delta$.
Similarly, the type-II error cannot be larger than $\delta/(1 -
\prior)$.

For the second inequality, consider a test that solves
$\nstar_\pfht\left(p,q, \frac{\delta}{2 \prior},
\frac{\delta}{2(1 - \prior)} \right)$.
Then the average error, where $p$ is chosen with probability
$\prior$, is at most $\prior (\delta/ (2 \prior)) + (1 - \prior)
\left( \delta/(2 (1 - \prior)) \right) = \delta$.

\end{proof}

Finally, the following proposition shows that the sample complexity 
is stable under mild changes in the problem parameters,
implying that both series of inequalities in \Cref{claim:reln-bayesian-prior-free} can be reversed, up to constants:
\begin{proposition}[Mild changes in prior and error probability]
\label{prop:mild-reduction}
Let $ \prior_1 \in(0,1/2]$, $\prior_2 \in (0,1/2]$, $\delta_1 \in
(0, \prior_1/4)$, and $\delta_2 \in (0,\prior_2/4)$.
Then
\begin{align}
\label{eq:mild-change-1}
\frac{\nstar_\bay(p,q,\prior_1, \delta_1)}{\nstar_\bay(p,q,\prior_2, \delta_2)}   \lesssim \max\left(1, \frac{\log\left( \frac{\prior_1}{\delta_1} \right)}{\log\left( \frac{\prior_2}{\delta_2} \right)}, \frac{\log\left( \frac{1}{\delta_1} \right)}{\log\left( \frac{1}{\delta_2} \right)}  \right).
\end{align}
In particular, if $\delta_1 = \poly(\delta_2)$ and $\delta_1/\prior_1 = \poly(\delta_2/\prior_2)$, the sample
complexities are related by multiplicative constants. Similarly, let $\alpha_1, \beta_1, \alpha_2, \beta_2 \in
(0,1/4]$.
Then
\begin{align}
\label{eq:mild-change-2}
\frac{\nstar_\pfht(p,q,\alpha_1, \beta_1)}{\nstar_\pfht(p,q,\alpha_2, \beta_2)}   \lesssim \max\left(1, \frac{\log\left( \frac{1}{\alpha_1} \right)}{\log\left( \frac{1}{\alpha_2} \right)}, \frac{\log\left( \frac{1}{\beta_1} \right)}{\log\left( \frac{1}{\beta_2} \right)}  \right)\,.
\end{align}
{In particular, if $\alpha_1 = \poly(\alpha_2)$ and
$\beta_1 = \poly(\beta_2)$, the sample complexities are
related by multiplicative constants.
}

\end{proposition}
We state the following corollary:
\begin{corollary}
\label{cor:equivalence-bayesian-pf}
Let $\prior \in (0,1/2]$ and $\delta\leq \prior/4$.
Then
$\nstar_\bay(p,q,\prior,\delta) \asymp \nstar_\pfht \left( p,q, \frac{\delta}{\prior},  \delta \right)$.
Similarly, let $\alpha, \beta \in (0,1/8)$ be such that $\beta \leq \alpha$.
Then
$\nstar_\pfht\left( p,q,\alpha, \beta \right) \asymp \nstar_\bay\left( p,q, \frac{\beta}{\alpha + \beta}, \frac{\alpha \beta}{\alpha + \beta} \right)$.
\end{corollary}
\begin{proof}[Proof of \Cref{cor:equivalence-bayesian-pf}]
By \Cref{claim:reln-bayesian-prior-free}, the Bayesian sample
complexity is sandwiched between the two following prior-free
sample complexities, so we obtain a bound on the ratio using
\Cref{prop:mild-reduction}:
\begin{align*}
\frac{\nstar_\pfht \left(p,q, \frac{\delta}{2 \prior}, \frac{\delta}{2 (1- \prior)}  \right)}{\nstar_\pfht\left( p,q, \frac{\delta}{\prior}, \frac{\delta}{1 -\prior} \right)}
&\lesssim \max\left( 1, \frac{\log(2\prior/\delta)}{\log(\prior/\delta)}, \frac{\log(2(1 -\prior)/\delta)}{\log((1- \prior)/\delta)} \right)\,.
\end{align*}
Both ratios above are at most a constant, since
$\log(1/x)/\log(1/2x) \lesssim 1$ when $x \in (0,1/4)$, and
$\delta/\prior \in (0,1/4)$ and $\delta/(1- \prior) \leq \delta/
\prior \leq 1/4$, since $\prior \leq 1/2$.

We now consider the second claim: By
\Cref{claim:reln-bayesian-prior-free}, it suffices to upper-bound
the ratio of the following Bayesian sample complexities, which
can be controlled using \Cref{prop:mild-reduction} (since the
priors satisfy $\frac{\beta}{\alpha + \beta} \leq 1/2$ and the
ratio of the error probability and prior is at most $1/4$):
\begin{align*}
\frac{\nstar_\bay \left(p,q, \frac{\beta}{\alpha + \beta}, \frac{\alpha \beta}{\alpha + \beta}  \right)}{\nstar_\bay\left( p,q, \frac{\beta}{\alpha+\beta}, \frac{2\alpha \beta}{\alpha+\beta} \right)}
&\lesssim \max\left( 1, \frac{\log(1/\alpha)}{\log(1/2\alpha)}, \frac{\log((\alpha+\beta)/\alpha \beta)}{\log((\alpha+\beta)/2\alpha \beta)} \right)\,.
\end{align*}
Both ratios above are at most a constant, since
$\alpha \in (0,1/4)$ and $ \alpha \beta/ (\alpha + \beta) \leq
\min(\alpha, \beta) \leq 1/4$.
\end{proof}

The following folklore result about boosting the probability of success 
using repetition, proved in \Cref{app:proof_of_cref_fact_boost} for completeness, will be crucial:
\begin{fact}[Repetition to boost success probability]
\label{fact:boost}
Consider a function $\phi:\cX^n \to \{0,1\}$ such that
$\P(\phi(X_1,\dots,X_n) \neq 0) \leq \tau$ for $\tau \leq 1/4$ when $X_1,\dots,X_n$
are sampled i.i.d.\ under $P$.\footnote{The constant $1/4$ can be improved to any fixed constant less
than $1/2$, at the cost of a bigger constant $2^{5}$ below.}
Consider the modified function $\phi':\cX^{nT} \to \{0,1\}$ that is defined
to be the majority of $\phi(X_1,\dots,X_n)$, $\dots$,$\phi(X_{(T-1)n
+ 1}\dots, X_{nT})$.
Then $\P(\phi'(X_1,\dots,X_{nT}) \neq 0)$ when
$X_1,\dots,X_{nT}$ are sampled i.i.d.\ from $P$ is at most
$\exp\left( - 2^{-5}
T \log(1/ \tau) \right)= \tau^{T/32}$.

Alternatively, for any $\tau' \leq \tau$ and any $T \geq
\frac{2^5 \log(1/\tau')}{\log(1/\tau)} $, the probability of
error (when the true distribution is $P$) of the boosted procedure is at most $\tau'$.
\end{fact}

\begin{proof}[Proof of \Cref{prop:mild-reduction}]
We begin with the proof of \Cref{eq:mild-change-2}.
Consider a test $\phi$ that solves $\cB_\pfht(p,q,\alpha_2,
\beta_2)$ with $n:= \nstar_\pfht(p,q,\alpha_2,\beta_2)$ samples.
We want to boost the test $k$ times in the sense of
\Cref{fact:boost}, using $nk$ samples in total, so that the
resulting type-I error is at most $\alpha_1$ and the type-II error
is at most $\beta_1$.
\Cref{fact:boost} then implies the result.

Similar arguments work for \Cref{eq:mild-change-1}.
Using \Cref{claim:reln-bayesian-prior-free}, it suffices to solve
$\cB_\pfht(p,q, (\delta_1/2 \prior_1), \delta_1/2)$ to solve
$\cB_\bay(p,q,\prior_1,\delta_1)$.
Consider a test $\phi$ that solves $\cB_\bay(p,q,\prior_2,
\delta_2)$ with $n:= \nstar_\bay(p,q,\prior_2,\delta_2)$ samples.
By \Cref{claim:reln-bayesian-prior-free} again, the type-I error of $\phi$ on $n$ samples is at most
$\delta_2/\prior_2$ and the type-II error is at most
$\delta_2/(1 -\prior_2) \leq 2 \delta_2$, both of which are
less than $1/4$.
We boost the test $k$ times so that the resulting type-I
error is at most $(\delta_1/2 \prior_1)$ and the type-II error is
at most $\delta_1/2$.
\Cref{fact:boost} then implies the result.
\end{proof}

\subsection{Probability divergences}
\label{sec:probability_divergences}

Our results crucially rely on $f$-divergences, which we formally define below:
\begin{definition}[$f$-divergence]
For a convex function $f:\R_+ \to \R$ with $f(1)= 0$, we use
$I_f(p,q)$ to denote the $f$-divergence between
$p$ and $q$, defined as $I_f(p,q):= \sum_{i}q_i f\left( p_i/q_i \right)$.\footnote{We use the following conventions~\cite{Sason18}: $f(0) = \lim_{t \to 0^+}f(t)$, $0 f(0/0) = 0$, and for $a>0$, $0f(a/0) = a \lim_{u \to \infty} f(u)/u$.}
\end{definition}

We will repeatedly use the following $f$-divergences and their properties:

\begin{definition}[$E_\gamma$-divergence]
\label{def:E-gamma}
For $\gamma \geq 1$, the $E_\gamma$-divergence, also known as
the hockey-stick divergence, is defined by $E_\gamma(p,q) :=
\sum_{i}(p_i - \gamma q_i)_+$.
\end{definition}
Observe that the $E_\gamma$-divergence is an $f$-divergence with
$f(x) = \max(0,x- \gamma)$.

\begin{proposition}[{Variational characterization of $E_\gamma$, see, e.g., \cite[Section VII]{SasVer16}}]
\label{prop:E_gamma_variational}
Let $p$ and $q$ be two distributions over $\cX$ with
$\sigma$-algebra $\cF$.
Then $E_\gamma(p,q) := \sup_{A \in \cF} \left( \P_p(A) - 
\gamma \P_q(A) \right)$.

\end{proposition}
As an application, we obtain the following exact characterization
of the Bayesian error probability:
\begin{proposition}[Exact characterization of error probability]
\label{prop:E-gamma}
Let $p$ and $q$ be two distributions over $\cX$, and let $\alpha
\in (0,0.5]$.
Then \begin{align*}
P_e^{(1)} =\min_{\phi:\cX \to \{p,q\}} \Big(\prior \P_{x \sim p} \left\{   \phi(x) \neq p\right\} + (1 - \prior) \P_{x \sim q} \left\{   \phi(x) \neq q\right\}\Big) 
&= \prior - \prior E_{\frac{1- \prior}{\prior}} \left( p,q \right)\,.
\\
&= \sum_i\min\left( \prior p_i, (1- \prior)q_i \right).
\end{align*}
In particular,
for any $n$, the minimum possible error probability for
the Bayesian simple binary hypothesis testing problem with prior
$\prior$, $P_e^{(n)}$, is $\prior\left(1 - E_{\frac{1-
\prior}{\prior}}\left(p^{\otimes n}, q^{\otimes
n}\right)\right)$.
\end{proposition}
Thus, the sample complexity for $\cB_\bay(p,q,\prior,\delta)$ is the
minimum $n$ such that $\prior\left(1 - E_{\frac{1-
\prior}{\prior}}\left(p^{\otimes n}, q^{\otimes n}\right)\right)$ is
at most $\delta$.
In particular, for $\prior=1/2$, the $E_1$-divergence becomes the
total variation distance, and the tight characterization of the
sample complexity for $\nstar_\bay(p,q,1/2,\delta)$ (for small
$\delta \leq 1/4$) in terms of the Hellinger divergence is obtained
by upper- and lower-bounding the total variation distance in terms of
the Hellinger divergence, coupled with the tensorization property of
the Hellinger divergence~\cite{Ziv02}.
Such an approach cannot be replicated for the $E_\gamma$-divergence.
Unlike most divergences such as the total variation distance, Hellinger
divergence, $\chi^2$-divergence, or Kullback--Leibler divergence, the
$E_\gamma$-divergence between $p$ and $q$ may be 0 even when $p \neq
q$.
This makes it impossible to upper-bound any such divergence by the
$E_\gamma$-divergence, necessitating a new approach.

Our results will use the following generalization of the Hellinger 
divergence (which reduces to the usual definition when $\lambda=1/2$):
\begin{definition}[$\cH_\lambda$-divergence]
\label{def:alpha-hel}
For $\lambda \in [0,1] $,
define the $\cH_\lambda(p,q)$-divergence to be
\begin{align*}
\cH_\lambda(p,q) :=   1 -  \sum_{i} p_i ^ \lambda q_i ^{1 - \lambda} = \sum_i q_i\left( 1 - \left( \frac{p_i}{q_i} \right)^ \lambda \right)\,.
\end{align*}
\end{definition}
The $\cH_\lambda$-divergence is an $f$-divergence with $f(x) = 1
- x^{\lambda}$, and is jointly convex in $(p,q)$ for $\lambda \in [0,1]$.
The $\cH_\lambda$ divergence generalizes the Hellinger divergence,
which is obtained (up to normalization) by setting $\lambda = 1/2$.
Moreover, the $\cH_\lambda$-divergence is closely related to the
R\'enyi divergence of order $\lambda$, via
$D_\lambda(p,q) = \frac{1}{\lambda} \log\left( 1 - \cH_\lambda(p,q)
\right)$.

\begin{proposition}[Tensorization of $\cH_\lambda$]
\label{prop:hel-tensor}
We have
$1 - \cH_\lambda(p^{\otimes n}, q^{\otimes n}) = \left(1
- \cH_\lambda(p,q) \right)^n$.

\end{proposition}
\begin{proof}
By expanding the definition of $\cH_\lambda(p^{\otimes n}, q^{\otimes n})$, we obtain
\begin{align*}
1 - \cH_\lambda(p^{\otimes n}, q^{\otimes n}) &= \sum_{x_1,\dots,x_n} \left( p_{x_1}\dots p_{x_n} \right)^{\lambda} \left( q_{x_1}\dots q_{x_n} \right)^{1-\lambda} \\
&= \sum_{x_1,\dots,x_n} \left( p_{x_1}^\lambda q_{x_1}^{1- \lambda}\right)\dots \left( p_{x_n}^\lambda  q_{x_n}^{1-\lambda}\right) \\
&=  \left(\sum_{x_1} p_{x_1}^\lambda q_{x_1}^{1- \lambda}\right)\dots \left( \sum_{x_n}p_{x_n}^\lambda  q_{x_n}^{1-\lambda}\right) \\
&= \left(1 - \cH_\lambda(p,q)  \right)^n.
\end{align*}
\end{proof}

We also define the following generalization of the Jensen--Shannon divergence:
\begin{definition}[Skewed $\js{\alpha}$-divergence]
\label{def:skew-symm}
For $\prior \in [0,1]$, we define the skewed Jensen--Shannon
divergence by
\begin{align}
\js{\prior}(p,q) := \prior\,\kl(p, \prior p + (1- \prior) q) + (1-\prior)\, \kl\,(q, \prior p + (1- \prior) q) \,.
\end{align}
\end{definition}
This divergence was called the asymmetric $\prior$-skew Jensen--Shannon divergence in \cite{Nie20}. Similar quantities have appeared previously in the literature, as well
\cite{Lin91,Ziv02,Nie11}
termed the $\prior$-Jensen--Shannon
divergence, which was defined to be $ K_\prior(p,q) +
K_\prior(q,p)$, where $K_\prior(p,q):= \kl(p, \prior p + (1-
\prior) q)$.
In comparison, $\js{\prior}(p,q) = \prior K_\prior(p,q) + (1- \prior)
K_{(1- \prior)}(q,p)$, thus is asymmetric.
Combining the facts that $K_\prior(\cdot,\cdot)$ is an
$f$-divergence~\cite{Lin91,Nie11}; the sum of two $f$-divergences is
an $f$-divergence; and the reversal of an $f$-divergence is another
$f$-divergence; we conclude that $\js{\prior}(p,q)$ is also an
$f$-divergence.

\begin{proposition}[Information-theoretic interpretation of skewed $\js{\prior}$-divergence]
\label{prop:skew-symm}
Let $\Theta$ be a $\Ber( \overline{\prior})$ random variable.
Let $X \sim p$ if $\Theta = 0$ and $X \sim q$ if $\Theta =
1$.
Then $I(\Theta; X) = \js{\prior}(p,q)$.

\end{proposition}
\begin{proof}
Let $P$ denote the joint distribution over $\theta$ and $X$.
Let $P_X$ denote the distribution of $X$ and let $P_\Theta$ denote
the distribution of $\Theta$. 
Then
\begin{align*}
I(X;\Theta) 
& = \kl(P, P_X\times P_\Theta) 
= \E_\Theta\left[\kl(P_{X|\Theta}, P_X)\right]\\
&= \prior \kl(p, P_X) + (1 -\prior) \kl(q, P_X),
\end{align*}
where the final result follows by noting that $P_X = \prior p + (1 -\prior) q$.
\end{proof}

For a random variable $Y \sim p_Y$, let
$H_{\text{ent}}(Y)$ 
 (also
written as $H(p_Y)$) be the Shannon entropy of
$Y$ defined as
$H_{\text{ent}}(Y) = -\sum_{y} p_Y(y) \log
p_Y(y)).
$
For $x \in [0,1]$, define
$\hbin(x)$ to be the Shannon
entropy of a $\Ber(x)$ random
variable.

\begin{fact}[Fano's inequality~\cite{CovTho06,PolWu23}]
\label{fact:fano}
Let $V \to X \to \widehat{V}$ be any Markov chain, where $V$
is a binary random variable..
Define $p_e:= \P(V \neq \widehat{V})$.
Then
\begin{align*}
\hbin(p_e) \geq \hent(X|V).
\end{align*}
\end{fact}

\begin{lemma}[Method of joint range~\cite{HarVaj11}] 
\label{lem:joint-range}
Let $D_{f_1}(\cdot,\cdot)$ and $D_{f_2}(\cdot,\cdot)$ be two $f$-divergences.
Suppose that for a constant $c >0$, we have
\begin{equation*}
D_{f_1}(\Ber(x),\Ber(y)) \leq
cD_{f_2}(\Ber(x),\Ber(y))
\end{equation*}
for all $x, y \in [0,1]$.
Then $D_{f_1}(p,q) \leq cD_{f_2}(p,q)$ for all probability
distributions $p$ and $q$.
\end{lemma}
\begin{proof}
The result from \cite{HarVaj11} states that $\cR$ is the convex hull of $\cR_2$, where
\begin{align*}
\cR &:= \{(u,v): u = D_{f_1}(p, q), v = D_{f_2}(p, q) \text{ for any $p$ and $q$ over any $\cX$} \} \text { and } \\
\cR_2 &:= \{(u,v): u = D_{f_1}(p, q), v = D_{f_2}(p, q) \text{ for any $p$ and $q$ over $\cX = \{0,1\}$} \}.
\end{align*}
Since the $f$-divergence inequality is given to hold for Bernoulli distributions, we know that $\cR_2 \subseteq \{(u, v): u \le cv\}$. Since the latter is a convex subset of $\real_+^2$, we conclude that the convex hull of $\cR_2$, which is the same as $\cR$, also lies within the same set. This concludes the proof.
\end{proof}

\section{Proof Sketch of \Cref{thm:main-result-intro-bay}}

In this section, we present a high-level sketch of the
proof of \Cref{thm:main-result-intro-bay}. 

\subsection{Linear error probability regime}
The formal statements of all results in this section and their complete proofs are in \Cref{sec:bayes-testing}.

We begin by focusing on the regime where $\delta$ is a small
fraction of $\prior$. Suppose $\delta \in (\prior/100,\prior/4)$.
We begin by presenting sample complexity lower bounds and upper bounds, respectively, using classical techniques. Finally, we present our main  technical result that establishes their equivalence.

A sample complexity lower bound is obtained via Fano's 
inequality:
\begin{proposition}[Lower bound, formally stated in \Cref{thm:lower-bound-bayes}]
\label{thm:lower-bound-bayes-intro}
$\nstar_\bay(p,q,\prior, \prior/4) \gtrsim \frac{\prior \log(1/\prior)}{I(\Theta;X_1)}$.
\end{proposition}
\begin{proof}[Proof Sketch] For $x\in(0,1)$, we use $\hbin(x)$ to 
denote the binary entropy of $x$ and $\hent(\cdot)$ to denote the entropy (and conditional entropy) of a random variable or a distribution.
Observe that $\Theta \to (X_1,\dots,X_n) \to \widehat{\Theta}$ is
a Markov chain.
Let $P_e^{(n)}$ be the probability of error, i.e., $\P(\Theta
\neq \widehat{\Theta})$.
By Fano's inequality (\Cref{fact:fano}), we have
\begin{align}
\label{eq:fano-lower-bd-intro}
\hbin(P_e^{(n)}) & \geq \hent(\Theta|X_1,\dots,X_n) = \hent(\Theta) - I(\Theta;
X_1{,}\dots{,}X_n) \notag \\
& = \hbin(\prior) - I(\Theta;
X_1{,}\dots{,}X_n).
\end{align}
The mutual information term is upper-bounded using the conditionally i.i.d.\ distribution of the $X_i$'s given $\Theta$, as
$I(\Theta; X_1,\dots,X_n) \le \sum_{i=1}^n I(\theta; X_i) = nI(\Theta; X_1)$.
Combining this inequality with
\Cref{eq:fano-lower-bd-intro}, we obtain
$n \geq
\frac{\hbin(\prior) - \hbin(P_e^{(n)})}{I(\Theta;X_1)}$.
As $P_e^{n} \leq \prior/4$, the desired inequality follows by noting that
$\hbin(\prior) - \hbin(\prior/4) \asymp \hbin(\prior)$.
\end{proof}
We note that a similar strategy for establishing a lower bound for the uniform prior was used in \cite[Theorem 4.27]{Ziv02}.
For a uniform prior and small constant error probability, this
gives an alternate proof of the $\frac{1}{\hel^2(p,q)}$ sample complexity lower bound, because $I(\Theta;X_1)
\asymp \hel^2(p,q)$ when $\Theta$ is uniform over $p$ and $q$~\cite{Ziv02}.
Turning our attention to an upper bound, we use the classical idea
of tensorization for the Hellinger family of $f$-divergences.
This tensorization leads directly to the following classical upper bound:

\begin{proposition}[Upper bound, formally stated in \Cref{thm:upper-bound-bayes}]
\label{thm:upper-bound-bayes-intro}
Set $\lambda = \frac{0.5 \log 2}{\log(1/\prior)}$, which lies in
$(0,0.5]$, and define $\bl := 1 - \lambda$.
Then $\nstar_\bay(p,q,\prior, \prior/4) \lesssim \left\lceil
\frac{2}{\cH_{\bl}(p,q)}\right\rceil$.
\end{proposition}

\begin{proof}[Proof Sketch]
Let $P_e^{(n)}$ be the minimum Bayesian error between $p$ and $q$, given $n$ samples.
By \Cref{prop:E-gamma}, we have $P_e^{(n)} = \sum_{x_1,\dots,x_n}
\min\left( \prior p^{\otimes n} (x_1,\dots,x_n), (1- \prior)
q^{\otimes n}(x_1,\dots,x_n) \right)$.
Since $\min(x,y) \leq x^{\lambda}y^{1 - \lambda}$ for any
nonnegative $x$, $y$ and $ \lambda \in (0,1)$, we obtain
\begin{align*}
P_e^{(n)} 
&\leq  \sum_{x_1,\dots,x_n}  \left(\prior p^{\otimes n}(x_1,\dots,x_n)\right)^{\bl} \left((1- \prior) q(x_1,\dots,x_n) \right)^{\lambda} \\
&= \prior^{\bl}(1 -\prior)^{\lambda} \left( 1 - \cH_{\bl}(p^{\otimes n}, q^{\otimes n} ) \right) \leq \prior^{\bl } \left( 1 -  \cH_{\bl}(p, q) \right)^n \tag{\Cref{prop:hel-tensor}}\\
&\leq \prior^{\bl} \left( e^{- \cH_{\bl}(p, q)} \right)^n\,. \tag{$1-x \leq e^{-x}$ for all $x$}
\end{align*}
By definition of $\lambda$, we have $\prior^{\bl} \lesssim
\prior$.
Thus, if $e^{- \cH_{\bl}(p, q) n}$ is sufficiently small, the whole error probability will be less than
$\prior/4$, yielding the desired upper bound on the sample
complexity.
\end{proof}

We remark that similar arguments have appeared classically in \cite{HelRav70}.

\paragraph{Main technical result:} We now arrive at the main technical contribution of our work.
Following the standard tools in our toolkit, we have obtained lower
and upper bounds that depend on two different quantities:
$I(\Theta;X_1)$ and $\cH_{\bl}(p,q)$, for $\lambda := \frac{0.5 \log
2}{\log(1/\prior)}$.
Our main contribution is to show that both these quantities
are equivalent up to constant factors for all $p$, $q$, and $\prior$.\footnote{The inequality
$\frac{\prior \log(1/\prior)}{I(\Theta;X_1)} \lesssim
\frac{1}{\cH_{\bl}(p,q)}$ follows directly from these results.}
Our proof technique relies on noting that $I(\Theta;X_1)$ is equal to the $f$-divergence $\js{\prior}(p,q)$  defined in \Cref{def:skew-symm}; see \Cref{prop:skew-symm}.

Thus, our task reduces to showing a tight relationship between the
$f$-divergences $\js{\prior}(p,q)$ and $\cH_{\bl}(p,q)$.
Framing this inequality in terms of $f$-divergences permits us to use
the powerful joint range method of \cite{HarVaj11} from \Cref{lem:joint-range}.
Thus, it remains to show that $\frac{\prior
\log(1/\prior)}{\js{\prior}(p,q)} \gtrsim \frac{1}{\cH_{\bl}(p,q)}$
for Bernoulli distributions $p$ and $q$.

\begin{lemma}[Relationship between $\js{}$ and $\cH$, formally stated in \Cref{conj:the-two-div}]
\label{conj:the-two-div-informal}
Let $p$ and $q$ be two Bernoulli distributions.
Let $\prior \in (0,1/2]$ and $\lambda = \frac{0.5 \log
2}{\log(1/\prior)}$.
Then $\frac{1}{\cH_{\bl}(p,q)} \lesssim \frac{ \prior
\log(1/\prior)}{ \js{\prior}(p,q) }$.
\end{lemma}
The proof of the above lemma is rather technical
and follows by careful case-by-case analyses.

\subsection{Sublinear error probability regime}
\label{sub:sublinear_failure_probability_regime}
The formal statements of all results in this section and their complete proofs are in \Cref{sec:reduction_between_problems}.

In this section, we sketch the argument for the proof of
\Cref{thm:main-result-intro-bay} in the regime $\delta \in
(\prior^2, \prior/100)$. Instead of characterizing the sample complexity in terms of an
$f$-divergence, we adopt a reduction approach and establish a
tight relationship between $\nstar(p,q,\prior,\delta)$ and
$\nstar(p,q,\prior',\delta')$.

Our starting point is the standard boosting lemma, which
states that if a decision procedure fails with
probability $\tau \leq 1/4$ with $n$ samples, then
taking the majority vote of $T $ (independent)
runs of the same decision procedure fails with
probability $\tau' \ll \tau$ if $T \asymp
\log(1/\tau')/\log(1/\tau)$;
see~\Cref{fact:boost}.
As a consequence, we immediately obtain the following result:

\begin{lemma}[Success amplification, formally stated in \Cref{lem:upper-bound-fre-tech}]
\label{lem:upper-bound-fre-tech-intro}
For all positive integers $T$ such that
$\max\left((\delta/\prior)^{1/T},\prior^{1/T}\right) \leq 1/8$, we have $
\nstar_\bay\left( p,q, \prior, \delta \right) \lesssim T \cdot
\nstar_\bay\left( p,q, \prior^{1/T} , \delta^{1/T} \right)$.
\end{lemma}
However, the above boosting procedure might be lossy; recall that
the optimal procedure remains the likelihood ratio test.
The following result establishes that the likelihood ratio test also
suffers the same sample complexity (up to constant factors):

\begin{lemma}[Error amplification, formally stated in \Cref{lem:bay-to-pfht-low-bd-amp}]
\label{lem:lower-bound-fre-tech-intro}
Under the same setting as in \Cref{lem:upper-bound-fre-tech-intro}, we have $\nstar_\bay\left( p,q, \prior, \delta \right) \gtrsim T
\nstar_\bay\left( p,q, \prior^{1/T} , \delta^{1/T} \right)$.
\end{lemma}
\begin{proof}[Sketch]
Consider the likelihood ratio test with $n := \nstar(p,q,\prior,\delta)$
samples.
Divide the $n$ samples into $T$ buckets, each of size
$n/T$; for simplicity of presentation, we assume that $n$ is divisible by $T$ for this proof sketch.
Let the $i^{\text{th}}$ bucket be $B_i$ and denote the likelihood ratio on
$B_i$ by $L_i := \log\left(\prod_{x \in B_i}^n
\frac{p_{x}}{q_{x}}\right)$.
The optimal likelihood ratio test $\phi^*$ takes $n$ samples
$x_1,\dots,x_n$ and returns $p$ if the joint likelihood
$L:=\sum_{i \in [T]} L_i$ is at least $\log((1 - \prior)/
\prior)$, and $q$ otherwise.

Consider the event $\cE_i$ over the points in the bucket $B_i$
(which has $m$ points), defined to be the domain where $L_i \leq
\log((1- \prior)/\prior)/T$.
Define $\cE:= \bigcap_{i=1}^{T} \cE_i$ and $\cE' =
\bigcap_{i=1}^{T} \cE_i^\complement$.
On the event $\cE$, we have $L < \log((1 - \prior)/\prior)$,
and on event $\cE'$, we have $L\geq \log((1 -
\prior)/\prior)$.

By the definition of the test, the (average) probability of error
of $\phi^*$ is
\begin{align*}
\P(\text{error}) &=  \prior \P(L < \log(\prior/\beta)|p) + (1 -\prior) \P(L > \log(\prior/\beta)|q) \\
&\geq \prior \P(\cE|\Theta = p) + (1 - \prior) \P(\cE'|\Theta = q) \\
&= \prior \left( \P(\cE_1|\Theta = p) \right)^T + (1 - \prior) \left( \P(\cE_1^\complement|\Theta = q) \right)^T \,,
\numberthis \label{eq:intro-low-boost-proof}
\end{align*}
where the last step uses the fact that the events are conditionally i.i.d. Consider the test $\phi$ on $m$ samples that outputs $p$ on
$\cE_1^\complement$ and $q$ on $\cE_1$.
Let $e_1$ and $e_2$ be its marginal failure probabilities,
defined as $\P(\cE_1|\Theta = p)$ and
$\P(\cE_1^\complement|\Theta = q)$, respectively.
By \Cref{eq:intro-low-boost-proof}, we have $\prior e_1^T + (1-
\prior) e_2^T \leq \delta $.
Thus, $e_1 \leq (\delta/\prior)^{1/T}$ and $e_2 \leq (2
\delta)^{1/T}$.
Suppose we use $\phi$ to solve
$\cB(p,q,\prior',\delta')$ for $\prior' = \prior^{1/T}$ and
$\delta' = 3\delta^{1/T}$.
Its average error probability will be at most $\prior'
e_1 + (1 - \prior') e_2 \leq 3 \delta^{1/T}$.
By the definition of sample complexity, we have $m \geq
\nstar(p,q,\prior',\delta')$, so
$\nstar \geq T \nstar(p,q,\prior',\delta')$.\footnote{For simplicity of presentation, we ignore the distinction
between the failure probabilities $\delta^{1/T}$ and $3 \delta^{1/T}$.}
\end{proof}

\section{Reduction: Error amplification and success amplification} %
\label{sec:reduction_between_problems}

In this section, we first show by a series of
simple reductions that we can reduce the problem
of computing $\nstar_\bay(p,q,\prior,\delta)$ for
any $\delta \leq \prior/4$ to computing
$\nstar_\bay(p,q,\prior',\prior'/10)$ for some
$\prior'$ appropriately defined, formalizing the discussion from \Cref{sub:sublinear_failure_probability_regime}.
Later, our main result will be to give an expression for
$\nstar_\bay(p,q,\prior,\prior/10)$.

\begin{proposition}[Self-reduction]
\label{prop:reduction}
Let $p$ and $q$ be two arbitrary distributions such that
$\nstar_\pfht(p,q,1/4,1/4) \geq 2$.\footnote{The condition
$\nstar_\pfht(p,q,1/4,1/4) \geq 2$ is guaranteed by $\hel^2(p,q,)
\leq \helconst$.}
\begin{enumerate}
\item (Prior-free testing problem)
Let $\alpha \in (0,1/8)$ and $\beta \in (0,1/4)$.
Then for all integers $T$ such that $\max\left(
(2\alpha)^{1/T}, (2\beta)^{1/T} \right)\leq 1/4$, we have
$\nstar_\pfht\left( p,q,\alpha,\beta \right) \asymp T
\cdot \nstar_\pfht\left( p,q, \alpha^{\frac{1}{T}},
\beta^{\frac{1}{T}}\right)\,$.

In particular, if $ \beta\leq \alpha < 1/8$,
let $T = \left\lfloor \frac{\log(1/ 2 \alpha)}{\log
4}\right\rfloor$, so that $\alpha^{\frac{1}{T}} \in \left[\frac{1}{16}, \frac{1}{4}\right]$.
Defining $\beta':= \beta^{1/T} $,
we obtain
\begin{align*}
\nstar_\pfht\left( p,q,\alpha,\beta \right) \asymp \log\left( \frac{1}{\alpha} \right) \nstar_\pfht\left( p,q, 1/4, \beta' \right)\,.
\end{align*}

\item (Bayesian testing problem)
Let $\prior \in (0,1/2]$ and $\delta \in
(0,\prior/8)$.
For any integer $T$ such that $
(\delta/\prior)^{\frac{1}{T}} \leq 1/8$ and
$(\prior)^{\frac{1}{T}} \leq 1/2$, we have
$\nstar_\bay\left( p,q,\prior,\delta \right) \asymp T
\cdot \nstar_\bay\left( p,q, \prior^{\frac{1}{T}},
\delta^{\frac{1}{T}}\right)$.%

In particular, if $\delta \in (\prior^2,\prior/8)$,
we can set $T = \left\lfloor \frac{\log(\prior/\delta)}{ \log
8} \right\rfloor$ and $\prior':= \prior^{\frac{1}{T}}$ (so
that $\delta^{\frac{1}{T}} \in \left[ \frac{\prior'}{64} , \frac{\prior'}{8} \right]$ and $\prior' \in
(0,1/2]$) to obtain
\begin{align*}
\nstar_\bay\left( p,q,\prior,\delta \right) \asymp \log\left( \frac{\prior}{\delta} \right)  \nstar_\bay\left( p,q, \prior', \prior'/4\right)\,.
\end{align*}
\end{enumerate}
\end{proposition}
The proof of the proposition is contained in
\Cref{SecPropReduction}.%

\paragraph{Error amplification:} We begin by deriving a lower bound on the prior-free sample
complexity in terms of an easier problem: a prior-free testing problem
on the same distributions, but with larger error probabilities.
\begin{lemma}[Error amplification]
\label{lem:self-reduction-freq-hp}
For any positive integer $T$, we have
\begin{align}
\nstar_\pfht(p,q,\alpha,\beta) \geq T \cdot \left(\nstar_\pfht \left(p,q,(2\alpha)^{1/T},(2 \beta)^{1/T}\right) -1 \right).
\end{align}
\end{lemma}

\begin{proof}[Proof of \Cref{lem:self-reduction-freq-hp}]
Exactly characterizing the sample complexity in the prior-free 
setting is difficult, because it is not clear what the optimal 
threshold should be for the likelihood ratio test.
Thus, we first reduce to Bayesian setting, using
\Cref{eq:reln-bayesian-prior-free-1} in
\Cref{claim:reln-bayesian-prior-free}:
\begin{align}
\label{eq:lem:self-reduction-freq-hp-1}
\nstar_\pfht(p,q,\alpha,\beta) \geq \nstar_\bay\left(p, q, \frac{\beta}{\alpha+ \beta}, \frac{2\alpha \beta}{\alpha + \beta}\right)\,.
\end{align}
Using the exact characterization of the optimal test for
the Bayesian setting (\Cref{fact:optimal-neyman-pearson}), we
establish the following error amplification result:
\begin{lemma}[Error amplification]
\label{lem:bay-to-pfht-low-bd-amp}
For any positive integer $T$, we have
\begin{align*}
\nstar_\bay\left(p, q, \frac{\beta}{\alpha+ \beta}, \frac{2\alpha \beta}{\alpha + \beta}\right) \geq T \cdot \left(\nstar_\pfht(p,q, (2\alpha)^{1/T}, (2\beta)^{1/T}) -1 \right)\,.
\end{align*}
\end{lemma}
The above lemma, combined with
\Cref{eq:lem:self-reduction-freq-hp-1}, completes the proof of
\Cref{lem:self-reduction-freq-hp}.
\end{proof}

We now provide the proof of
\Cref{lem:bay-to-pfht-low-bd-amp}.

\begin{proof}[Proof of \Cref{lem:bay-to-pfht-low-bd-amp}]
Consider the problem $\cB_\bay\left(p,q, \frac{\beta}{\alpha+ \beta}, \frac{2 \alpha \beta}{\alpha + \beta}\right)$ with sample complexity $\nstar := \nstar_\bay\left(p,q, \frac{\beta}{\alpha+ \beta}, \frac{2 \alpha \beta}{\alpha + \beta}\right)$. Given any positive integer $T$, let $n$ be the smallest integer, at least $\nstar$, which is divisible by $T$.
Observe that $ T\lceil \nstar/T\rceil \geq n \geq \nstar$. By \Cref{fact:optimal-neyman-pearson}, an optimal test $\phi^*$
takes $n$ samples $x_1,\dots,x_n$ and returns $p$ if
$L:=\log(\prod_{i=1}^n p_{x_i}/q_{x_i})$ is at least
$\log(\alpha/\beta)$, and $q$ otherwise. Since $n\geq \nstar$, the average error probability must be at most $2 \alpha \beta/ (\alpha + \beta)$.

Let $\alpha' = (2\alpha)^{1/T}$ and $\beta' = (2\alpha)^{1/T}$.
We will now show that unless $n \geq T \cdot
\nstar_\pfht(p,q,\alpha',\beta')$, the (average) probability of
error of the test $\phi^*$ is strictly larger than $2\alpha
\beta/(\alpha + \beta)$, which gives the desired claim.

Divide the entire dataset into $T$ buckets $B_1,\dots,B_{T}$,
each of size $n/T$ (such that the partition is independent
of the data). Define $L_i$ to be the likelihood ratio on the bucket $B_i$,
i.e., $L_i = \log(\prod_{x \in B_i} \left( p_{x}/q_x \right))$.
We have $L = \sum_{i=1}^{T} L_i$.
Consider the event $\cE_i$ over the points in the bucket $B_i$
(which has at most $m := n/T$ points) defined to be the domain
where $L_i \leq \log(\alpha/\beta)/T$.
Define $\cE:= \bigcap_{i=1}^{T} \cE_i$ and $\cE' =
\bigcap_{i=1}^{T} \cE_i^\complement$.
On the event $\cE$, we have $L < \log(\alpha/\beta)$, and on
event $\cE'$, we have $L\geq \log(\alpha/\beta)$.

By the definition of the test $\phi^*$, the (average) probability of error
of $\phi^*$ is as follows:
\begin{align*}
\P(\text{error}) &= 
\frac{\beta}{\alpha + \beta} \P(L < \log(\alpha/\beta)|p) +
\frac{\alpha}{\alpha + \beta} \P(L > \log(\alpha/\beta)|q) \\
&> \frac{\beta}{\alpha + \beta} \P(\cE|p) + \frac{\alpha}{\alpha + \beta} \P(\cE'|q) \\
&= \frac{\beta}{\alpha + \beta} \P\left(\bigcap_{i=1}^{T} \cE_i\Big|p\right) + \frac{\alpha}{\alpha + \beta} \P\left(\bigcap_{i=1}^{T} \cE_i^\complement\Big|q\right) \\
&= \frac{\beta}{\alpha + \beta} \left( \P(\cE_1|p) \right)^T + \frac{\alpha}{\alpha + \beta} \left( \P(\cE_1^\complement|q) \right)^T \,,
\end{align*}
where in the last step, we use the fact that events are i.i.d., conditioned on the underlying measure.

To show that the probability of error is strictly larger than
$2\alpha \beta/(\alpha + \beta)$, it suffices to show that either
$\left( \cP(\cE_1|p) \right)^T > 2\alpha$ or $\left(
\cP(\cE_1^\complement|q) \right)^T > 2\beta$ unless $n$ is as
large, as claimed.
Consider the test $\phi_{\cE_1}$ on $m=n/T$ points that
outputs $p$ on the event $\cE_1$ and outputs $q$ on the event
$\cE_1^\complement$.
Thus, if $\cP(\cE_1|p) \leq (2\alpha)^{1/T} \leq \alpha'$ and
$\cP(\cE_1^\complement|q) \leq (2\beta)^{1/T} = \beta'$, the
test $\phi_{\cE_1}$ solves the hypothesis testing problem
$\cB_\pfht(p,q, \alpha', \beta')$ with $m$ samples.
By definition of the sample complexity, we have $m \geq \nstar_\pfht(p,q,
\alpha', \beta')$.
Combined with the definition of $n$ and $m$, we obtain $T \lceil \nstar/T \rceil \geq n \geq T \nstar_\pfht(p,q, \alpha', \beta') $, so $(\nstar/T) + 1 \geq \lceil \nstar/T \rceil  \geq \nstar_\pfht(p,q, \alpha', \beta')$.
\end{proof}

\paragraph{Accuracy amplification:}
We now turn our attention to the upper bound.
To mirror the style of our lower bound, which was error
amplification, we will establish the upper bound using success
amplification.

\begin{lemma}[Success amplification]
\label{lem:upper-bound-fre-tech}
For all positive integers $T$ such that $\max\left(\alpha^{1/T},\beta^{1/T}\right) \leq 1/4$,
we have
\begin{align}
\nstar_\pfht\left( p,q, \alpha, \beta \right) \lesssim  T \cdot \nstar_\pfht\left( p,q, \alpha^{1/T} , \beta^{1/T} \right)  \,.
\end{align}
\end{lemma}

\begin{proof}[Proof of \Cref{lem:upper-bound-fre-tech}]
Let $\alpha' = \alpha^{1/cT}$ and $\beta' = \beta^{1/cT}$, 
where $c = 1/32$ is from \Cref{fact:boost}, both of which are less
than $1/4$ by assumption.

Consider the test $\phi^*$ achieving the sample complexity
$n:=\nstar_\pfht(p,q, \alpha', \beta')$.
Design the test $\phi'$ that randomly divides the data of $nT$ samples
into $T$ buckets and runs $\phi^*$ on each
bucket, then outputs the majority of the $T$ individual
outputs.

By \Cref{fact:boost}, the probability of error of the boosted
procedure when the distribution is $p$ (respectively, $q$) is at
most $(\alpha')^{T/32} = (\alpha^{32 /T})^{T/32} = \alpha$
(respectively, $\beta$).
Thus, the sample complexity is at most $T \cdot
\nstar_\pfht(p,q,\alpha^{1/cT}, \beta^{1/cT})$.
The desired result will follow if we show that
$\nstar_\pfht(p,q,\alpha^{1/cT}, \beta^{1/cT})$ is at most a
constant times $\nstar_\pfht(p,q,\alpha^{1/T}, \beta^{1/T})$,
which follows from \Cref{prop:mild-reduction}, since $\alpha^{1/T}
\leq1/4$ and $\beta^{1/T} \leq 1/4$.
This completes the proof.
\end{proof}

\subsection{Proof of \Cref{prop:reduction}}
\label{SecPropReduction}

\begin{proof}[Proof of \Cref{prop:reduction}]
We first start with the prior-free setting.
Combining
\Cref{lem:upper-bound-fre-tech,lem:self-reduction-freq-hp}, we
obtain the following (for all $T \in \N$ such that $\max\left(
\alpha^{1/T}, \beta^{1/T} \right) \leq 1/4$, which holds under the
assumptions on $T$ from \Cref{prop:reduction}):
\begin{align}
\label{eq:proof-prop-red-1}
T \cdot \left(\nstar_\pfht \left(p,q,(2\alpha)^{1/T},(2 \beta)^{1/T}\right)-1
\right) \leq \nstar_\pfht(p,q,\alpha,\beta) \lesssim T \cdot  \nstar_\pfht \left(p,q,\alpha^{1/T},\beta^{1/T}\right)\,.
\end{align}
Since $\max\left((2\alpha)^{1/T},(2 \beta)^{1/T} \right)
\leq 1/4$, monotonicity of the sample complexity implies that
$$\nstar_\pfht \left(p,q,(2\alpha)^{1/T},(2
\beta)^{1/T}\right) \geq \nstar_\pfht \left(p,q,1/4,1/4\right)
\geq 2,$$
which further implies that the left-hand expression
in \Cref{eq:proof-prop-red-1} is further lower-bounded by \\
$0.5 T
\nstar_\pfht \left(p,q,(2\alpha)^{1/T},(2 \beta)^{1/T}\right)$.
Thus, we obtain
\begin{align*}
0.5 T \cdot \nstar_\pfht \left(p,q,(2\alpha)^{1/T},(2 \beta)^{1/T}\right) \leq \nstar_\pfht(p,q,\alpha,\beta) \lesssim T \cdot  \nstar_\pfht \left(p,q,\alpha^{1/T},\beta^{1/T}\right)\,.
\end{align*}
Note that $\max\left((2\alpha)^{1/T},(2
\beta)^{1/T},\alpha^{1/T},\beta^{1/T}\right) \leq 1/4$ and
$\log(1/x)/\log(1/2x) \lesssim 1$ for $x \leq 1/4$.
Combining these facts with the application of
\Cref{prop:mild-reduction} to the display equation above, we
conclude that the upper and lower bounds are within constants of
each other, so
\begin{align*}
\nstar_\pfht(p,q,\alpha,\beta) \asymp T \cdot  \nstar_\pfht \left(p,q,\alpha^{1/T},\beta^{1/T}\right)\,.
\end{align*}

We now explain the special choice of $T$: we set $T =
\left\lfloor \frac{\log(1/2 \alpha)}{\log 4}\right\rfloor$, which
is at least $1$, since $ \alpha\leq 1/8$.
We will now show that $\max ( (2 \alpha)^{1/T}, (2 \beta)^{1/T} )
\leq 1/4$, and additionally that $(2 \alpha)^{1/T} \in
[1/16,1/4]$.
Thus, we obtain
\begin{equation*}
(2 \beta)^{1/T} \leq (2 \alpha)^{1/T} =
e^{ \frac{\log(2 \alpha)}{T} } \leq e^{ \frac{\log(2 \alpha)}{
\log(1/2 \alpha) / \log 4 } } = e^{- \log 4} = 1/4.
\end{equation*}
Similarly, it can be shown that
\begin{equation*}
(2 \alpha)^{1/T} = e^{
\frac{\log(2 \alpha)}{T} } \geq e^{ \frac{\log(2 \alpha)}{
\log(1/2 \alpha) /(2 \log 4)} } = e^{-2 \log 4} = 1/16,
\end{equation*}
where we
use the fact that $T \geq 1$, so $T \geq \frac{\log(1/2
\alpha)}{2\log 4}$.
Thus, the result established above specializes to
$\nstar_\pfht(p,q,\alpha,\beta) \asymp T \nstar(p,q,
(\alpha)^{\frac{1}{T}} , (\beta)^{1/T} )$, and the latter is
equivalent to (up to constants) $\nstar_\pfht(p,q, 1/4,
\beta^{1/T})$ by \Cref{prop:mild-reduction}.
{We conclude that
\begin{align*}
\nstar_\pfht\left( p,q,\alpha,\beta \right) \asymp \log\left( \frac{1}{\alpha} \right) \nstar_\pfht\left( p,q, 1/4, \beta' \right)\,,
\end{align*}
proving the first part of \Cref{prop:reduction}.}

{We now prove the second part of the proposition by translating
the results from the prior-free setting to the Bayesian setting using \Cref{cor:equivalence-bayesian-pf}.}
For any positive integer $T$ such that $\left( \delta/\prior
\right)^{\frac{1}{T}} \leq 1/8$, we have
\begin{align*}
\nstar_\bay\left( p,q,\prior,\delta \right) 
&\asymp \nstar_\pfht\left( p,q, \frac{\delta}{\prior} , \delta \right) \tagm{ \text{by \Cref{cor:equivalence-bayesian-pf}, since $\delta \leq \prior/4$}}\\
& \asymp T \cdot \nstar_\pfht \left( p,q, \left( \frac{\delta}{\prior} \right)^{\frac{1}{T}}, \delta^{1/T} \right) \tagm{ \text{by part one, since $ (2\delta/\prior)^{1/T} \leq 2(\delta/\prior)^{1/T} \leq 1/4$}}\\
& \asymp T \cdot \nstar_\bay\left( p,q, \frac{ \delta^{\frac{1}{T}}}{\delta^{\frac{1}{T}} + \left( \delta/\prior \right)^{\frac{1}{T}}}, \frac{ \delta^{\frac{1}{T}}}{\delta^{\frac{1}{T}} + \left( \delta/\prior \right)^{\frac{1}{T}}} \cdot \left( \frac{\delta}{\prior} \right)^{\frac{1}{T}}  \right) \tagm{ \text{by \Cref{cor:equivalence-bayesian-pf}, since $(\delta/\prior)^{1/T} \leq 1/8$}} \\
& = T \cdot \nstar_\bay\left( p,q, \frac{ \prior^{\frac{1}{T}}}{ \prior^{\frac{1}{T}} + 1}, \frac{ \delta^{\frac{1}{T}}}{ \prior^{\frac{1}{T}} + 1}  \right) \\
&\asymp T \cdot \nstar_\bay\left( p,q, \prior^{\frac{1}{T}}, \delta^{\frac{1}{T}}  \right),
\end{align*}
where the last equivalence follows from
\Cref{prop:mild-reduction}, as follows: (i) both priors are less
than $1/2$, (ii) the ratio of the error probability to the prior is
$(\delta/\prior)^{\frac{1}{T}}$ is at most $1/4$, and (iii)
$\log(1/\delta_1)/\log(1/\delta_2) \asymp 1$ for $ \delta_1 =
\delta^{\frac{1}{T}}$ and $\delta_2 = \delta_1/(1 +
\prior^{\frac{1}{T}})$, since $\delta_1/2\leq \delta_2 \leq
\delta_1 \leq 1/4$.

We now explain the special choice of $T$: we set $T =
\left\lfloor \frac{\log(\prior/\delta)}{\log 8} \right\rfloor$,
which is at least $1$, and we also took
$\prior' = \prior^{\frac{1}{T}}$ and $\delta' =
\delta^{1/T}$.
In particular, $\frac{\delta'}{\prior'} =
\left(\frac{\delta}{\prior}\right)^{ \frac{1}{T}} \leq
e^{\frac{-\log(\prior/\delta) \log 8}{\log(\prior/\delta)}} \leq
1/8$.
Similarly, we see that $\frac{\delta'}{\prior'} \geq
\frac{1}{64}$.
Moreover, we have $\prior' \in (0,1/2)$, since $T \leq
\frac{\log(1/\prior)}{\log 8}$ because $\delta \geq \prior^2$, implying that
$\prior^{\frac{1}{T}} \leq e^{\frac{\log(\prior) \log
8}{ \log(1/\prior)}} < 1/2$.

Specializing the result above to this choice of $T$, we
obtain
\begin{align*}
\nstar_\bay(p,q,\prior,\delta) \asymp T \nstar_\bay\left( p, q, \prior', \delta' \right) \asymp T \nstar_\bay\left( p, q, \prior', \prior'/4 \right),
\end{align*}
where the last equivalence uses
\Cref{prop:mild-reduction}, since $\delta' \in (\prior'/64,
\prior'/8)$.
\end{proof}

\section{Bayesian testing problem: Linear regime} %
\label{sec:bayes-testing}

In this section, we establish upper and lower bounds for
the Bayesian simple binary hypothesis testing problem $\cB_\bay(p,q,\prior,\delta)$
when $\delta$ is a small constant fraction of $\prior$.
Observe that the following result recovers the case $\prior=1/2$ by 
noting that $\js{1/2}(p,q) \asymp \hel^2(p,q)$:

\begin{theorem}[Characterization of sample complexity of $\nstar_\bay$ in linear error probability regime]
\label{thm:main-bayes}
For any $\prior \in (0,0.5]$, we have
\begin{align*}
\left \lceil \frac{3}{16}\frac{\prior\log(1/\prior)}{\js{\prior}(p,q)} \right\rceil
 \leq \nstar_\bay(p,q,\prior, \prior/4) \lesssim  
\left\lceil \frac{380 \prior \log(1/\prior) }{\js{\prior}(p,q)} \right\rceil\,.
\end{align*}
In particular, if the condition $\nstar_\bay(p,q,\prior, \prior/4) > 1$ holds\footnote{This mild regularity condition necessarily happens when either (i) $\nstar_\pfht(p,q,1/4,1/4) > 1$, (ii) $\hel^2(p,q) \leq \helconst$, or (iii) $\js{\prior}(p,q) < (3/16) \prior \log(1/\prior)$, or (iv) $\cH_{1 - \frac{0.5 \log 2}{\log(1/\prior)}} (p,q) \leq \frac{1}{1024}$. See also the discussion in \Cref{sec:preliminaries}.}, we have
\begin{align*}
\nstar(p,q,\prior,\prior/4) \asymp \frac{\prior \log(1/\prior)}{\js{\prior}(p,q)}  \asymp \frac{1 }{ \cH_{1 - \frac{0.5 \log 2}{\log(1/\prior)}} (p,q)}\,. 
\end{align*}

\end{theorem}

In fact, our techniques give upper and
lower bounds on $\nstar_\bay(p,q,\prior,
\prior(1- \gamma))$ for all sufficiently small $\gamma$
(not just a large enough constant), which we will state in \Cref{app:weak-detection}.

\begin{proof}
We establish this result in three steps: First, we show a lower
bound in terms of $\js{\prior}(p,q)$ (cf.
\ \Cref{thm:lower-bound-bayes}); second, we show an upper bound in 
terms of $\cH_{\lambda}(p,q)$ for a particular value of $\lambda$ (cf.\ \Cref{thm:upper-bound-bayes});
and finally, we show that the upper and lower bounds are within
constants of each other (cf.\ \Cref{conj:the-two-div}).

We start with the lower bound for general $\prior$:
\begin{proposition}[Lower bound for $\nstar_\bay$]  
For any $\prior \in (0,0.5]$ and $\gamma \in (0,1)$, we have
\label{thm:lower-bound-bayes}
\begin{align*}
\nstar_\bay(p,q,\prior, \prior(1 - \gamma)) \geq \left \lceil \frac{ \prior \gamma \log\left( \frac{1 - \prior}{\prior} \right) + \prior^2\gamma^2}{\js{\prior}(p,q)} \right\rceil \,.
\end{align*}
In particular, for $\gamma = 3/4$, we have $\nstar_\bay(p,q,\prior,
\prior/4) \geq \left\lceil \frac{3}{16} \frac{\prior \log(1/\prior)}{
\js{\prior}(p,q)} \right\rceil$.
\end{proposition}
\begin{proof}[Proof of \Cref{thm:lower-bound-bayes}]
Let $\Theta, \widehat{\Theta}, X_1,\dots,X_n$ be from \Cref{prob:bay-intro}.
Since $\widehat{\Theta}$ is a function of $X_1,\dots,X_n$, 
it follows that $\Theta \to (X_1,\dots,X_n) \to \widehat{\Theta}$ is a
Markov chain.
By Fano's inequality (\Cref{fact:fano}), we have
\begin{align}
\label{eq:fano-lower-bd}
\hbin(P_e^{(n)}) \geq \hent(\Theta|X_1,\dots,X_n)\,.
\end{align}

We now calculate upper and lower bounds on $I(\Theta;
X_1,\dots,X_n)$ using two different expressions.
First, using i.i.d.\ distribution of $X_i$'s conditioned on $\Theta$, we obtain: 
\begin{align*}
I(\Theta; X_1,\dots,X_n) &= \hent(X_1,\dots,X_n) - \hent(X_1,\dots,X_n| \Theta) \\
&\leq \left(\sum_i \hent(X_i)\right) - \hent(X_1,\dots,X_n| \Theta) \\
&= n \hent(X_i) - \hent(X_1,\dots,X_n| \Theta) \\
&= \left(n \hent(\prior p + (1- \prior)q)\right)   \\
&\qquad-\left( \prior \hent(X_1,\dots, X_n| \Theta = p) + (1 -\prior)\hent(X_1,\dots,X_n| \Theta=q) \right)  \\
&= \left(n \hent(\prior p + (1- \prior)q)\right) - \left( \prior n \hent(p) + (1 -\prior)n \hent(q) \right)  \\
&= n \left( \hent(\prior p + (1- \prior)q) - \prior \hent(p) - (1- \prior) \hent(q) \right) \\
&= n I(X;\theta) = n \js{\prior}(p,q), \numberthis \label{eq:mi-upper-bd}
\end{align*}
where the last equality is by \Cref{prop:skew-symm}.
Combining this inequality with the fact that $$I(\Theta;
X_1,\dots,X_n) = \hent(\Theta) -
\hent(\Theta|X_1,\dots,X_n),$$
we obtain
\begin{align}
\hent(\Theta|X_1,\dots,X_n) \geq \hent(\Theta) - n \js{\prior}(p,q)\,.
\end{align}
Combining this with \Cref{eq:fano-lower-bd}, we obtain
\begin{align*}
\hbin(P_e^{(n)}) \geq \hent(\Theta) - n \js{\prior}(p,q),
\end{align*}
which implies that $n \geq \frac{\hbin(\prior) -
\hbin(P_e^{(n)})}{\js{\prior}(p,q)}$, using the fact that
$\hent(\Theta) = \hbin(\prior)$.
Finally, $P_e^{(n)} \leq \prior(1 - \gamma)$ implies
$\hbin(P_e^{(n)}) \leq \hbin\left(\prior\left( 1 - \gamma
\right)\right)$ by monotonicity of the binary entropy on
$[0,1/2]$.
The strong concavity of the binary entropy
function $\hbin(\cdot)$\footnote{Indeed, the first and second
derivates of $\hbin(x)$ are $\log\left( \frac{1-x}{x} \right)$
and $\frac{1}{x (1-x)}$, respectively.
} implies that
\begin{equation*}
\hbin(\prior) - \hbin(\prior(1 - \gamma))\geq  \prior \gamma \nabla \hbin(\prior) + \prior^2\gamma^2= \prior \gamma \log\left( \frac{1- \prior}{\prior} \right) + \prior^2 \gamma^2.
\end{equation*}
Thus, we obtain the desired lower bound of $n \geq \frac{\prior \gamma \log\left( \frac{1- \prior}{\prior} \right) + \prior^2\gamma^2}{\js{\prior}(p,q)}$. Since $n \in \N$, we also obtain the lower bound with the ceiling operation.

We now simplify the expression for $\gamma = 3/4$.
We first observe that $\prior \log( \ba /\prior) \geq 0.5
\prior \log(1/\prior)$ for $\prior \leq 1/3$.
Thus, we have $\prior \gamma \log\left( \frac{1- \prior}{\prior}
\right) + \prior^2\gamma^2 \geq (3/8)\prior \log(1/\prior)$
for $\prior \leq 1/3$.
For larger $\prior \in (1/3,1/2]$, we note that
$\prior^2 \gamma^2 \geq \prior \log(1/\prior) \frac{\prior
\gamma^2}{ \log 2} \geq \prior \log (1/\prior) (3/16)$, since
$\prior \geq 1/3$, $\gamma = 3/4$, and $\log 2 \leq 1$.

\end{proof}

We now establish a family of upper bounds in terms of
$\cH_\lambda$.

\begin{proposition}[Upper bounds]
\label{thm:upper-bound-bayes}
For any $\lambda \in (0,1)$ and $\prior \in (0,1/2)$, we have
\begin{align*}
\nstar_\bay(p,q,\prior, \prior(1 - \gamma)) \leq \left\lceil  \frac{ \lambda \log\left(\frac{\ba}{\prior}\right) + \log\left( \frac{1}{\bg} \right)}{\cH_{\bl}(p,q)} \right\rceil.
\end{align*}
In particular, for $\gamma = 3/4$, setting $\lambda =
\frac{0.5 \log 2}{\log(1/\prior)}$, which lies in $(0,0.5]$,
we
have
\begin{align}
\label{eq:bayesian-prior-gamma}
\nstar_\bay(p,q,\prior,  \prior/4) \leq \left\lceil \frac{2}{\cH_{\bl}(p,q)}\,\right\rceil,
\end{align}
recovering the case of uniform prior ($\prior = 0.5$)
with constant error probability (since $\lambda = 0.5$).%
\end{proposition}
\begin{proof}
Recall that $P_e^{(n)}$ denotes the minimum Bayesian error between
$p^{\otimes n}$ and $q^{\otimes n}$, which can be evaluated using \Cref{prop:E-gamma}.
Since $\min(x,y) \leq x^{\lambda}y^{1 - \lambda}$ for any $x, y \ge 0$ and $ \lambda \in (0,1)$, we
have
\begin{align*}
P_e^{(n)} &= \sum_{x_1,\dots,x_n} \min\left( \prior p^{\otimes n} (x_1,\dots,x_n), (1- \prior) q^{\otimes n}(x_1,\dots,x_n) \right) \\
&\leq  \sum_{x_1,\dots,x_n}  \left(\prior p^{\otimes n}(x_1,\dots,x_n)\right)^{\bl} \left((1- \prior) q(x_1,\dots,x_n) \right)^{\lambda} \\
&= \prior^{\bl}(1 -\prior)^{\lambda} \left( 1 - \cH_{\bl}(p^{\otimes n}, q^{\otimes n} ) \right) \\
&= \prior^{\bl }\ba^{\lambda} \left( 1 -  \cH_{\bl}(p, q) \right)^n \tag{\Cref{prop:hel-tensor}}\\
&\leq \prior^{\bl} \ba^{\lambda} \left( e^{- \cH_{\bl}(p, q)} \right)^n\,. \tag{$1-x \leq e^{-x}$ for all $x$}
\end{align*}
For $n \geq\left \lceil\frac{\log(\prior^{\bl} \ba^{ \lambda}/\delta)}{\cH_\lambda(p,q)}\right\rceil$, the right-hand side is at most $\delta$.
Finally, taking $\delta = \prior (1 - \gamma) = \prior \bg$,
the previous upper bound specializes to
\begin{align*}
\left\lceil\frac{\log(\prior^{\bl} \ba^{\lambda}/ \prior \bg)}{\cH_{\bl}(p,q)}\right\rceil 
= \left\lceil\frac{\log \left(\left( \frac{\ba}{\prior} \right)^{\lambda} \frac{1}{\bg} \right)}{\cH_{\bl}(p,q)}\right\rceil 
\leq \left\lceil\frac{ \lambda \log \left(\frac{1}{\prior} \right) + \log\left(\frac{1}{\bg} \right)}{\cH_{\bl}(p,q)}\right\rceil\,. 
\end{align*}
The specialization to $\gamma = 3/4$ follows by noting that
if we set $\lambda = \frac{0.5 \log 2}{\log(1/\prior)}$,
the numerator is upper-bounded by $0.5 \log(2) + \log 4 \leq
2$.

\end{proof}
Combining \Cref{thm:lower-bound-bayes,thm:upper-bound-bayes},
we obtain, for all $\prior \in (0,0.5]$, that
\begin{align}
\left\lceil \frac{3}{16} \frac{ \prior \log \left( \frac{1}{ \prior} \right)}{\js{\prior}(p,q)} \right\rceil \leq \nstar_\bay\left( p,q,\prior, \frac{\prior}{4} \right) \leq \left\lceil\frac{2}{\cH_{\bl}(p,q)}\right\rceil  \,,
\end{align}
whenever $\lambda = \frac{0.5 \log 2}{ \log(1/\prior)}$.
To show matching bounds on the sample complexity, we need to
establish the reverse inequality.
The following result shows that this is in fact true, 
and is the main technical result of the paper:
\begin{lemma}[Relationship between $\js{}$ and $\cH$]
\label{conj:the-two-div}
Let $p$ and $q$ be two distributions.
Let $\prior \in (0,1/2]$ and $\lambda \in (0,1/2]$.
Then
\begin{align}
\label{eq:conj-the-two-div}
\js{\prior}(p,q) \leq \frac{32 e^{2 \lambda \log(1/\prior)} \prior}{\lambda} 
\cH_{\bl}(p,q).
\end{align}
In particular, for $\lambda = \frac{0.5 \log
2}{\log(1/\prior)}$, we have
\begin{align*}
\left\lceil \frac{3}{16} \frac{ \prior \log \left( \frac{1}{ \prior} \right)}{\js{\prior}(p,q)} \right\rceil \leq \left\lceil\frac{2}{\cH_{\bl}(p,q)}\right\rceil  \leq \left\lceil \frac{256}{\log 2}\frac{  \prior \log(1/\prior)}{ \js{\prior}(p,q) } \right\rceil.
\end{align*}
\end{lemma}

We defer the proof of \Cref{conj:the-two-div} to the
end of this section and first focus on showing that
\Cref{thm:lower-bound-bayes,thm:upper-bound-bayes,conj:the-two-div}
suffice to prove \Cref{thm:main-bayes}.%

To this end, we obtain the following lower bound from
\Cref{thm:lower-bound-bayes}:
\begin{align}
\label{eq:main-thm-1}
\nstar(p,q,\prior,\prior/4) \geq \left\lceil \frac{3}{16} \frac{  \prior \log(1/\prior)}{\js{\prior}(p,q)} \right\rceil\,. 
\end{align}
We now establish the upper bound using
\Cref{thm:upper-bound-bayes}: for $\lambda = 0.5 \frac{\log
2}{\log(1/\prior)} \in (0,0.5]$, we have
\begin{align}
\label{eq:main-thm-2}
\nstar_\bay(p,q,\prior,\prior/4) \leq \left\lceil \frac{2}{\cH_{\bl}(p,q)} \right\rceil\, \leq 
\left\lceil \frac{ 256 \prior \log(1/\prior)}{ \log 2 \js{\prior}(p,q) } \right\rceil,
\end{align}
which combined with \Cref{conj:the-two-div}, gives the
desired result.%

The second statement in \Cref{thm:main-bayes} follows by noting
that if a positive integer $x\geq 2$ satisfies that $c_1 y \leq x
\leq \lceil c_2y \rceil $ for some $y, c_1,c_2 > 0$,
then $\lceil c_1 y \rceil \leq x \leq 2c_2 y$.
This is because if $x\geq 2$, then $c_2 y > 1$, so $\lceil
c_2y \rceil \leq 2 c_2 y$.
\end{proof}

\subsection{Proof of \Cref{conj:the-two-div}}

\begin{proof}[Proof of \Cref{conj:the-two-div}]
Our starting point is to note that both $\js{\prior}$ and $\cH_{\lambda}$ are $f$-divergences.
By the method of joint range~\cite{HarVaj11},
it suffices to establish \Cref{eq:conj-the-two-div} when
$p$ and $q$ are Bernoulli distributions (cf.\ \Cref{lem:joint-range}).
By abusing notation in this proof, we use $p$ and $q$ to
represent the biases of these Bernoulli distributions.

We set $\lambda = \frac{r}{\log(1/\prior)} \in (0,1/2]$ for some
$ r \in (0, 0.5\log(1/\prior)]$ in this proof.
Without loss of generality, we will assume that $p \leq 0.5$.
Going forward, we fix $p$ and $\prior$ to be arbitrary and take
$q$ to be a variable.
Define $j(q) := \js{\prior}(\Ber(p), \Ber(q))$ and $h(q):=
(\prior e^{2r}/r) \log(1/\prior) \cH_{\bl} (\Ber(p), \Ber(q))$.
Thus, our goal is to show that $j(q) \lesssim h(q)$.

In fact, we will establish a stronger result: there exists an
absolute constant $C_*$, independent of $p, q$, and $\prior$,
such that $C_*h(q) - j(q)$ is a nonnegative convex function; $C_* = 32$ suffices.
We argue in three steps:
\begin{enumerate}[label=(\roman*)]
\item $ch(p) - j(p)$ is equal to $0$, for all constants $c$.
\item $p$ is a stationary point of $ch(q) -  j(q)$, for all constants $c$.
\item There is an absolute constant $C_*$ such that $\frac{d^2}{dq^2} (C_* h(q) - j(q)) \geq 0$, for all $q$.
\end{enumerate}

The first claim is easy to verify since both expressions are
$f$-divergences.
The second claim is an immediate calculation, as shown below, and
the details are deferred to
\Cref{app:omitted_details_bayes_testing}.
\begin{restatable}[First derivatives of $h(\cdot)$ and $j(\cdot)$]{claim}{ClFirstDerivatives}
\label{cl:first-derivative}
\begin{align*}
\frac{d}{dq} \js{\prior}(\Ber(p), \Ber(q)) &= \ba \log\left( \frac{q}{\bq}  \frac{\prior \bp + \ba \, \bq}{\prior p + \ba q} \right)\\
\frac{d}{dq} \cH_{\bl} (\Ber(p), \Ber(q)) &= - \lambda p^{\bl}q^{-\bl} + \lambda \bp^{\bl} \bq^{-\bl}.
\end{align*}

\end{restatable}
Thus, the remaining proof is dedicated to showing the last claim
about the second derivatives,
stating that the second derivate of $j(q)$ is less than a
multiple of the second derivative of $h(q)$.
We now state the main technical part of this proof, with the proof of the lemma
deferred below.
\begin{lemma}
\label{lem:conj-two-div-hessian}
For any $p \in [0,0.5]$, $q \in (0,1)$,
$\prior \in [0,1/2]$, and $\lambda = \frac{r}{\log(1/\prior)}$, for some $r\geq 0$ such that $\lambda \in (0,0.5]$, we have
\begin{align}
\label{eq:main-technical-result}
\frac{d^2 }{d q^2} \js{\prior}(\Ber(p), \Ber(q)) \leq \frac{32 e^{2r}\prior}{r}  \cdot  \frac{d^2 }{d q^2} \left( \log(1/\prior) \cH_{\bl} (\Ber(p), \Ber(q)) \right)\,.
\end{align}

\end{lemma}
As argued above, \Cref{lem:conj-two-div-hessian} completes the
proof of \Cref{eq:conj-the-two-div}; the conclusion for $q$ being either exactly $0$ or $1$ follows from continuity of $j(\cdot)$ and $h(\cdot)$.
\end{proof}
We now provide the proof of \Cref{lem:conj-two-div-hessian}:
\begin{proof} [Proof of \Cref{lem:conj-two-div-hessian}]
We begin by stating \Cref{cl:hessian}, which calculates the second derivatives of these functions. Its proof is deferred to \Cref{app:omitted_details_bayes_testing}.

\begin{restatable}[Second derivatives of $h(\cdot)$ and $j(\cdot)$]{claim}{ClSecondDerivatives}
\label{cl:hessian}
Let $\delta = p - q$. Then
\begin{align}
\label{eq:hessian-lhs}
\frac{d^2 }{d q^2} \js{\prior}(\Ber(p), \Ber(q))
&=
\frac{\prior\ba\left(   \prior p \bp + \ba p \bq + \ba q \bp - \ba q \bq\right)}{q \bq (\prior p + \ba q)\left( \prior \bp + \ba \,\bq \right)}
= \frac{\prior \ba}{q \bq} \left( 1 + \frac{\ba \delta \left( \bq - q - \prior \delta \right)}{\left( q + \prior \delta\right)\left(\bq - \prior \delta  \right)} \right),
\end{align}
and if $\lambda = r/\log(1/\prior)$  and $\lambda \in (0,1/2]$, we have
\begin{align}
\frac{d^2 }{d q^2} \left(\log(1/\prior)  \cH_{\bl} (\Ber(p), \Ber(q)) \right)  &= \lambda \bl  \log(1/\prior) \left( \frac{1}{q} \left( \frac{p}{q} \right)^{\bl} + \frac{1}{\bq} \left( \frac{\bp}{\bq} \right)^{\bl} \right)\,\nonumber\\
\label{eq:hessian-rhs}
&\geq  \frac{r}{2} \left( \frac{1}{q} \left( \frac{p}{q} \right)^{\bl} + \frac{1}{\bq} \left( \frac{\bp}{\bq} \right)^{\bl} \right)\,.
\end{align}
\end{restatable}

In the sequel, the following approximation, proved in \Cref{app:proof-of-cl-linear-vs-nearly-linear}, will be especially useful:
\begin{restatable}{claim}{ClLinearVsNearlyLinear}
\label{cl:linear-vs-nearly-linear}
For all $x \geq 0$, $\prior \in (0,1/2]$, and $\lambda_* \in \left[0,\frac{r}{\log(1/ \prior)}\right]$, for some $r\geq 0$ such that $\lambda_* \in (0,0.5]$, we have
$1 + \frac{x}{1 + \prior x} \leq e^{2r} \left( 1 +  x \right)^{\bl_*}$.

\end{restatable}

We will establish \Cref{eq:main-technical-result}
using a case analysis.

\paragraph{Case 1: $q \leq p \leq  1/2$.}
In this case, $\delta:= p-q$ is nonnegative, and we use the
last expression in \Cref{eq:hessian-lhs}.
Additionally, we use the following bounds that can be easily verified: (i) $\bq - \prior
\delta \ge \bq - 0.5(0.5-q) \geq 0.5$, (ii) $|\bq - q - \prior
\delta| \leq 1$, and (iii) $\ba \delta \leq \delta$.
Starting with \Cref{eq:hessian-lhs} in \Cref{cl:hessian}, we
obtain
\begin{align*}
\frac{d^2 }{d q^2} \js{\prior}(\Ber(p), \Ber(q)) 
&= \frac{\prior \ba}{q \bq} \left( 1 + \frac{\ba \delta \left( \bq - q - \prior \delta \right)}{\left( q + \prior \delta\right)\left(\bq - \prior \delta  \right)} \right) \\
&\leq \frac{\prior \ba}{q \bq} \left( 1 + \frac{ \ba \delta \left| \bq - q - \prior \delta \right|}{\left( q + \prior \delta\right)\left(\bq - \prior \delta  \right)} \right) \tag{$\delta\geq 0$} \\
&\leq 2\frac{\prior}{q} \left( 1 + \frac{ 2\delta }{( q + \prior \delta)} \right) \tag{$\bq - \prior 	\delta\geq 0.5$ and $\bq \geq 0.5$}\\
&\leq 4\frac{\prior}{q} \left( 1 + \frac{(\delta/q) }{(1 + \prior (\delta/q))} \right) \\
&\leq 4 e^{2r}\frac{\prior}{q} \left( 1 + \frac{\delta}{q} \right)^{\bl} \tag{\Cref{cl:linear-vs-nearly-linear} and definition of $\lambda$}\\ 
&= 4 e^{2r} \prior\frac{1}{q} \left(\frac{p}{q} \right)^{\bl} \\ 
&\leq \frac{8 e^{2r} \prior}{r}  \frac{d^2 }{d q^2} \left( \log(1/\prior) \cH_{\bl} (\Ber(p), \Ber(q)) \right) ,
\end{align*}
where the last line uses \Cref{eq:hessian-rhs} in
\Cref{cl:hessian}.

\paragraph{Case 2: $q \geq 1/2$.}
Let $\tau = q - p = - \delta$, which we know is positive (since
$p \leq 1/2$).
Using the final expression in \Cref{eq:hessian-lhs} in
\Cref{cl:hessian}, we start with the following expression for the
second derivative:
\begin{align*}
\frac{d^2 }{d q^2} \js{\prior}(\Ber(p), \Ber(q)) 
&= \frac{\prior \ba}{q \bq} \left( 1 + \frac{\ba \tau \left( q - \prior \tau - \bq \right)}{\left( q - \prior \tau\right)\left(\bq + \prior \tau  \right)} \right)\,\\
&\leq \frac{\prior \ba}{q \bq} \left( 1 + \frac{\ba \tau \left|\bq + \prior \tau - q \right|}{\left( q - \prior \tau\right)\left(\bq + \prior \tau  \right)} \right) \tag{$\tau\geq 0$}\\
&\leq \frac{2\prior }{\bq} \left( 1 + \frac{4\tau}{\bq + \prior \tau  } \right) \tag{$q\geq 1/2$, $\ba \leq 1$, $q - \prior \tau \geq 1/4$ }\\
&\leq \frac{8\prior }{\bq} \left( 1 + \frac{\tau/\bq}{1 + \prior \tau/\bq  } \right) \\
&\leq \frac{8e^{2r} \prior}{\bq} \left( 1 + \tau/\bq \right)^{\bl} \tag{\Cref{cl:linear-vs-nearly-linear} and definition of $\lambda$}\\
&= 8 e^{2r} \prior \frac{1}{\bq} \left( \frac{\bp}{\bq}\right)^{\bl}\\
&\leq \frac{16e^{2r} \prior}{r} \frac{d^2 }{d q^2} \left( \log(1/\prior) \cH_{\bl} (\Ber(p), \Ber(q)) \right) ,
\end{align*}
where the last line uses \Cref{eq:hessian-rhs} in
\Cref{cl:hessian}.

\paragraph{Case 3: $p \leq q \leq  1/2$.}
In this case, we start with the first equality in
\Cref{eq:hessian-lhs} in \Cref{cl:hessian}:
\begin{align*}
\frac{d^2 }{d q^2} \js{\prior}(\Ber(p), \Ber(q)) 
&=
\frac{\prior\ba\left(   \prior p \bp + \ba p \bq + \ba q \bp - \ba q \bq\right)}{q \bq (\prior p + \ba q)\left( \prior \bp + \ba\, \bq \right)} \\
&\leq  \frac{\prior\ba\left(   \prior p \bp + \ba p \bq + \ba q q \right)}{q \bq (\prior p + \ba q)\left( \prior \bp + \ba\, \bq \right)} \tag{$\bp - \bq = q - p \leq q$} \\
&\leq  \frac{\prior \left(   \prior p \bp + \ba p \bq + \ba q^2 \right)}{q (1/2) (q/2)\left( 1/2 \right)} \tag{$\ba\geq 0.5$, $\bp \geq 0.5$, $\bq \geq 0.5$}\\
&\le  \frac{8 \prior\left(   \prior p \bp + \ba p \bq + q^2 \right)}{q^2}  \\
&\leq  8 \prior +  \frac{16 \prior p}{q^2} \tag{$\prior p \bp + \ba p \bq \leq 2 p$} \\
&\leq  8 \prior  \left(\frac{\bp}{\bq}\right)^{\bl} +  \frac{16 \prior p}{q^2} 
&& \left(\text{using $\bp \geq \bq$ and $\bl \geq 0$} \right)\\
&\leq 8 \prior \frac{1}{\bq} \left(\frac{\bp}{\bq}\right)^{\bl} + 16 \frac{\prior}{q} \left( \frac{p}{q} \right)^{\bl} && \left( \text{using $p \leq q$ and $\bl \in (0,1]$} \right) \\
&\leq 16\prior \left( \frac{1}{\bq} \left(\frac{\bp}{\bq}\right)^{\bl} + \frac{\prior}{q} \left( \frac{p}{q} \right)^{\bl} \right)  \\
&\leq \frac{32  \prior}{r}  \frac{d^2 }{d q^2} \left(\log(1/\prior) \cH_{\bl} (\Ber(p), \Ber(q)) \right) \tag{\Cref{eq:hessian-rhs}}\,.
\end{align*}
Thus, for each value of $q$, the desired inequality
holds, completing the proof of \Cref{lem:conj-two-div-hessian}.%
\end{proof}

\subsection{Proofs of \Cref{cl:first-derivative,cl:hessian}}
\label{app:omitted_details_bayes_testing}
\ClFirstDerivatives*

\begin{proof}
Since
\begin{equation*}
\js{\prior}\left( \Ber(p), \Ber(q) \right) = 
\prior \kl(\Ber(p), \Ber(\prior p + \ba a))) + \ba \kl \left( \Ber(q), \Ber(\prior p + \ba q) \right),
\end{equation*}
we obtain
\begin{align*}
\js{\prior}\left( \Ber(p), \Ber(q) \right) &= \prior p \log\left(\frac{p}{ \prior p + \ba q} \right) + \prior \bp \log\left( \frac{\bp}{ \prior \bp + \ba \bq} \right) \\
&\qquad +   \ba q \log\left(\frac{q}{ \prior p + \ba q} \right) + \ba \bq \log\left( \frac{\bq}{ \prior \bp + \ba \bq } \right),
\end{align*}
implying that
\begin{align*}
\frac{d}{dq}  \js{\prior}\left( \Ber(p), \Ber(q) \right) &=  -\prior p \frac{1}{\prior p + \ba q} \ba - \prior \bp  \frac{1}{\prior p + \ba \bq} (-\ba) \\
&\qquad+ \ba \log\left( \frac{q}{\prior p + \ba q}  \right) + \ba q \frac{\prior p + \ba q}{q} \left( \frac{(\prior p + \ba q) - \ba q }{(\prior p + \ba q)^2} \right) \\
&\qquad - \ba \log\left( \frac{\bq}{\prior \bp + \ba \bq}  \right) + \ba \bq \frac{\prior \bp + \ba \bq}{\bq} \left( \frac{-(\prior \bp + \ba \bq) +  \ba q }{(\prior \bp + \ba \bq)^2} \right)\\
&= -\frac{\prior \ba p}{\prior p + \ba q}  + \frac{\prior \ba p}{ \prior \bp + \ba \bq}    + \ba \log\left( \frac{q}{\bq}  \frac{\prior \bp + \ba \bq}{\prior p + \ba q} \right) +  \frac{\prior \ba p  }{\prior p + \ba q}  - \frac{\prior \ba \bp}{\prior \bp + \ba \bq}\\
&= \ba \log\left( \frac{q}{\bq}  \frac{\prior \bp + \ba \bq}{\prior p + \ba q} \right)\,.
\end{align*} 
Next, observe that $\cH_{\bl} (\Ber(p), \Ber(q)) = 1 - p^{\bl} q^{\lambda} -  \bp^{\bl} \bq^{\lambda}$.
Thus, the first derivative is
\begin{align*}
\frac{d}{dq} \cH_{\bl} (\Ber(p), \Ber(q)) = - \lambda p^{\bl}q^{-\bl} + \lambda \bp^{\bl} \bq^{-\bl}.
\end{align*}
\end{proof}

\ClSecondDerivatives*
\begin{proof}
Starting with $\js{\prior}$, we obtain the following, using the expression for the first derivative in \Cref{cl:first-derivative}:
\begin{align*}
\frac{d^2 }{d q^2} \js{\prior}(\Ber(p), \Ber(q)) &= \ba \frac{d^2 }{d q^2} \left( \log q - \log(\bq) + \log(\prior p + \ba \bq) - \log(\prior p + \ba q) \right)\\
&= \ba\left( \frac{1}{q} + \frac{1}{\bq} - \frac{\ba}{\prior p + \ba q} - \frac{\ba}{\prior \bp + \ba \bq} \right)\\
&= \ba \left( \frac{1}{q \bq} - \ba\left(  \frac{\prior \bp + \ba \bq + \prior p + \ba q}{\left( \prior p + \ba q \right)\left( \prior \bp + \ba \bq \right)}\right) \right) \\
&= \ba \left( \frac{1}{q \bq} - \frac{\ba}{\left( \prior p + \ba q \right)\left( \prior \bp + \ba \bq \right)}\right) \\
&= \ba \left( \frac{ \prior^2 p \bp + \prior \ba p\bq + \prior \ba q \bp + (\ba)^2 q\bq - \ba q \bq}{q \bq \left( \prior p + \ba q \right)\left( \prior \bp + \ba \bq \right)} \right) \\
&=  \prior \ba \left( \frac{ \prior p \bp + \ba p\bq +  \ba q \bp - \ba q \bq}{q \bq \left( \prior p + \ba q \right)\left( \prior \bp + \ba \bq \right)} \right) \,.
\end{align*}
We now show the derivation for the alternate expression given in the lemma. First, the numerator $\prior p \bp + \ba p\bq +  \ba q \bp - \ba q \bq$ simplifies to $p- \prior p^2 - 2 \ba pq + \ba q^2$.
Similarly, the expression $(q + \prior \delta)(\bq - \prior \delta) + \ba \delta(\bq - q - \prior \delta)$ simplifies to $p- \prior p^2 - 2 \ba pq + \ba q^2$.

Calculating the second derivative using the expression in \Cref{cl:first-derivative}, we obtain
\begin{align*}
\frac{d^2}{dq^2} \cH_{\bl} (\Ber(p), \Ber(q)) =   \lambda \bl p^{\bl}q^{-\bl-1} + \lambda \bl \bp^{\bl} \bq^{-\bl-1} = \lambda \bl \left( \frac{1}{q} \left( \frac{p}{q} \right)^{\bl} + \frac{1}{\bq} \left( \frac{\bp}{\bq} \right)^{\bl} \right) \,,
\end{align*}
which is always positive.
Since $\lambda =\frac{r}{\log(1/\prior)} \in (0,1/2)$, we have $\bl \geq 1/2$, so
\begin{align*}
\frac{d^2}{dq^2} \left(\log(1/\prior)\cH_{\bl} (\Ber(p), \Ber(q))\right) &\geq \log(1/\prior) \lambda \bl \left( \frac{1}{q} \left( \frac{p}{q} \right)^{\bl} + \frac{1}{\bq} \left( \frac{\bp}{\bq} \right)^{\bl} \right) \\
&\geq 0.5 \gamma \left( \left( \frac{p}{q} \right)^{\bl} + \frac{1}{\bq} \left( \frac{\bp}{\bq} \right)^{\bl} \right).
\end{align*}
\end{proof}

\subsection{Proofs of \Cref{thm:main-result-intro-bay,thm:prior-free-simple-hyp-testing-intro}}
\label{app:proof-main-result-intro-bay}

We restate \Cref{thm:main-result-intro-bay} below:
\ThmMainResultIntroBay*
\begin{proof}
The first regime follows from \Cref{thm:main-bayes}, the second regime follows from \Cref{prop:reduction}, and the final regime follows from \Cref{eq:bound-literature}.
\end{proof}

We now restate \Cref{thm:prior-free-simple-hyp-testing-intro}:
\ThmPriorFreeSBHT*
The proof of \Cref{thm:prior-free-simple-hyp-testing-intro} follows from \Cref{cor:equivalence-bayesian-pf}.

\section{Weak detection}
\label{app:weak-detection}
In this section, we consider the weak detection regime of the Bayesian hypothesis testing problem $\cB_\bay(p,q,\prior,\delta)$.
Recall that $\delta \geq \prior$ is the vacuous regime, and our 
results in \Cref{thm:main-result-intro-bay} give tight sample 
complexity bounds for $\delta \leq \prior/4$.
Here, we discuss the remaining regime of $\delta \in (\prior/4, \prior)$ with the special focus on $\delta \to \prior$ from below.
In particular, we parametrize $\delta$ to be $\delta = (1 - \gamma) \prior$, for $\gamma$ small enough.

As a starting point, we consider the uniform prior regime $\prior=1/2$.
Recall that in this vanilla setting, the sample complexity of strong detection $(\delta/\prior \leq 1/4)$ is characterized by the Hellinger divergence.
In the weak detection regime, however,  the following example suggests that the Hellinger divergence does not  characterize the  sample complexity: it could be either $\frac{\gamma^2}{\hel^2(p,q)}$ or $\frac{\gamma}{\hel^2(p,q)}$.

\begin{example}
\label{example:uniform-prior}
Let $c$ be a small enough absolute constant.
For any $\gamma \in (0, 1/4)$,
and $\epsilon \in (0, c \gamma^2 )$, there exist distributions $p$, $q$, $p'$, and $q'$ such that $\hel^2(p,q) \asymp \hel^2(p',q') \asymp \epsilon$,  such that  $\nstar_\bay\left( p,q,0.5, 0.5(1 - \gamma) \right) \asymp \frac{\gamma^2}{\hel^2(p,q)}$, but $\nstar_\bay\left( p',q',0.5, 0.5(1 - \gamma) \right) \asymp \frac{\gamma}{\hel^2(p',q')}$.
That is, the dependence on $\gamma$ is different.
\end{example}
\begin{proof}
We will use \Cref{prop:E-gamma} with $\prior = 1/2$ and thus $E_\gamma$ divergence simply becomes the total variation distance.

For $\delta \leq \gamma$ to be decided, 
let $p = \cN(\eps,1) $ and $q = \cN(-\delta,1)$, where $\cN(\mu,\sigma^2)$ denotes the univariate Gaussian distribution with mean $\mu$ and variance $\sigma^2$.
Let $u \in \R^n$ denote the all-ones vector.
Using spherical symmetry of isotropic Gaussians, we observe that 
\begin{align*}
  \dtv(p^{\otimes n}, q^{\otimes n}) &= 
  \dtv(\cN(\delta u, I), \cN(-\delta u, I)) \\
  &=\dtv(\cN(\sqrt{n}\delta, 1), \cN(-\sqrt{n}\delta, 1)) \\
  &\asymp \min (\sqrt{n} \delta, 1),
\end{align*}
where the last claim follows from \cite[Theorem 1.3]{DevMR18}.
In the context of \Cref{prop:E-gamma},
setting $(1 - \dtv(p^{\otimes n}, q^{\otimes n}))/2 = 1/2 - \gamma$, we see that $\nstar_\bay(p,q,1/2,(1 - \gamma)/2) =: \nstar$ has to be the smallest $n$ such that $\dtv(p^{\otimes n}, q^{\otimes n}) \asymp \gamma$.
 Therefore, the sample complexity $\nstar$ satisfies that $\min (\sqrt{\nstar} \delta, 1) \asymp \gamma$ and thus $\nstar \asymp \gamma^2/\hel^2(p,q)$, where $\delta^2 \asymp \dtv^2(p,q) \asymp \hel^2(p,q) \asymp \epsilon$.

  Let $p' = \Ber(1)$ and $q' = \Ber(1 - \epsilon)$.
  Then 
  \begin{align*}
    \dtv({p'}^{\otimes n}, {q'}^{\otimes n}) &= 
    1 - (1 - \epsilon)^n\,.
  \end{align*}
  Standard calculations shows that for $n \leq 1/ \epsilon$, the above expression is equivalent to $n \epsilon$ (up to constants).
  Setting $(1 - \dtv({p'}^{\otimes n}, {q'}^{\otimes n}))/2 = 1/2 - \gamma$, we
 see that the sample complexity is the smallest $n$ such that $\dtv({p'}^{\otimes n}, {q'}^{\otimes n}) \asymp \gamma$, which is then equivalent to $tn \asymp \gamma$, and hence $\nstar\left( p',q', 1/2, (1- \gamma)/2 \right) \asymp \gamma/ \epsilon \asymp \gamma/\hel^2(p',q')$, which also satisfies the constraint that $n \leq 1 /\epsilon$ used earlier in the approximation.
\end{proof}

Thus, the uniform prior regime already exhibits sufficiently different behavior in the weak detection regime.
In the weak detection regime, our bounds on the sample complexity from the previous section for general prior are as follows:

\begin{theorem}
\label{thm:bayes-all-regimes}
For any $\prior \in (0,0.5]$, $\gamma \in (0,1)$, and $\lambda \in (0,0.5]$ such that $\lambda \log(\ba/\prior) \leq \log(1/\bg)$, we have
\begin{align*}
\left\lceil \frac{\prior \gamma \log\left( \frac{1 - \prior}{\prior} \right) + \prior^2\gamma^2}{\js{\prior}(p,q)}\right\rceil  \lesssim \nstar_\bay(p,q,\prior,\prior(1 - \gamma)) \lesssim  
\left\lceil \frac{64 \prior\log(1/\bg) e^{2 \lambda \log(1/\prior) }  }{ \lambda \js{\prior}(p,q)} \right\rceil.
\end{align*}
Let $\gamma \in(0,c)$ for a small enough absolute constant $c$. Then the result above implies the following simplified expressions:
\begin{enumerate}
\item (Almost uniform prior) For $\prior = 1/2 - \eta$, for $\eta \leq 1/4$, we can take $\lambda \asymp \min(1, \gamma/ \eta)$ and use $\js{\prior}(p,q) \asymp \hel^2(p,q)$ for $\prior\in[1/4,1/2]$ to obtain
\begin{align}
\label{eq:weak-testing-almost-uni}	
\left\lceil \frac{\gamma\max(\gamma,\eta)}{\hel^2(p,q)} \right\rceil \lesssim \nstar_\bay(p,q,\prior,\prior(1- \gamma)) \lesssim \left\lceil\frac{\max(\gamma,\eta) }{\hel^2(p,q) }\right\rceil\,.
\end{align}

\item (Prior going to 0) When $\prior \leq 1/4$, by taking $\lambda \asymp \gamma/\log(1/\prior)$, we obtain
\begin{align}
\label{eq:weak-testing-small-prior}	
\left\lceil \frac{\prior \gamma \log(1/\prior)}{\js{\prior}(p,q)} \right\rceil \lesssim \nstar_\bay(p,q,\prior,\prior(1 -\gamma)) \lesssim \left\lceil\frac{\prior \log(1/\prior)}{\js{\prior}(p,q)}\right\rceil\,.
\end{align}
\end{enumerate}
Moreover, there exist distributions $p$ and $q$ such that the upper bounds in \Cref{eq:weak-testing-small-prior,eq:weak-testing-almost-uni} are achieved (cf.\ \Cref{ex:upper-bounds-tight-weak-detection}).
\end{theorem}
In both cases, we see that the upper and lower bounds
differ by a factor of $\gamma$.
Perhaps surprisingly, when $\prior$ is far from uniform (much further
than $\gamma$), the upper bound on the sample complexity is
independent of $\gamma$.
That is, \textit{extremely weak detection (when $\gamma\to 0$) seems to be as hard as
strong detection (when $\gamma = 9/10$).
}
We give explicit examples where this indeed happens:

\begin{example}[Upper bounds for sample complexity are tight in the worst case]
\label{ex:upper-bounds-tight-weak-detection}
For every $\epsilon \in (0,1/4)$, there exist distributions $p$ and $q$
such that $\js{\prior}(p,q) \asymp \prior \epsilon$.
Moreover, for every  $\prior \in (0,1/2) $ and $\gamma \in (0,1/4)$ 
such that $ \frac{\log\left(\frac{\ba}{\prior \bg}  \right)}{\epsilon} \geq 1$,
the sample complexity satisfies
$\nstar\left( p,q,\prior,\prior(1 - \gamma) \right) \asymp \frac{\log(\frac{\ba}{\prior \bg})}{\epsilon}$.
In particular, 
the following holds:
\begin{enumerate}
\item If $\prior = 1/2 - \eta$, then
$\nstar_\bay(p,q, \prior, \prior(1 - \gamma)) \asymp \frac{\max(\gamma,\eta)}{\js{\prior}(p,q)} \asymp \frac{\max(\gamma,\eta)}{\hel^2(p,q)}$.
Thus, the upper bound in \Cref{eq:weak-testing-almost-uni} is tight.

\item If $\prior \leq 1/4$, then $
\nstar_\bay(p,q, \prior, \prior(1 - \gamma)) \asymp \frac{\prior \log(1/\prior)}{\js{\prior}(p,q)}$. 
Thus, the upper bound in \Cref{eq:weak-testing-small-prior} is tight.

\end{enumerate}
\end{example}
\begin{proof}
 Let $p = \Ber(0)$ and $q = \Ber(\epsilon)$.
Then a quick calculation shows that $\js{\prior} \asymp \prior \epsilon$:
\begin{align*}
\js{\prior}(p,q) &= \prior \kl\left(\Ber(0), \Ber(\ba \epsilon)\right) + \ba
\kl\left(\Ber(\epsilon),\Ber(\ba \epsilon) )\right) \\
&= \prior \log\left(\frac{1}{1 - \ba \epsilon}\right) +
\ba \cdot \overline{\epsilon} \log\left( \frac{\overline{\epsilon}}{1 - \ba \epsilon}\right) + \ba \epsilon \log\left( \frac{1}{\ba} \right)\,.
\end{align*}
Since $\prior\leq 1/2$ and $\epsilon\leq 1/2$, the first term is of the order $\prior \epsilon$,\footnote{Note that $\log\left( \frac{1}{1 - \ba \epsilon} \right) = \log\left( 1 + \frac{\ba \epsilon}{1 - \ba \epsilon} \right) \asymp \frac{\ba \epsilon}{1 - \ba \epsilon} \asymp \epsilon$.} so it implies the desired lower bound on $\js{\prior}(p,q)$. To upper-bound the second term, we note that the second term is negative, while the third term is at most $\prior \epsilon$, since $\ba\epsilon \log(1/\ba) = \ba \epsilon \log\left( 1 + \prior/\ba \right) \leq \epsilon \prior$. Thus, we have $\js{\prior}(p,q) \asymp \prior \epsilon$.

Moreover, the average error of hypothesis testing with $n$
samples is equal to (cf.\ \Cref{prop:E-gamma})
\begin{align*}
P_e &= \sum_{x_1,\dots,x_n} \min\left( \prior p^{\otimes n} (x_1,\dots,x_n), (1- \prior) q^{\otimes n}(x_1,\dots,x_n) \right) \\ 
&= \min\left( \prior, (1- \prior) (1 - \epsilon)^n \right),
\end{align*}
which is less than $\prior(1 - \gamma)$ if and only if $(1-
\prior) (1 - \epsilon)^n \leq \prior(1 - \gamma)$, which happens only if
$n \asymp \frac{\log(\frac{\ba}{\prior \bg})}{\log(1/
\overline{\epsilon})}$.
Since $\epsilon \leq 1/2$, the denominator $\log(1/\overline{\epsilon})$ is equivalent to $\epsilon$, up to constants.

We evaluate the resulting expression $n \asymp \frac{\log(\frac{\ba}{\prior \bg})}{\epsilon}$ in the different regimes of $\prior$.
For $\prior \leq 1/4$, the numerator is equivalent to $\log(1/\prior)$, so $n \asymp
\frac{\log(1/\prior)}{ \epsilon} \asymp
\frac{\prior\log(\frac{1}{\prior})}{\js{\prior}(p,q)}$.
In the other regime of $\prior = 1/2 - \eta$ for $\eta \in (0,1/4)$, the numerator $\log(\frac{\ba}{\prior
\bg})$ satisfies
\begin{equation*}
\log(\frac{\ba}{\prior
\bg}) = \log\left( 1 + \frac{4 \eta}{1 - 2 \eta} \right) + \log\left( 1 + \frac{\gamma}{\bg} \right) \asymp \eta + \gamma \asymp \max(\eta,\gamma).
\end{equation*}
Thus, the sample complexity is $\frac{\max(\eta,\gamma)}{ \js{\prior}(p,q)} \asymp \frac{\max(\eta,\gamma)}{ \hel^2(p,q)}$.

\end{proof}

\section{Distributed hypothesis testing} 
\label{sec:distributed_hypothesis_testing}

In this section, we provide the proofs of our results in \Cref{sub:distributed-simple-binary}.

\subsection{Hypothesis testing under communication constraints}
\label{app:hypothesis-testing-under-communication-constraints}
We provide the proof of \Cref{thm:comm-hyp-testing}, recalled below:

\ThmCommHypTesting*

We shall use the following result from \cite{BhaNOP21}:

\begin{theorem}[\cite{BhaNOP21}]
\label{thm:bhanop21}
Let $p$ and $q$ be two distributions over a domain $\cX$ of size $k$.
Let $\Theta$ be a binary random variable over $\{p,q\}$ with $\P(\Theta = p) = \prior$. Let $X|\Theta = \theta$ be distributed as $\theta$.
Let $D \in \N$ be a communication budget greater than $1$. 
Then there exists an algorithm $\cA$ that takes as input $p$, $q$, 
$\prior$, and $D$ and outputs a function $f: \cX \to [D]$, in time $\poly(k, |D|)$, such that
\begin{align*}
I\left(\Theta; f(X)\right) \gtrsim I(\Theta;X) \cdot \min\left( 1,   \frac{D}{\log\left( 1/ I(\Theta;X) \right)} \right).
\end{align*}
\end{theorem}

\begin{proof}[Proof of \Cref{thm:comm-hyp-testing}]
Let $X$ be distributed as $X_1$ in \Cref{prob:bay-intro}.

We will choose $f:\cX \to [D]$ shortly.
Let $p'$ and $q'$ be the distributions of $f(Y)$ for $Y$ distributed as $p$ and $q$, respectively.
The server solves $\cB_\bay(p',q',\prior,\delta)$, and as a result, requires $\nstar_\bay(p',q',\prior,\delta)$ samples.
Since $\nstar_{\bay, \comm}(p,q,\prior,\delta,D) \leq \nstar_\bay(p',q',\prior,\delta)$, in the remainder of the proof, we will upper-bound $\nstar_\bay(p',q',\prior,\delta)$.

First consider the regime $\delta \in (\prior/100,\prior/4]$.
Let $f: \cX \to [D]$ be the quantizer from \Cref{thm:bhanop21} with $\prior$, which can be computed in polynomial-time using \Cref{thm:bhanop21}.
Then \Cref{thm:main-result-intro-bay} implies that 
\begin{align*}
 \nstar_\bay(p',q',\prior,\delta) &\asymp \frac{\prior \log(1/\prior)}{I(\Theta; f(X_1))} \\
&\lesssim \frac{\prior \log(1/\prior)}{I(\Theta;X)}  \max\left( 1, \frac{\log(1/I(\Theta;X))}{D} \right) \tag*{(using \Cref{thm:bhanop21})} \\
&\asymp \nstar_\bay\left( p,q,\prior,\delta \right) \max \left( 1, \frac{\log\left(\frac{\nstar_\bay\left( p,q,\prior,\delta \right)}{\prior \log(1/\prior)}\right)}{D} \right) \tag*{(using \Cref{thm:main-result-intro-bay}).}
\end{align*}

Now consider the sublinear regime $\delta \in (\prior^2, \prior/100)$.
Following the notation of \Cref{thm:main-result-intro-bay}, let $f: \cX \to [D]$ be the quantizer from \Cref{thm:bhanop21} with $\prior = \prior'$.
Then \Cref{thm:main-result-intro-bay} implies that 
\begin{align*}
 \nstar_\bay(p',q',\prior,\delta) &\asymp \log\left( \prior/\delta \right) \nstar_\bay\left( p',q', \prior', \prior'/4 \right)\\
&\lesssim \log\left( \prior/\delta \right) \nstar_\bay\left( p,q,\prior',\prior'/4 \right) \max\left( 1, \frac{ \log\left( \frac{ \nstar_\bay\left(p,q,\prior',\prior'/4  \right)}{\prior'\log(1/\prior')} \right)}{D} \right) \\
&\asymp \nstar_\bay\left( p,q,\prior,\delta \right) \max\left( 1, \frac{ \log\left( \frac{ \nstar_\bay\left(p,q,\prior,\delta  \right)}{\log(\prior/\delta) \prior'\log(1/\prior')} \right)}{D} \right) \\
&\asymp \nstar_\bay\left( p,q,\prior,\delta \right) \max\left( 1, \frac{ \log\left( \frac{ \nstar_\bay\left(p,q,\prior,\delta  \right)}{ \prior'\log(1/\prior)} \right)}{D} \right) \\
&\leq \nstar_\bay\left( p,q,\prior,\delta \right) \max\left( 1, \frac{ \log\left( \frac{ \nstar_\bay\left(p,q,\prior,\delta  \right)}{ \prior\log(1/\prior)} \right)}{D} \right).
\end{align*}
The final regime $\delta\leq \prior^2$ follows analogously by noting that $\hel^2(p,q) \asymp I(\Theta';X)$ when $\Theta'$ is uniform over $\{p,q\}$.
Choosing the $f$ from \Cref{thm:bhanop21} for the uniform prior, we obtain
\begin{align*}
\nstar(p',q',\prior,\delta) &\lesssim \frac{\log(1/\delta)}{I(\Theta'; f(X))} \lesssim \frac{\log(1/\delta)}{I(\Theta;X)} \max\left( 1, \frac{\log(1/I(\Theta';X))}{D} \right) \\
&\lesssim \nstar_\bay\left( p,q,\prior,\delta \right) \max\left( 1, \frac{\log\left( \frac{\nstar_\bay\left( p,q,\prior,\delta \right)}{\log(1/\delta)} \right)}{D} \right).
\end{align*}

\end{proof}

\subsection{Hypothesis testing under local differential privacy}
\label{app:hypothesis-testing-under-local-differential-privacy}

In this section, we prove \Cref{thm:computational-cost-privacy}, restated below:

\ThmCompCostPrivacy*

For an $\epsilon$-LDP mechanism $\cM:\cX \to \cY $ and a distribution $p$ over $\cX$, we use $\cM(p)$ to denote the distribution of $\cM(X)$ for $X \sim p$. 
We shall use the following result from \cite{PenAJL23}:

\begin{theorem}[{\cite[Corollary 1.18]{PenAJL23}}]
\label{thm:PenAJL23}
Let $p$ and $q$ be two distributions on $[k]$.
For any $\ell \in  \N$, let $\cC_{\ell}$ be the set of all $\epsilon$-LDP mechanisms that maps $\cX$ to $[\ell]$.
Let $\cA = \{ \left( \cM(p), \cM(q) \right): \cM \in \cC_\ell\}$.
Let $g(\cdot,\cdot): \cA \to \R$ be a jointly quasi-convex function.
Then, there is an algorithm that takes as input $p, q, \epsilon, \ell$, and returns an $\epsilon$-LDP mechanism $\cM$ that maximizes $g(\cM(p),\cM(q))$ over $\cM \in \cC_\ell$.
Moreover, the algorithm returns in time polynomial in $k^{\ell^2}$ and $2^{\ell^3 \log l}$. 
\end{theorem}

\begin{proof}[Proof of \Cref{thm:computational-cost-privacy}]
Observe that the optimal $\epsilon$-LDP mechanism $\cM'$ is the one that minimizes $\nstar_\bay\left( \cM'(p), \cM'(q),\prior,\delta \right)$.
Let $\nstar_{\bay,\priv} := \nstar_{\bay,\priv}\left( p,q,\prior,\delta,\epsilon \right)  =  \nstar_\bay\left( \cM'(p), \cM'(q),\prior,\delta \right) $ be the private sample complexity.
Letting $\cY'$ be the range of $\cM$, we will further quantize it using a (deterministic) function $f: \cY' \to [\ell] $. 
If we choose the best quantizer, \Cref{thm:comm-hyp-testing} implies that the resulting sample complexity $n'$ is at most
\begin{align*}
n' \lesssim \nstar_{\bay,\priv} \cdot \max\left( 1, \frac{ \log\left( \frac{\nstar_{\bay,\priv}} {\prior} \right)  }{l} \right)\,,
\end{align*}
yielding the desired bound on the sample complexity.
We now show the runtime guarantees using \Cref{thm:PenAJL23}.
Let $\cC_\ell$ be the set of all $\epsilon$-LDP mechanisms that map $\cX \to [\ell]$, which in particular includes $f(\cM')$ analyzed above.
Observe that for all the regimes of $\delta$,
minimizing $\nstar_\bay(\cM(p), \cM(q),\prior,\delta)$ over 
$\cM \in \cC_\ell$ is equivalent (up to constant factors) to
 maximizing a jointly convex function $a \cH_{\lambda}( \cM(p), \cM(q))$ for some $\prior > 0$ and $\lambda \in (0,1)$ depending on the parameters $\prior$ and $\delta$;
 the joint convexity follows because $\cH$ is a jointly convex function. 
Let $\cM''$ be the $\epsilon$-LDP mechanism returned by \Cref{thm:PenAJL23} for maximizing the jointly convex function above. Thus, the resulting sample complexity is less than a constant multiple of $n'$.
\end{proof}

\section{Sequential hypothesis testing}
\label{section:sequential}

The sequential hypothesis testing problem, investigated by
Wald~\cite{Wald45}, considered a hypothesis testing setting where the statistician observes samples sequentially and chooses to make a decision after they have gathered enough evidence. 
In contrast to the prior-free hypothesis testing from \Cref{prob:pfht-intro}, the number of samples is not fixed in advance---and indeed, depends stochastically on the observed samples---which leads to the \emph{expected} sample complexity as a metric of interest. We formally define the problem below:
\begin{definition}\label{def: seq}
Let $p$ and $q$ be two distributions over a discrete domain $\cX$. %
A sequential hypothesis test $\phi$ consists of a stopping rule $T$ that is a stopping time (i.e., the event $\{T=t\}$ depends only on the samples observed up to time $t$) and a decision $\phi(X_1, \dots, X_T) \in \{p, q\}$. We say that $\phi$ solves the simple sequential hypothesis testing problem with type-I error $\alpha$, type-II error $\beta$, and expected sample complexity $n$, if
\begin{align*}
\P_p(\phi(X_1, \dots, X_T) \neq p) \le \alpha, \quad \P_q(\phi(X_1, \dots, X_T) \neq q)\le \beta, \quad{ and } \quad \max\{\E_p T, \E_q T\} \le n.
\end{align*}
We use $\cB_\seq(p,q,\alpha,\beta)$ to denote this problem, and define its sample complexity $\nstar_\seq(p,q,\alpha,\beta)$ to be the smallest $n$ such that there exists a test $\phi$ which solves $\cB_\seq(p,q,\alpha,\beta)$.
\end{definition}

The counterpart to the Neyman--Pearson test in the sequential setting is Wald's sequential probability ratio test (SPRT). This test updates the log-likelihood ratio $\Lambda_t$ according to the rule $\Lambda_t = \sum_{i=1}^t \log \frac{p(x_i)}{q(x_i)}$, where $x_i$ is the $i$-th observation. The test stops and declares $p$ when $\Lambda_t$ hits an upper threshold $A$, or stops and declares $q$ when it hits a lower threshold $B$. The thresholds $A$ and $B$ are computed based on the desired type-I and type-II errors, for example, setting $A = \log \frac{1-\beta}{\alpha}$ and $B = \log \frac{\beta}{1-\alpha}$ suffices to get the desired errors. Wald and Wolfowitz~\cite{WaldWolf48} showed that the SPRT (with optimal thresholds) is \emph{optimal}, in the following sense. Let $\phi$ be the SPRT with stopping rule $T$ and errors $\alpha$ and $\beta$, and let $\phi'$ be any with stopping rule $T'$ and errors $\alpha' \le \alpha$ and $\beta' \le \beta$. Then we must necessarily have $\E_p [T'] \ge \E_p[T]$ and $\E_q [T'] \ge \E_q[T]$.%

If the Type-I error is at most $\alpha$ and the Type-II error is
at most $\beta$, Wald~\cite{Wald45} showed the following lower bounds on the
expected sample size:
\begin{align*}
\E_p T \ge \frac{\kl(\alpha, \overline \beta)}{\kl(p, q)} \quad \text{and} \quad \E_q T \ge \frac{\kl(\beta, \overline \alpha)}{\kl(q, p)},
\end{align*}
where we used $\kl(x, y)$ to denote $\kl(\Ber(x),
\Ber(y))$. This gives the sample complexity lower bound
\begin{align}\label{eq:wald-lb}
\nstar_\seq(p,q,\alpha,\beta) \geq \max \left\{\frac{\kl(\alpha, \overline \beta)}{\kl(p, q)}, \frac{\kl(\beta, \overline \alpha)}{\kl(q, p)} \right\}.
\end{align}
Since any fixed sample-size algorithm is also a valid sequential algorithm, we also have the simple upper bound
\begin{align}\label{eq: wald ub}
\nstar_\seq(p,q,\alpha,\beta) \leq \nstar_\pfht(p, q, \alpha, \beta).
\end{align}

The upper and lower bounds need not match. A simple example is when the two distributions have mismatched supports, leading to unbounded KL divergences. 
\begin{example} \label{example: seq}
Consider two distributions $p = (0, 1/2, 1/2)$ and
$q = (1/2, 1/2, 0)$.
The fixed-size sample complexity to get errors of both types at most $\delta$ is $\Theta(\log(1/\delta)/\hel^2(p, q)) \asymp \log(1/\delta)$. The lower bound on the sequential sample complexity from \Cref{eq:wald-lb} does not depend on $\delta$ and equals 0, since the KL divergence is infinite. 
More interestingly, the \emph{exact} sequential sample complexity also does not depend on $\delta$. 
In fact, it is possible to get 0 errors of both types by taking $2$ samples in expectation. 
\end{example}
This example shows that the sample complexity behaviour between the fixed sample-size case and the sequential case can be quite different. It is unclear if the expected sample size for the optimal sequential test is characterized by any $f$-divergence in general. In what follows, we address some special cases of interest.

\subsection{Constant error regime}

When $\alpha$ and $\beta$ are constants (say $1/8$), the sample complexity of sequential hypothesis testing turns out to be identical to that of the prior-free hypothesis testing problem. This is summarized in the following theorem:

\begin{theorem}\label{thm: seq_constant}
Let $p$ and $q$ be as in \Cref{def: seq}. Then $\nstar_\seq(p, q, 1/8, 1/8) \asymp \frac{1}{\hel^2(p,q)}$.
\end{theorem}
\begin{proof}
Clearly, $\nstar_\seq(p, q, 1/8, 1/8) \le \nstar_\pfht(p, q, 1/8, 1/8) \asymp \frac{1}{\hel^2(p,q)}$ since a fixed-size test is also a valid sequential test. It remains to show the reverse inequality. Suppose there exists a sequential test with sample complexity $N$ that attains errors of both types at most $1/8$. By Markov's inequality, the probability that  the sample-size of this test exceeds $8 N$ is at most $1/8$. Thus, we may construct a fixed size test as follows: take $8N$ samples and declare the result of the sequential test if there is one, otherwise declare 0. By the union bound, the probabilities of type-I and type-II errors of this test are most 1/4. This means $8N \ge \nstar_\pfht(p, q, 1/4, 1/4) \asymp \frac{1}{\hel^2(p,q)}$, and so $N \gtrsim \frac{1}{\hel^2(p,q)}$. Since this is true for any sequential test with errors bounded by 1/8, we conclude that $\nstar_\seq(p,q,1/8,1/8) \gtrsim \frac{1}{\hel^2(p,q)}$. This concludes the proof.
\end{proof}

This proof technique fails when the desired errors are much smaller than a constant, say $\delta$, because the application of Markov's inequality leads to a lower bound that is $\asymp \delta \nstar_\pfht(p, q, \delta, \delta)$, which is far away from the upper bound of $\nstar_\pfht(p, q, \delta, \delta)$. 

\subsection{Bounded likelihood ratio}
One way to avoid the situation in \Cref{example: seq} is considering distribution pairs that have bounded likelihood ratios. This ensures that the KL divergence is finite, and the lower bound~\eqref{eq:wald-lb} is non-trivial. The assumption of bounded log-likelihood simplifies the relationships between various $f$-divergences, as noted in the following lemma:
\begin{lemma}
\label{lemma: bounded_likelihood}
Let $p$ and $q$ be two distributions supported on a discrete space
$\cX$.
Suppose there exists a constant $C > 1$ such that
\[
\frac{1}{C} \le \frac{p(x)}{q(x)} \le C,
\]
for all $x \in \cX$.
Then the following hold:
\begin{enumerate}
\item
$\hel^2(p, q) \asymp \kl(p \| q)$,
\item
For $\alpha \le 1/2$, $\js{\alpha}(p, q) \asymp \alpha \hel^2(p, q)$,
\end{enumerate}
\end{lemma}

\begin{proof}[Proof of \Cref{lemma: bounded_likelihood}]
The proof of the first part follows from Lemma 5 in \cite{BirMas98}, where it is shown that if $p(x)/q(x) \le C$ then
\[
\hel^2(p, q) \le \kl(p, q) \le \hel^2(p, q) (2 + \log C).
\]
Note that the lower bound holds even without the bounded likelihood assumption.

The second part follows by expressing $\js{\alpha}(p, q)$ in
terms of the KL-divergence, and using part 1 above:
\begin{align*}
\js{\alpha}(p,q) &= \alpha \kl(p, \alpha p + \overline \alpha q) + \overline 
\alpha \kl(q, \alpha p + \overline \alpha q)\\
&\stackrel{(a)}\asymp \alpha \hel^2(p,\alpha p + \overline \alpha q) + \overline \alpha \hel^2(q, \alpha p + \overline \alpha q)\\
&\stackrel{(b)}\le 2\alpha \overline \alpha \hel^2(p, q)\\
&\asymp \alpha \hel^2(p, q).  
\end{align*}
Here, $(a)$ follows from part 1, and $(b)$ follows from
the convexity of $\hel^2$ in its arguments.
To show the other direction,
\begin{align*}
\js{\alpha}(p,q) &= \alpha \kl(p, \alpha p + \overline \alpha q) + \overline 
\alpha \kl(q, \alpha p + \overline \alpha q)\\
&\ge \alpha \kl(p, \alpha p + \overline \alpha q)\\
&\ge \alpha \hel^2(p, \alpha p + \overline \alpha q).
\end{align*}
Now observe that $\hel^2(p, \alpha p + \overline \alpha q)$ is
a monotonically decreasing function of $\alpha$ for fixed $p$ and
$q$.
Since $\alpha \le 1/2$, we have
\begin{align*}
\hel^2(p, \alpha p + \overline \alpha q) &\ge \hel^2(p,  (p + q)/2)\\
&\ge 0.29 \hel^2(p, q),
\end{align*}
where the last line follows from Lemma 5 in
\cite{BirMas98}.
This leads us to conclude
\begin{align*}
\js{\alpha}(p, q) \gtrsim \alpha \hel^2(p, q),
\end{align*}
which completes the proof of part 2.
\end{proof}

Let $\alpha, \beta < 1/4$. Using \Cref{lemma: bounded_likelihood}, it is easy to see that the lower bound from \eqref{eq:wald-lb} is simply
\begin{align*}
\nstar_\seq(p, q, \alpha, \beta) \ge \frac{\log(1/\min\{\alpha,
\beta\})}{\hel^2(p, q)}.
\end{align*}
The upper bound \eqref{eq: wald ub} gives 
\begin{align*}
\nstar_\seq(p, q, \alpha, \beta) &\le \nstar_\pfht(p, q, \alpha, \beta)\\
&\le \nstar_\pfht(p, q, \min\{\alpha, \beta\}, \min\{\alpha, \beta\})\\
&\asymp \frac{\log(1/\min\{\alpha,
\beta\})}{\hel^2(p, q)}\,.
\end{align*}
These two bounds match, giving  a tight characterization of the worst-case expected
sample size under the bounded-likelihood assumption, summarized in the following theorem:

\begin{theorem}[Sample complexity for sequential hypothesis testing under the bounded likelihood-ratio condition]\label{thm: seq_bdd}
Let $p$ and $q$ be two distributions on a discrete set $\cX$ that satisfy the bounded likelihood-ratio assumption. Then $\nstar_\seq(p,q,\alpha, \beta) \asymp \frac{\log(1/\min\{\alpha,
\beta\})}{\hel^2(p, q)}.$
\end{theorem}

\begin{remark}
Our main result on Bayesian sample complexity identifies the
sample complexity $\nstar_\bay(p, q, \prior, \delta)$ for
different regimes of $\delta \le \prior/4$ for $\prior \in (0,1/2]$.
Using \Cref{lemma: bounded_likelihood}, it can be easily
verified that under the bounded likelihood-ratio assumption, the sample complexity is given by
\[
\nstar_\bay(p, q, \prior, \delta) \asymp \frac{\log(1/\delta)}{\hel^2(p, q)},
\]
that is, the trivial upper bound is tight.
A similar conclusion holds for the prior-free setting, where we obtain that the sample complexity for $\alpha, \beta \in (0,1/8]$ is
\[
\nstar_\pfht(p, q, \beta, \alpha) \asymp \frac{\log(1/\min\{\alpha, \beta\})}{\hel^2(p, q)}.
\]
\end{remark}

\section{Hypothesis testing with erasures}
\label{section:abstention}

Hypothesis testing with erasures (also called \emph{abstention} or \emph{rejection}) allows the statistician to declare an output in the range $\{p, q, \star\}$, where the $\star$ indicates insufficient evidence to make a decision. 
A decision of $\star$ leads to an error event called the \emph{erasure error}, whereas an incorrect decision (that is not an erasure) is called an \emph{undetected error}. 
When erasures are not allowed, the usual type-I and type-II errors fall under the undetected error category. The rationale behind the $\star$ option is that is some settings, the cost of making an undetected error is significantly higher than that of making an erasure error. 
Thus, the statistician is willing to incur some probability of erasure error to reduce the probability of undetected errors. 
This added flexibility can allow for tests with small undetected errors using potentially much fewer samples than tests without the $\star$ option. %

To give a concrete example, consider the prior-free hypothesis testing problem with distributions $p$ and $q$, and suppose we require the type-I and type-II errors to be at most $\beta$.
Without erasures, the sample complexity is $\asymp \frac{\log(1/\beta)}{\hel^2(p,q)}$. 
Now suppose the statistician is allowed to declare $\star$ with probability at most a small constant, say $1/10$, under either hypothesis, while keeping probabilities of undetected errors at most $\beta$ under both hypotheses. 
Will the sample complexity of this new statistical problem be meaningfully smaller than $\frac{\log(1/\beta)}{\hel^2(p,q)}$? We answer this question in the affirmative, showing that the new sample complexity is $\asymp \max\{\nstar_\pfht(p, q, 1/10, \beta), \nstar_\pfht(p, q, \beta, 1/10)\}$. This quantity can be arbitrarily smaller than $\frac{\log(1/\beta)}{\hel^2(p,q)}$, as shown in \Cref{example: erasure}. We generalise this special case and identify the sample complexity in the prior-free and Bayesian settings with erasure (to be defined later) in almost all error regimes. 

\subsection{Background on simple hypothesis testing with erasures}\label{section: background_erasure}

In what follows, we denote the simple binary hypothesis testing problem without erasures as the \emph{standard setting} and with erasures as the \emph{erasure setting}.
We briefly describe known results in the erasure setting. The problem of trading off erasures versus undetected errors was studied in Forney~\cite{For68} for binary (as well as $M$-ary) hypothesis testing. 
Consider a (Bayesian) hypothesis $\theta$ taking values in $\{p,q\}$ and a \emph{single} observation $x$ that comes from $p$ or $q$. 
Let $\P(x,\theta)$ be the joint probability of choosing hypothesis $\theta$ and observing $x$ under $\theta$. 
For simplicity, let $\overline \theta$ denote the second hypothesis. There are three error events to consider: an \emph{error}, denoted by $\cE_\tot$, an \emph{undetected error}, denoted by $\cE_\und$, and an \emph{erasure error} denoted by $\cE_\era$. 
An error is when the estimate is not the correct hypothesis; i.e., it is either $\overline \theta$ or $\star$. 
An undetected error is when estimate is the incorrect hypothesis $\overline \theta$. An erasure error is when the estimate is the erasure $\star$. 
Observe that if the test declares $p$ in $R_p \subseteq \cX$, $q$ in $R_q \subseteq \cX$, and $\star$ in $(R_p \cup R_q)^c$, then the probabilities of these errors are given by
\begin{align*}
\P(\cE_\tot) &= \sum_{\theta \in \{p,q\}} \sum_{x \in R_\theta^c} \P(x,\theta),\\
\P(\cE_\und) &= \sum_{\theta \in \{p,q\}} \sum_{x \in R_\theta} \P(x,\overline \theta), \quad \text{ and }\\
\P(\cE_\era) &= \P(\cE_\tot) - \P(\cE_\und).
\end{align*}
Forney~\cite{For68} considers the problem of optimally trading off $\P(\cE_\tot)$ versus $\P(\cE_\und)$, as this captures the tradeoff between declaring erasures and undetected errors. Using the Neyman--Pearson lemma, Forney~\cite{For68} established the optimality of a certain family of tests, where optimality means no other test can simultaneously better both $\P(\cE_\tot)$ and $\P(\cE_\und)$. This family of tests is characterised by thresholds $\eta > 1$ as follows:
\begin{align}
\phi(x) = 
\begin{cases}
p &\text{ if } \frac{p(x)}{q(x)} \ge \eta,\\
q &\text{ if } \frac{q(x)}{p(x)} \ge \eta,\\
\star &\text{ otherwise.}
\end{cases}
\end{align}
Since $\eta > 1$, the decision regions $R_p$ and $R_q$ are disjoint. %
We were unable to find a complete analysis of the trade-offs between $\P(\cE_\tot)$ and $\P(\cE_\und)$ in the asymptotic regime (similar to the Stein and Chernoff regimes in the standard setting). For completeness, we have included this in \Cref{app: erasure}. This analysis may be of independent interest. 

We now formally define the simple binary hypothesis testing problem with erasures. 

\begin{definition}[Bayesian simple binary hypothesis testing with erasures]\label{def: bayesian_erasure}
Let $p$ and $q$ be two distributions over a discrete domain $\cX$, and let $\prior, \delta \in (0,1)$ and $\lambda > 0$.
We say that a 
test $\phi: \cup_{n=0}^\infty \cX^n \to \{p,q, \star\}$
\textit{solves the Bayesian simple binary hypothesis testing with erasures problem} of distinguishing $p$ versus $q$, under the prior
$(\prior, 1- \prior)$, with sample complexity $n$ and $\lambda$-weighted error $\delta$, if for $X := (X_1,\dots,X_n)$, we have
\begin{align}
\label{eq: bayesian_erasure}
\prior \cdot \P_p \left\{   \phi(X) \neq p\right\} + (1 - \prior) \cdot \P_{q} \left\{   \phi(X) \neq q\right\} + \lambda \P(\phi(X) = \star) \leq \delta.
\end{align}
Here, we use $\P_p$ and $\P_q$ to denote $\P_{X \sim p^{\otimes n}}$ and $\P_{X \sim q^{\otimes n}}$, respectively.
Let $\cB_\baye(p,q,\prior,\lambda,\delta)$  denote this problem and define its sample complexity $\nstar_\baye(p,q,\prior, \lambda,\delta)$
to be the smallest $n$ such that there exists a test $\phi$ which
solves $\cB_\baye(p,q,\prior,\gamma, \delta)$. \end{definition}

\begin{remark}[Remark on the notation]
When $p$ and $q$ are clear from context, we shall also use the notation $\cB_\baye(\prior, \lambda, \delta)$ and $\nstar_\baye(\prior, \lambda, \delta)$. For a given test $\phi$ that uses $n$ samples, we shall use the notation $\cE_\tot^n$, $\cE_\und^n$, and $\cE_\era^n$ to denote the error event, the undetected error event, and the erasure error event, respectively.
For a test $\phi$, we shall also use the shorthand $\P^\phi(\cdot)$, $\P_p^\phi(\cdot)$ and $\P_q^\phi(\cdot)$ to denote the probabilities of the error events for $\phi$. With this notation, we can rewrite~\eqref{eq: bayesian_erasure} as
\begin{align*}
\lambda \P^\phi(\cE_\tot^n) + \overline \lambda \P^\phi(\cE_\und^n) \leq \delta, \quad \text{ or } \quad \P^\phi(\cE_\und^n) + \lambda \P^\phi(\cE_\era^n) \leq \delta.
\end{align*}
Note that $\lambda$ could be arbitrarily large; however, we show in \cref{lemma: np_erasure} that when $\lambda \ge 1/2$, the optimal decision rule can disregard erasures altogether. Intuitively, this is because the since the cost of declaring an erasure is too high. Thus, for all intents and purposes, we may restrict to $\lambda < 1/2$, which is also the more practical regime.
\end{remark}

We also consider a prior-free version, where 
all possible errors are analysed separately.
\begin{definition}[Prior-free simple binary hypothesis testing with erasures]
\label{def: pf_erasure}
Let $p$, $q$, and $\phi$ be as in \Cref{def: bayesian_erasure}. We say that $\phi$ \emph{solves the simple binary
hypothesis testing with erasures problem} of $p$ versus $q$ with type-I error
$\alpha$, type-II error $\beta$, erasure errors $\delta_0$ and $\delta_1$, and sample complexity $n$, if
\begin{align}
\P_p^\phi(\phi(X) = q) \le \alpha, \quad \P_q^\phi(\phi(X) = p) \le \beta, \quad \P_p^\phi(\phi(X) = \star) \le \delta_0, \quad \text{ and } \P_q^\phi(\phi(X) = \star) \le \delta_1.
\end{align}
We use $\cB_\pfe(p,q,\alpha,\beta, \delta_0, \delta_1)$ to denote this
problem, and define its sample complexity $\nstar_\pfe(p,q,\alpha,\beta, \delta_0, \delta_1)$ to be the smallest $n$ such that
there exists a test $\phi$ which solves
$\cB_\pfe(p,q,\alpha,\beta, \delta_0, \delta_1)$. When the context is clear, we shall drop $p$ and $q$ and use the notation $\cB_\pfe(\alpha, \beta, \delta_0, \delta_1)$ and $\nstar_\pfe(\alpha, \beta, \delta_0, \delta_1)$.
\end{definition}

\subsection{Prior-free hypothesis testing with erasures}

Our main contribution in this section is to find a formula for $\nstar_\pfe(\alpha, \beta, \delta_0, \delta_1)$ when the parameters lie in $(0,1/8]$. %

\begin{theorem}[Sample complexity of simple binary hypothesis testing with erasures]\label{thm: sc_pf_erasure}
\label{thm:prior-free-erasure}
Let $p$ and $q$ be distributions on a finite alphabet $\cX$. Let $\alpha$, $\beta$, $\delta_0$, and $\delta_1$ all lie in $(0,1/8]$. Then the prior-free sample complexity of the binary hypothesis testing problem with erasures is given by
\begin{align*}
\nstar_\pfe(\alpha, \beta, \delta_0, \delta_1) \asymp \max\{\nstar_\pfht(\alpha, \beta+\delta_1), \nstar_\pfht(\alpha+\delta_0, \beta)\},
\end{align*}
where we used the shorthand $\nstar_\pfht(\alpha, \beta) := \nstar_\pfht(p,q,\alpha,\beta)$ to denote the sample complexity of the standard binary hypothesis testing problem (without erasures) with type-I error at most $\alpha$ and type-II error at most $\beta$.
\end{theorem}

\begin{remark}[Benefit of erasure]
\label{rem:benefit-erasure}
We noted earlier that the sample complexity of simple binary hypothesis testing when the type-I and type-II errors are at most $\beta$ is $\asymp \frac{\log(1/\beta)}{\hel^2(p,q)}$. If erasure probabilities of $\delta_0 = \delta_1 = 1/10$ are allowed, the above result shows that the new sample complexity is roughly $\max(\nstar_\pfht(0.1,\beta), \nstar_\pfht(\beta,0.1))$ $\asymp \frac{\beta\log(1/\beta)}{\min\{\js{\beta}(p, q), \js{\beta}(q,p)\}}$. 
This can be much lower than the sample complexity without erasures. 

\begin{example}\label{example: erasure}
Let $p = (1-\epsilon, \epsilon, 0)$ and $q = (1-\epsilon, 0, \epsilon)$. Observe that $\hel^2(p,q) \asymp \epsilon$, and so $\nstar_\pfht(p,q, \beta, \beta) \asymp \frac{\log(1/\beta)}{\hel^2(p,q)}$. Now suppose we are allowed to declare an erasure with some fixed probability like $1/10$. Using the above result, the sample complexity for getting type-I and type-II errors at most $\beta$ will be characterized by ${\min\{\js{\beta}(p, q), \js{\beta}(q,p)\}}$. Observe that 
\begin{align*}
\js{\beta}(p,q) &= H(\beta p + \bar \beta q) - \beta H(p) - \bar \beta H(q)\\
&= H(\beta p + \bar \beta q) - H(p)\\
&= \left[\bar \epsilon \log \frac{1}{\bar \epsilon} + \beta \epsilon \log \frac{1}{\beta\epsilon} + \bar\beta \epsilon \log \frac{1}{\bar\beta \epsilon}\right] - \left[\epsilon \log \frac{1}{\epsilon} + \bar \epsilon \log \frac{1}{\bar \epsilon}\right]\\
&= \epsilon \beta \log \frac{1}{\beta} + \epsilon \bar \beta \log \frac{1}{\bar \beta}\\
&\asymp \epsilon \beta \log \frac{1}{\beta}.
\end{align*}
By symmetry, $\js{\beta}(q,p) = \js{\beta}(p,q)$. Thus, the sample complexity with erasures is
\begin{align*}
\nstar_\pfe(p,q,\beta, \beta, 0.1) \asymp \frac{\beta \log(1/\beta)}{\epsilon \beta \log \frac{1}{\beta}} = \frac{1}{\epsilon}.
\end{align*}
Observe that this is much favourable compared to $\log(1/\beta)/\epsilon$, which is the sample complexity without erasures.
\end{example}
\end{remark}

\begin{proof}[Proof of \Cref{thm:prior-free-erasure}]
We start with the lower bound.
Suppose a test $\phi$ solves $\cB_\pfe(\alpha, \beta, \delta_0, \delta_1)$ with sample complexity $\nstar_\pfe(\alpha, \beta, \delta_0, \delta_1)$. We can modify this test to a new test $\hat \phi$ that outputs a decision in the set $\{p,q\}$ by always declaring $p$ when the original decision is $\phi(X) = \star$. This new test is prior-free hypothesis testing problem \emph{without} erasures with type-I and type-II errors controlled by
\begin{align*}
\P_p(\hat \phi(X) = q) \le \alpha, \quad \text{ and } \quad \P_q(\hat \phi(X) = p) \le \beta+\delta_1.
\end{align*}
Thus, we must have the lower bound 
\begin{align*}
\nstar_\pfe(\alpha, \beta, \delta_0, \delta_1) \ge \nstar_\pfht(\alpha, \beta + \delta_1).
\end{align*}
Using a similar argument but declaring $\hat \phi(X) = q$ when $\phi(X) = \star$, we conclude the lower bound
\begin{align*}
\nstar_\pfe(\alpha, \beta, \delta_0, \delta_1) \ge \max \{ \nstar_\pfht(\alpha, \beta + \delta_1), \nstar_\pfht(\alpha + \delta_0, \beta) \}.
\end{align*}
We now turn our attention to the upper bound.
Observe that any test in the standard setting that controls type-I and type-II errors by $\alpha$ and $\beta$, respectively, is also a valid test for the erasure setting. This gives the upper bound 
\begin{align*}
\nstar_\pfe(\alpha, \beta, \delta_0, \delta_1) \le \nstar_\pfht(\alpha, \beta). 
\end{align*}
Now if either $\delta_0 \le \alpha$ or $\delta_1 \le \beta$, the lower and upper bounds are tight. To see this, suppose $\delta_0 \le \alpha$, then 
$$\nstar_\pfht(2\alpha, \beta) \le \nstar_\pfht(\alpha+\delta_0, \beta) \le \nstar_\pfht(\alpha, \beta). $$ 
Since $\max\{2\alpha, \alpha, \beta\} \le 1/4$, we may apply \Cref{prop:mild-reduction} to conclude that the sample complexity is $\asymp \nstar_\pfht(\alpha, \beta)$. A similar argument also works for $\delta_1 \le \beta$. 

Thus, we may restrict our attention to the case with $\delta_0 > \alpha$ and $\delta_1 > \beta$.  (This is also the more relevant regime in practice.) To show a matching upper bound, we will construct a test with sample complexity $\nstar_\pfht(\alpha/3, (\beta + \delta_1)/3) +  \nstar_\pfht((\alpha + \delta_0)/3, \beta/3)$ that achieves the necessary bounds on the prior-free problem with erasures. Our test is the following: 
\begin{enumerate}
\item
Let $\phi_1$ be a binary hypothesis test without erasures that attains type-I and type-II errors of at most $\alpha/3$ and $(\beta + \delta_1)/3$, respectively, using $N_1:=\nstar_\pfht(\alpha/3, (\beta + \delta_1)/3)$ samples. Sample $N_1$ points and let the first decision $D_1 = \phi_1(X_1, \dots, X_{N_1})$.
\item
Let $\phi_2$ be a binary hypothesis test without erasures that attains type-I and type-II errors of at most $(\alpha+\delta_0)/3$ and $\beta/3$, respectively, using $N_2:=\nstar_\pfht((\alpha+\delta_0)/3, \beta/3)$ samples. Sample $N_2$ points and let the second  decision $D_2 = \phi_2(X_{N_1+1}, \dots, X_{N_1+N_2})$.
\item
Our final output $D := \phi(X_1, \dots, X_{N_1+N_2})$ is given as follows. If $D_1 = D_2 = p$, declare $D = p$. If $D_1 = D_2 = q$, declare $D = q$. Otherwise, declare $\star$.
\end{enumerate}
 We now evaluate all the errors for this test. The undetected error under $p$ is bounded as
\begin{align*}
\P_p(D = q) = \P_p(D_1 = q)\P_p(D_2 = q)
\le \frac{\alpha}{3} \cdot \frac{\alpha+\delta}{3}
\le \frac{\alpha}{3}
\le \alpha.
\end{align*}
The undetected error under $q$ is bounded as
\begin{align*}
\P_q(D = p) &= \P_q(D_1 = p)\P_q(D_2 = p)
\le \frac{\beta+\delta_1}{3} \cdot \frac{\beta}{3}
\le \frac{\beta}{3}
\le \beta.
\end{align*}
The erasure error under $p$ is bounded as
\begin{align*}
\P_p(D = \star) &= \P_p(D_1 = p)\P_p(D_2 = q) + \P_p(D_1 = q)\P_p(D_2 = p)\\
&\le 1 \cdot \frac{\alpha+\delta_0}{3} + \frac{\alpha}{3} \cdot 1
= \frac{2\alpha+\delta_0}{3}
\le \delta_0. 
\end{align*}
Similarly, the erasure error under $q$ is bounded as
\begin{align*}
\P_q(D = \star) &= \P_q(D_1 = p)\P_q(D_2 = q) + \P_q(D_1 = q)\P_q(D_2 = p)\\
&\le \frac{\beta+\delta_1}{3}\cdot 1 + 1 \cdot \frac{\beta}{3} = \frac{2\beta + \delta_1}{3} \le \delta_1.
\end{align*}

Thus, our test achieves the desired errors using $N_1+N_2$ samples, which is easily seen to be
\begin{align*}
N_1+N_2 &= \nstar_\pfht(\alpha/3, (\beta + \delta_1)/3) + \nstar_\pfht((\alpha+\delta_0)/3, \beta/3)\\
&\asymp \nstar_\pfht(\alpha, \beta + \delta_1) + \nstar_\pfht(\alpha+\delta_0, \beta) \tag{by \Cref{prop:mild-reduction}}\\
&\asymp \max \{\nstar_\pfht(\alpha, \beta + \delta_1), \nstar_\pfht(\alpha+\delta_0, \beta) \}. 
\end{align*}
This concludes the proof.
\end{proof}

\subsection{Bayesian hypothesis testing with erasures}

The sample complexities of the Bayesian and prior-free hypothesis testing problems in the erasure setting satisfy a relation similar to that in the standard setting noted in \Cref{claim:reln-bayesian-prior-free}. Throughout this section, we shall assume without loss of generality that $\prior \le 1/2$.

\begin{claim}[Relation between Bayesian and prior-free sample complexities in the erasure setting]
\label{claim: bayesian_priorfree_erasure}
For any $\prior, \beta, \delta_0, \delta_1 \in (0,1)$ and $\lambda > 0$, we have
\begin{align}\label{eq: erasure_pf_b}
\nstar_\pfe\left(\frac{\delta}{\prior}, 
\frac{\delta}{\overline \prior}, \frac{\delta}{\prior\lambda}, \frac{\delta}{\lambda \overline \prior}\right) \le \nstar_\baye\left(\prior, \lambda, \delta \right) \le \nstar_\pfe\left(\frac{\delta}{4\prior}, 
\frac{\delta}{4\overline \prior}, \frac{\delta}{4\prior\lambda}, \frac{\delta}{4\lambda \overline \prior}\right).
\end{align}
\end{claim}

\begin{proof}
Note that the $\lambda$-weighted Bayes error of a test $\phi$ being at most $\delta$ means
\begin{align*}
\prior \P_p(\phi(X) = q) + \overline \prior \P_q(\phi(X) = p) + \lambda \prior \P_p(\phi(X) = \star) + \lambda \overline \prior \P_q(\phi(X) = \star) \le \delta.
\end{align*}
The sample complexity lower bound follows by noting that any test in the Bayesian setting whose $\lambda$-weighted error is at most $\delta$ must have each of the four terms on the left hand side bounded by $\delta$, as well. The sample complexity upper bound follows by noting that any prior-free test that has its undetected and erasure error probabilities bounded by the parameters on the right has each of the four terms at most $\delta/4$, and so its $\lambda$-weighted Bayes error is at most $\delta$.
\end{proof}

This result, combined with \Cref{thm: sc_pf_erasure} and \Cref{prop:mild-reduction} immediately yields the following corollary:

\begin{corollary}[Bayes sample complexity in a limited regime]\label{cor: baye_limited}
Without loss of generality, assume $\prior \le 1/2$. Let $\lambda > 0$ and let $\delta \le \lambda\prior/8$. Then $\nstar_\baye(\prior, \lambda, \delta) \asymp \nstar_\pfe\left(\frac{\delta}{\prior}, 
\frac{\delta}{\overline \prior}, \frac{\delta}{\prior\lambda}, \frac{\delta}{\lambda \overline \prior}\right)$.
\end{corollary}
\begin{proof}
We may apply \Cref{prop:mild-reduction} and \Cref{thm: sc_pf_erasure} to conclude that the lower and upper bounds in expression~\eqref{eq: erasure_pf_b} are within constants of each other if $$\max \left\{\frac{\delta}{\prior}, 
\frac{\delta}{\overline \prior}, \frac{\delta}{\prior\lambda}, \frac{\delta}{\lambda \overline \prior}\right\} \le 1/8,$$
which holds when $\delta \le \lambda\prior/8$.
\end{proof}

This result, though useful, is not completely satisfactory because of the limited range of $\delta$. Note that the Bayesian problem is only worth solving with $\delta \le \max\{\prior, \lambda\}$, since the constant predictors $\phi(X) = q$ and $\phi(X) = \star$ achieve errors of $\prior$ and $\lambda$, respectively. Analogous to the standard setting, the regime of interest in the erasure setting should therefore be $\delta \in (0, \min\{\prior, \lambda\}/4]$. This regime may be much larger that $(0, \prior\lambda/8]$, for example, when $\lambda \asymp \prior$.

In what follows, we develop an alternate approach to solve the sample complexity problem in the larger regime of interest. Our next lemma identifies the optimal Bayes test.

\begin{lemma}[Optimal Bayes test]\label{lemma: np_erasure}
Let $p$ and $q$ be two distributions on a finite set $\cX$, and let the prior distribution on the hypotheses be $(\prior, \overline \prior)$ for some $\prior \le 1/2$. For any 
single-sample
test $\phi: \cX \to \{p, q, \star\}$, consider the $\lambda$-weighted Bayes error $\lambda \P^\phi(\cE_\tot) + \overline \lambda \P^\phi(\cE_\und)$, which equals $\lambda \P^\phi(\cE_\era) + \P^\phi(\cE_\und)$. Then the following tests minimise the $\lambda$-weighted Bayes error:
\begin{enumerate}[label=(\roman*)]
\item
When $\lambda \ge 1/2$, the optimal test $\phi$ declares no erasures and is the same as the optimal Bayes test in the standard setting. To be precise,
\begin{align*}
\phi(x) = 
\begin{cases}
p \quad &\text{ if } \frac{\prior p(x)}{\overline \prior q(x)} \ge 1\\
q \quad &\text{ if } \frac{\prior p(x)}{\overline \prior q(x)} < 1.
\end{cases}
\end{align*}

\item
When $\lambda < 1/2$, the following test $\phi$ is optimal:
\begin{align*}
\phi(x) = 
\begin{cases}
p \quad &\text{ if } \frac{\prior p(x)}{\overline \prior q(x)} \ge \frac{\overline \lambda}{\lambda}\\
q \quad &\text{ if } \frac{\prior p(x)}{\overline \prior q(x)} \le \frac{\lambda}{\overline \lambda}\\
\star &\text{ otherwise.}
\end{cases}
\end{align*}
\end{enumerate} 
\end{lemma}
\begin{proof}
It is enough to look at deterministic hypothesis tests. 

\paragraph{Proof of (i):} Consider any test $\phi: \cX \to \{p, q, \star\}$ with decision regions $R_p$, $R_q$, and $\cE_\era$, where the test declares $p$, $q$, and $\star$ respectively. We can write the Bayes error as
\begin{align*}
&\lambda \prior \P_p(\cE_\era) + \lambda \overline \prior \P_q(\cE_\era) + \prior \P_p(R_q) + \overline \prior \P_q(R_p)\\
& \ge \frac{ \prior \P_p(\cE_\era) +  \overline \prior \P_q(\cE_\era)}{2} + \prior \P_p(R_q) + \overline \prior \P_q(R_p) \tag{$\lambda \ge 1/2$}\\
&\ge \min\{\prior \P_p(\cE_\era), \overline \prior \P_q(\cE_\era) \} + \prior \P_p(R_q) + \overline \prior \P_q(R_p).
\end{align*}
Suppose $\prior \P_p(\cE_\era) \le \overline \prior \P_q(\cE_\era)$. We show that there is a test without erasure that achieves the lower bound above. Consider the test $\tilde \phi: \cX \to \{p, q\}$ that declares $p$ on $R_p$ and $q$ on $R_q \cup \cE_\era$. The $\lambda$-weighted Bayes error for this test is easily seen to be
\begin{align*}
\prior \P_p(R_q) + \overline \prior \P_q(R_p) + \prior \P_p(\cE_\era).
\end{align*}
A similar argument can be made for when $\overline \prior \P_q(\cE_\era) \le \prior \P_p(\cE_\era)$. This shows that for every test with erasures, there exists a test without erasures that is as good or better. Thus, the optimal test can be the best test without erasures, which is simply the maximum-a-posteriori rule. This completes the proof of part (i).

\paragraph{Proof of (ii):}
Suppose a test declares $p$ on a region $R_p \subseteq \cX$, $q$ on a region $R_q \subseteq \cX$ and $\star$ elsewhere, its weighted error term can be expressed as
\begin{align}\label{eq: bayes_two_terms}
\lambda \P(\cE_\tot) + \overline \lambda \P(\cE_\und) &= \left( \lambda \sum_{x \in R_p^c} \prior p(x) + \overline \lambda \sum_{x \in R_p} \overline \prior q(x) \right) + \left( \lambda \sum_{x \in R_q^c} \overline \prior q(x) + \overline \lambda \sum_{x \in R_q} \prior p(x) \right).
\end{align}
We can separately minimise both terms, first one over $R_p$ and second one over $R_q$. If the resulting minimisers are disjoint, then that means we have identified an optimal test. To minimise the first bracketed expression, it is clear that one should set $x \in R_p$ if $\overline \lambda \overline \prior q(x) \le \lambda \prior p(x)$. A similar argument works for the second bracketed expression as well, yielding the rule $x \in R_q$ if $\lambda \overline \prior q(x) \ge \overline \lambda \prior p(x)$ . Since $\lambda < 1/2$, the regions $R_p$ and $R_q$ are disjoint, and hence $\phi$ is a valid test with erasures. This concludes the proof of part (ii).
\end{proof}

\Cref{lemma: np_erasure} shows that we can focus on $\lambda < 1/2$, since the sample complexity in the erasure setting when $\lambda \ge 1/2$ will be identical to that in the standard setting. To tackle the $\lambda < 1/2$ case, we first show an alternate characterisation of the $\lambda$-weighted Bayes error.

\begin{lemma}\label{lemma: bayesian_erasure_error}
Let $p$ and $q$ be supported on a finite set $\cX$. Let the prior distribution on the hypotheses be $(\prior, 1-\prior)$ with $\prior \le 1/2$, and let $\lambda \in (0,1/2)$. Then the following equality holds:
\begin{align*}
\min_{\phi}\P^\phi(\cE_\und) + \lambda \P^\phi(\cE_\era) = \sum_{x \in \cX} \min\{\lambda \prior p(x), \overline \lambda \overline \prior q(x) \} + \sum_{x \in \cX} \min\{\overline \lambda \prior p(x),  \lambda \overline \prior q(x) \},
\end{align*}\
where the minimum is taken over all tests $\phi: \cX \to \{p, q, \star\}$.
\end{lemma}
\begin{proof}
Consider the equality in equation~\eqref{eq: bayes_two_terms}. As $\lambda < 1/2$, the minimum of the left hand side over all tests $\phi$ is the same as the minimum of each of the two terms on the right hand side over the possible decision sets $R_p$ and $R_q$. It is easy to check that 
\begin{align*}
\min_{R_p} \left( \lambda \sum_{x \in R_p^c} \prior p(x) + \overline \lambda \sum_{x \in R_p} \overline \prior q(x) \right) =  \sum_{x \in \cX} \min\{\lambda \prior p(x), \overline \lambda \overline \prior q(x) \},
\end{align*}
and
\begin{align*}
\min_{R_q} \left( \lambda \sum_{x \in R_q^c} \overline \prior q(x) + \overline \lambda \sum_{x \in R_q} \prior p(x) \right) = \sum_{x \in \cX} \min\{\overline \lambda \prior p(x),  \lambda \overline \prior q(x) \}.
\end{align*}
\end{proof}

\begin{remark}
Suppose $P_e(\prior, p,q)$ is the optimal Bayes error for a simple binary hypothesis testing problem \emph{without} erasures, where the prior is $(\prior, \overline \prior)$ and the two distributions are $p$ and $q$ from \Cref{lemma: np_erasure}.
Then the right hand side in the above lemma may be thought of as the weighted error probabilities of two separate binary hypothesis testing problems without erasures:
\begin{align}\label{eq: bayesian_two_tests}
\min_\phi \P^\phi(\cE_\und) + \lambda \P^\phi(\cE_\era) &= (\lambda \prior + \overline \lambda \overline \prior)P_e\left(\frac{\lambda \prior}{\lambda \prior + \overline \lambda \overline \prior}, p , q \right) + (\overline \lambda \prior + \lambda \overline \prior)P_e\left(\frac{\overline \lambda \prior}{\overline \lambda \prior + \lambda \overline \prior}, p , q \right).
\end{align}
\end{remark}

\begin{theorem}\label{thm: erasure_sc_bayesian}
Let $p$ and $q$ be distributions on a finite priorbet $\cX$. Let the prior distribution on the hypotheses be $(\prior, 1-\prior)$ where $\prior \le 1/2$. Let $\lambda >0$ and $\delta \in (0,1)$. The sample complexity of the simple binary hypothesis testing with erasures problem $\cB_\baye(p,q, \prior, \lambda, \delta)$ is given by:
\begin{enumerate}
\item
When $\lambda \ge 1/2$, the sample complexity $\nstar_\baye(p, q, \prior, \lambda, \delta) \asymp \nstar_\bay(p,q,\prior,\delta)$.
\item
Let $\prior< 1/16$ and $\lambda < 1/16$. Define $I_1 = (0, \prior\lambda/4]$ and $I_2 = [2\prior\lambda, \min\{\prior, \lambda\}/8]$.  
\begin{enumerate}
\item When $\prior \le \lambda$, the sample complexity is given by
\begin{align*}
\nstar_\baye(p, q, \prior, \lambda, \delta) \asymp
\max\{\nstar_\bay(p, q, \prior', \delta'), \nstar_\bay(p, q, \prior'', \delta'')\} \quad \text{ if } \delta \in I_1 \cup I_2,
\end{align*}
where $\prior' = \frac{\lambda \prior}{\lambda \prior + \overline \lambda \overline \prior}$, $\delta' = \frac{\delta}{\lambda \prior + \overline \lambda \overline \prior}$, $\prior'' = \frac{\overline \lambda \prior}{\overline \lambda \prior + \lambda \overline \prior}$ and $\delta'' = \frac{\delta}{\overline \lambda \prior + \lambda \overline \prior}$.

\item When $\lambda \le \prior$, the sample complexity is given by
\begin{align*}
\nstar_\baye(p,q,\prior, \lambda, \delta) \asymp \max\{ \nstar_\bay(p,q,\prior', \delta'), \nstar_\bay(q,p,\prior'', \delta'')\}
 \quad \text{ if } \delta \in I_1 \cup I_2,
\end{align*}
where $\prior' = \frac{\lambda \prior}{\lambda \prior + \overline \lambda \overline \prior}$, $\delta' = \frac{\delta}{\lambda \prior + \overline \lambda \overline \prior}$, $\prior'' = \frac{\lambda \overline \prior}{\overline \lambda \prior + \lambda \overline \prior}$ and $\delta'' = \frac{\delta}{\overline \lambda \prior + \lambda \overline \prior}$. (Note that the $\prior''$ here is different from above to ensure that $\prior'' \le 1/2$. Note also that the second terms within the maximisation has the roles of $p$ and $q$ switched.)

\end{enumerate}
\end{enumerate}
\end{theorem}

\begin{remark}
A peculiar point in result $2$ above is that we can identify the sample complexity for all $\delta \le \min\{\prior, \lambda\}/8$, except those in a small range around $\prior\lambda$, specifically, $(\prior\lambda/4, 2\prior\lambda)$. This may not be an artefact of the analysis. The analysis relies on applying \Cref{prop:mild-reduction} which is only valid when $\delta \le \prior/4$ in the $\cB_\bay(p, q, \prior, \delta)$ problem. We believe the proposition is not always valid outside this regime, since the sample complexity $\nstar_\bay(p, q, \prior, \delta)$ exhibits unexpected behaviour when $\delta$ is close to $\prior$. This may lead to some discontinuity in the sample complexity when $\delta \in (\prior\lambda/4, 2\prior\lambda)$.
\end{remark}

\begin{proof}
The $\lambda \ge 1/2$ case is immediate using part 1 of \Cref{lemma: np_erasure}. We consider the two cases:
\begin{enumerate}
\item
$\prior \le \lambda < 1/16$.
\item
$\lambda < \prior < 1/16$.
\end{enumerate}

\paragraph{Case I ($\prior \le \lambda$):} First note that when $\prior \le \lambda < 1/2$, we have
\begin{align*}
\lambda \prior < \overline \lambda\overline \prior, \quad \text{ and } \overline \lambda \prior < \lambda \overline \prior. 
\end{align*}
Using the error probability interpretation in equation~\eqref{eq: bayesian_two_tests}, we see that if $\nstar_\baye(\prior, \lambda, \delta) = n$, then
\begin{align}\label{eq: erasure_error}
(\lambda \prior + \overline \lambda \overline \prior)P_e^{(n)}\left(\frac{\lambda \prior}{\lambda \prior + \overline \lambda \overline \prior}, p , q \right) + (\overline \lambda \prior + \lambda \overline \prior)P_e^{(n)}\left(\frac{\overline \lambda \prior}{\overline \lambda \prior + \lambda \overline \prior}, p , q \right) \le \delta,
\end{align}
and so
\begin{align*}
P_e^{(n)}\left(\frac{\lambda \prior}{\lambda \prior + \overline \lambda \overline \prior}, p , q \right) &\le \frac{\delta}{\lambda \prior + \overline \lambda \overline \prior} \quad \text{ and }\\
P_e^{(n)}\left(\frac{\overline \lambda \prior}{\overline \lambda \prior + \lambda \overline \prior}, p , q \right) &\le \frac{\delta}{\overline \lambda \prior + \lambda \overline \prior}.
\end{align*}
This gives us a lower bound on $n$, specifically,
\begin{align}\label{eq: erasure_lower}
n \ge \max \left\{\nstar_\bay\left(p, q, \frac{\lambda \prior}{\lambda \prior + \overline \lambda \overline \prior}, \frac{\delta}{\lambda \prior + \overline \lambda \overline \prior}\right), \nstar_\bay\left(p, q, \frac{\overline \lambda \prior}{\overline \lambda \prior + \lambda \overline \prior}, \frac{\delta}{\overline \lambda \prior + \lambda \overline \prior} \right) \right\}.
\end{align}
Furthermore, we can also get an upper bound on $n$ by making each of the two terms on the left hand side of \eqref{eq: erasure_error} at most $\delta/2$, giving the bound
\begin{align}\label{eq: erasure_upper}
n \le \max \left\{\nstar_\bay\left(p, q, \frac{\lambda \prior}{\lambda \prior + \overline \lambda \overline \prior}, \frac{\delta}{2(\lambda \prior + \overline \lambda \overline \prior)}\right), \nstar_\bay\left(p, q, \frac{\overline \lambda \prior}{\overline \lambda \prior + \lambda \overline \prior}, \frac{\delta}{2(\overline \lambda \prior + \lambda \overline \prior)} \right) \right\}.
\end{align}
If we can show that the lower bound~\eqref{eq: erasure_lower} matches the upper bound~\eqref{eq: erasure_upper} up to constants, we will be done.

Consider the two intervals: $I_1 = (0, \prior\lambda/4]$ and $I_2 = [2\prior\lambda, \prior/8]$.  (Note that $I_2$ is a valid interval since $\lambda < 1/16$.) We claim that the sample complexity $n$ is determined within constants for $\delta \in I_1 \cup I_2$. This resolves the sample complexity question for almost the entire range of interest, i.e. $(0, \prior/8]$, except for a small range when $\delta \in (\prior\lambda/4, 2\prior\lambda)$. 

Observe that \Cref{prop:mild-reduction} states that $\nstar_\bay(p, q, \prior, \delta) \asymp \nstar_\bay(p, q, \prior, \delta/2)$ as long as $\delta \le \prior/4$. Let $\prior' = \frac{\lambda \prior}{\lambda \prior + \overline \lambda \overline \prior}$, $\delta' = \frac{\delta}{\lambda \prior + \overline \lambda \overline \prior}$, $\prior'' = \frac{\overline \lambda \prior}{\overline \lambda \prior + \lambda \overline \prior}$ and $\delta'' = \frac{\delta}{\overline \lambda \prior + \lambda \overline \prior}$. If $\delta \in I_1$, it is easy to check that $\delta' \le \prior'/4$ and $\delta'' \le \prior''/4$. Hence, we conclude that
\begin{align*}
\nstar_\bay\left(p, q, \prior', \delta'\right) &\asymp  \nstar_\bay\left(p, q, \prior', \delta'/2\right), \quad \text{ and }\\
\nstar_\bay\left(p, q, \prior'', \delta''\right) &\asymp  \nstar_\bay\left(p, q,  \prior'', \delta''/2\right).
\end{align*}
Thus, the sample complexity when $\delta \in I_1$ is simply 
\begin{align*}
\nstar_\baye(p, q, \prior, \lambda, \delta) \asymp \max\{\nstar_\bay(p, q, \prior', \delta'), \nstar_\bay(p, q, \prior'', \delta'')\},
\end{align*}
and both terms on the right hand side can be computed (up to constants) using \Cref{thm:main-result-intro-bay}.

If $\delta \in I_2$, observe that $\nstar_\bay(p, q, \prior', \delta') = \nstar_\bay(p, q, \prior', \delta'/2) = 0$, because $\delta' \ge 2\prior'$. Furthermore, we have
\begin{align*}
\nstar_\bay(p, q, \prior'', \delta'') \asymp \nstar_\bay(p, q, \prior'', \delta''/2)
\end{align*}
as long as $\delta'' \le \prior''/4$, which happens when
\begin{align*}
\delta \le \overline \lambda\prior/4.
\end{align*}
This is easily satisfied since $\delta \le \prior/8$ and $\overline \lambda > 1/2$. This shows that for $\delta \in I_2$, the sample complexity is
\begin{align*}
\nstar_\baye(p, q, \prior, \lambda, \delta) \asymp \nstar_\bay(p, q, \prior'', \delta'').
\end{align*}

\paragraph{Case 2 ($\lambda < \prior$):} 
First note that when $\lambda \le \prior < 1/2$, we have
\begin{align*}
\lambda \prior < \overline \lambda\overline \prior, \quad \text{ and }  \lambda \overline \prior < \overline \lambda  \prior. 
\end{align*}
Using the error probability interpretation in equation~\eqref{eq: bayesian_two_tests}, we see that if $\nstar_\baye(\prior, \lambda, \delta) = n$, then
\begin{align}\label{eq: erasure_error2}
(\lambda \prior + \overline \lambda \overline \prior)P_e^{(n)}\left(\frac{\lambda \prior}{\lambda \prior + \overline \lambda \overline \prior}, p , q \right) + (\overline \lambda \prior + \lambda \overline \prior)P_e^{(n)}\left(\frac{\lambda \overline \prior}{\overline \lambda \prior + \lambda \overline \prior}, q , p \right) \le \delta.
\end{align}
Note that the second term has $p$ and $q$ exchanged to ensure the argument denoting the prior probability in $P^n_e(\cdot)$ is at most half. This means
\begin{align*}
P_e^{(n)}\left(\frac{\lambda \prior}{\lambda \prior + \overline \lambda \overline \prior}, p , q \right) &\le \frac{\delta}{\lambda \prior + \overline \lambda \overline \prior} \quad \text{ and }\\
P_e^{(n)}\left(\frac{ \lambda \overline \prior}{\overline \lambda \prior + \lambda \overline \prior}, q , p \right) &\le \frac{\delta}{\overline \lambda \prior + \lambda \overline \prior}.
\end{align*}
Using a similar argument as in Case 1, we get a lower bound
\begin{align}\label{eq: erasure_lower_bound2}
n &\ge \max \left\{\nstar_\bay\left( p,q,\frac{\lambda \prior}{\lambda \prior + \overline \lambda \overline \prior}, \frac{\delta}{\lambda \prior + \overline \lambda \overline \prior}\right), \nstar_\bay\left(q, p, \frac{ \lambda \overline \prior}{\overline \lambda \prior + \lambda \overline \prior}, \frac{\delta}{\overline \lambda \prior + \lambda \overline \prior} \right) \right\} \nonumber \\
&= \max\{ \nstar_\bay(p,q,\prior', \delta'), \nstar_\bay(q,p,\prior'', \delta'')\}.
\end{align}
and an upper bound
\begin{align}\label{eq: erasure_upper_bound2}
n &\le \max \left\{\nstar_\bay\left(p, q, \frac{\lambda \prior}{\lambda \prior + \overline \lambda \overline \prior}, \frac{\delta}{2(\lambda \prior + \overline \lambda \overline \prior)}\right), \nstar_\bay\left(q, p, \frac{ \lambda \overline \prior}{\overline \lambda \prior + \lambda \overline \prior}, \frac{\delta}{2(\overline \lambda \prior + \lambda \overline \prior)} \right) \right\} \nonumber \\
&= \max\{ \nstar_\bay(p,q,\prior', \delta'/2), \nstar_\bay(q,p,\prior'', \delta''/2)\}.
\end{align}
We now consider the intervals $I_1 = (0, \prior\lambda/4]$ and $I_2 = [2\prior\lambda, \lambda/8]$. Using essentially identical arguments as in Case 1, we note that when $\delta \in (0, \prior\lambda/4]$, the sample complexity is
\begin{align*}
\nstar_\baye(p, q, \prior, \lambda, \delta) \asymp \max\{ \nstar_\bay(p,q,\prior', \delta'), \nstar_\bay(q,p,\prior'', \delta'')\}.
\end{align*}
For the second interval, $\nstar_\bay(p,q,\prior', \delta') = \nstar_\bay(p,q,\prior', \delta'/2) = 0$, and so only the second term matters. Here, we are able to apply \Cref{prop:mild-reduction} because 
\begin{align*}
\delta \le \lambda \overline \prior/4,
\end{align*}
which is true since $\delta \le \lambda/8 \le \lambda \overline \prior/4.$ Thus, when $\delta \in I_2$, we get the result
\begin{align*}
\nstar_\baye(p, q, \prior, \lambda, \delta) \asymp \nstar_\bay(q,p, \prior'', \delta'').
\end{align*}

This concludes the proof of all cases.
\end{proof}

\section{Discussion}
We studied simple binary hypothesis testing from the sample complexity perspective, which is a non-asymptotic regime that had received scant attention in prior literature. We identified the sample complexity (up to universal multiplicative constants) in the Bayesian and prior-free formulations of the problem, in terms of $f$-divergences from the Hellinger and Jensen--Shannon families. We also addressed the sample complexity of Bayesian and prior-free when erasures are permitted, and we made some observations about the sample complexity in the sequential setting. There are many interesting research directions that emerge from our work. Characterising the sample complexity in the weak detection regime (cf. \Cref{sub:results-large-failure}) remains an open problem. Indeed, even in the simple setting of  $p = \cN(-1,1)$ and $q = \cN(1,1)$ and $\prior < 1/2$, the sample complexity (computed using a computer) demonstrates unusual behaviour with respect to $\gamma$. In particular, it does depend on $\gamma$ (i.e., the upper bound in \Cref{thm:bayes-all-regimes} is loose), but the dependence on $\gamma$ is super-linear (i.e., the lower bound in \Cref{thm:bayes-all-regimes} is also loose). The Hellinger characterisation of the sample complexity of simple binary hypothesis testing appears in Le Cam's bound in minimax optimal estimation~\cite{Lecam73}, in rates of convergence for locally minimax optimal estimation~\cite{DonL91}, and in communication complexity lower bounds in computer science theory~\cite{BravGMNW16}. It would be interesting to see if the skewed $\js{\prior}$ characterization of the sample complexity presented in this work has similar applications. Finding the sample complexity in the sequential setting is another interesting direction to pursue.

\section*{Acknowledgements}
Part of this work was done when VJ and PL were at semester-long programs at the Simons Institute for the Theory of Computing and the Simons--Laufer Mathematical Sciences Institute. VJ and PL thank these institutions for their hospitality.

\printbibliography

\section{Auxiliary results} %
\label{app:omitted_details_pfht}

\subsection{Proof of \Cref{fact:boost}}
\label{app:proof_of_cref_fact_boost}
We will use the following form of Bennet's inequality:
\begin{fact}[{Concentration of binomial: Bennet's inequality~\cite[Theorem 2.9.2]{Vershynin18}}]
\label{fact:binomial}
Let $X_1,\dots,X_m$ be $m$ i.i.d.\ Bernoulli random variables with bias $p > 0$.
For $u \geq 0$, we have
\begin{align*}
\P\left\{ \left( \frac{1}{m} \sum_{i=1}^m X_i \right)  \geq p (1 + u) \right \} &\leq \P\left\{ \left( \frac{1}{m} \sum_{i=1}^m X_i \right)  \geq mp (1 + u(1-p))   \right \} \\
&\leq \exp \left(  -mp (1-p) f(u)\right),
\end{align*}
where $f(u) = (1+u) \log (1+u) - u$.
In particular, if the bias $p$ satisfies $p \leq 0.5$, and $t \geq 1$, we have\footnote{This follows by noting that $(1-p) \geq 0.5$ and $(1+u) \log(1+u) - u \geq 0.5 u \log(1+u)$. Indeed, the derivative of the second expression is $0.5\left( \log(1+u) - \frac{u}{u+1}\right)$, which is nonnegative. }
\begin{align}
\label{eq:bennet-simple}
\P\left\{ \left( \frac{1}{m} \sum_{i=1}^m X_i \right)  \geq t p\right \}  \leq \exp \left(  - (1/4) mp (t-1) \log t\right) \,.
\end{align}
\end{fact}

We provide the proof of \Cref{fact:boost} below for completeness:
\begin{proof}[Proof of \Cref{fact:boost}]
Since the buckets are independent and $\tau \leq 1/4$, the probability that the fraction of the buckets reporting wrong outcomes is at least $1/2$ can be upper-bounded using \Cref{eq:bennet-simple} (by taking $m = T$, $p = \tau$, $t = \frac{0.5}{\tau}$) 
by 
\begin{align*}
 &\exp\left( - (1/4)(T)(\tau)((1/2\tau) - 1)  \log\left( \frac{1}{2\tau} \right)  \right) \leq \exp\left( - (1/4) T (1/4 \tau) ((1/2)\log(1/\tau)) \right) \\
 & \quad \leq \exp(- (T \log (1/\tau))/2^5)\,.
 \end{align*}
\end{proof}

\subsection{Proof of \Cref{cl:linear-vs-nearly-linear}}
\label{app:proof-of-cl-linear-vs-nearly-linear}
\ClLinearVsNearlyLinear*
\begin{proof}
We consider two cases.
For all $x$ such that $\prior x \leq 1$, we have
\begin{align*}
1 + \frac{x}{1 + \prior x} \leq 1 + x = \left( 1 + x \right)^{ \bl_*} \left( 1 + x \right)^{\lambda_*} \leq \left( 1 + x \right)^{\bl_*} \left( 1 + \frac{1}{\prior} \right)^{\lambda_*}.
\end{align*}
Thus, it remains to control the overhead term $\left( 1 +
\frac{1}{\prior} \right)^{\lambda_*}$, which we will show is
at most $e^{2r}$.
We first rewrite it as $\exp\left( \lambda_* \log\left( 1 +
\frac{1}{\prior} \right)\right)$.
For the exponent, we have
\begin{equation*}
\lambda_* \log(1+ 1/\prior)
\leq \frac{r}{\log(1/\prior)} \cdot 2\log(1/\prior) \leq
2r\,,
\end{equation*}
where we use that $\log((1 +\prior)/\prior) \leq 2 \log(1/\prior)$.
Thus, we have $1 + \frac{x}{1 + \prior x} \leq e^{2r} (1 +
x)^{\bl_*}$.

We now consider the case when $x > 1/\prior$.
Following similar steps, we arrive at the desired conclusion:
\begin{align*}
1 + \frac{x}{1 + \prior x} \leq 1 + \frac{x}{\prior x}=  1 + \frac{1}{\prior} =  \left(1 + \frac{1}{\prior}\right)^{\bl_*} \left(  1+ \frac{1}{\prior}\right)^{\lambda_*} \leq e^{2r} \left(1 + \frac{1}{\prior} \right)^{\bl_*} \leq e^{2r}(1+ x)^{\bl_*}.
\end{align*}
\end{proof}

\section{Asymptotic analyses of the errors in hypothesis testing with erasures}\label{app: erasure}

In this section, we analyse how the probabilities of undetected error and erasure error behave in the asymptotic regime of sample size $n$ tending to infinity. We could not find these results in the literature and they may be of independent interest.
 
There are three different asymptotic regimes we consider: the Stein regime, the Chernoff regime, and the Bayes regime. The regime names and their analyses closely follow the proof techniques in~\cite{PolWu23} for the simple binary hypothesis testing setting without erasures, which we shall refer to as the \emph{standard setting}. The simple binary hypothesis testing with erasures will be referred to as the \emph{erasure setting}.

We recap known results in the standard setting:
\begin{enumerate}
\item
\emph{Stein regime:} When the type-I error is at most $\epsilon$, the type-II error decays exponentially fast with exponent $\kl(p,q)$. When the type-II error is at most $\epsilon$, the type-I error decays exponentially fast with exponent $\kl(q,p)$.
\item \emph{Chernoff regime:} Let $A$ and $B$ be the exponential decay rates of the type-I and the type-II errors respectively. Then the set of all achievable $(A,B)$ is 
\begin{align*}
\{(A_\lambda, B_\lambda) : A_\lambda = \kl(p_\lambda,q), B_\lambda = \kl(p_\lambda, p), \lambda \in [0,1]\}
\end{align*}
where $p_\lambda$ is the tilted distributed $p_\lambda(x) \propto p(x)^{\overline \lambda}q(x)^\lambda$.
\item \emph{Bayes regime:} When the prior distribution on the hypothesis is $(\prior, 1-\prior)$ for $\prior \in (0,1)$, the exponential rate of decay of the Bayes error (with is a weighted sum of the type-I and type-II errors) is given by the Chernoff information between $p$ and $q$, which equals
\begin{align*}
\mathrm{CI}(p,q) = -\inf_{\lambda \in [0,1]} \log \sum_{x \in \cX} p(x)^{1-\lambda}q(x)^\lambda.
\end{align*}
Chernoff information also equals $A_{\lambda^*}$ from the Chernoff regime where $\lambda^*$ is such that $A_{\lambda^*} = B_{\lambda^*}$.
\end{enumerate}

In the erasure setting, if the prior is $(\prior, 1-\prior)$, the two errors of interest are:
\begin{align*}
\P(\cE_\era) &= \P(\phi(X_1,\dots, X_n) = \star) = \prior \P_p(\phi(X_1,\dots, X_n) = \star) + \overline \prior \P_q(\phi(X_1,\dots, X_n) = \star),
\end{align*}
and
\begin{align*}
\P(\cE_\und) &= \prior \P_p(\phi(X_1,\dots, X_n) = q) + \overline \prior \P_q(\phi(X_1,\dots, X_n) = p).
\end{align*}
The questions to be addressed in the three asymptotic regimes are as follows:
\begin{enumerate}
\item
\emph{Stein regime:} When $\P(\cE_\era)$ is at most $\epsilon$, what is the exponential rate of decay of $\P(\cE_\und)$?
\item
\emph{Chernoff regime:} When both $\P(\cE_\era)$ and $\cP(\cE_\und)$ are decaying exponentially fast, say with rates $A$ and $B$, what are all pairs of achievable rates $(A,B)$?
\item \emph{Bayes regime:} What is the exponential rate of decay of $\lambda \P(\cE_\era) +  \P(\cE_\und)$, which corresponds to the weighted sum $\lambda \P(\cE_\tot) + \overline \lambda \P(\cE_\und)$? 
\end{enumerate}

Thus, the main difference between the standard and the erasure settings is that the role of type-I and type-II errors is taken up by the erasure and undetected errors. The latter two errors are not symmetric, leading to a few differences between the two settings. First, it is necessary to specify the prior distribution on the hypotheses to define the errors in the erasure setting, and so the Stein and Chernoff regimes are not ``prior-free'' as in the standard setting. And second, in the Stein regime, it only makes sense to bound the erasure error probability and ask for the decay rate of the undetected error probability; the counterpart of bounding the undetected error probability and studying the decay rate of the erasure error probability does not make sense. This is because we can simply refrain from declaring erasures, thus ensuring an erasure error probability of 0, while making the undetected error probability arbitrarily small for all large $n$ (note that it decays like $e^{-n(\mathrm{CI}(p,q)+o(1))}$).

\subsection{Stein regime in the erasure setting}\label{subsection: stein}

\begin{theorem}\label{thm: stein_erasure}
Let $p$ and $q$ be two distributions on a finite domain $\cX$. Let the prior distribution on the hypothesis $\theta$ be $(\prior, \overline \prior)$ for $\prior \in (0,1)$. Consider a (possibly randomised) hypothesis test $\phi$ that takes as input $n$ i.i.d.\ samples $X_1, \dots, X_n$ distributed according to $p^{\otimes n}$ if $\theta = p$ and  $q^{\otimes n}$ if $\theta = q$, and produces an estimate $\phi(X_1, \dots, X_n) \in \{p, q, \star\}$. Let $\cE_\und^n$ be the event of an undetected error and $\cE_\era^n$ be the event on an erasure error when using $n$ observations. The $\epsilon$-optimal exponent in Stein's regime is defined as
\begin{align*}
V_\epsilon := \sup\{E ~:~ \exists n_0, \forall n > n_0, \exists \phi \text{ such that } \P(\cE_\era^n) < \epsilon, \cP(\cE_\und^n) < e^{-nE}\}. 
\end{align*}
Then the Stein exponent is given by
\begin{align*}
V_\epsilon = 
\begin{cases}
\min\{ \kl(p,q), \kl(q,p) \} \quad &\text{ if } \epsilon \in (0, \prior),\\
\max\{ \kl(p,q), \kl(q,p) \} \quad &\text{ if } \epsilon \in (\prior,1).
\end{cases}
\end{align*}
\end{theorem}

Observe that in the Stein regime, if erasures are not allowed then the probability of undetected error decays like $e^{-n(\mathrm{CI}(p,q)+o(1))}$. Since $\mathrm{CI}(p,q) \le \min\{\kl(p,q), \kl(q,p)\}$, it is evident that the flexibility of declaring erasures boosts the error exponent even for very small $\epsilon$. As $\epsilon$ increases beyond $\prior$, we observe a threshold phenomenon wherein the exponent jumps from the minimum of the two KL divergences to their maximum. It is interesting to note that this phenomenon has no counterpart in the standard setting.

We begin by collecting structural results for the erasure setting.
\begin{lemma}\label{lemma: np_erasure_app}
Let $p$ and $q$ be two distributions on a finite set $\cX$. Let $0 < \lambda_1 < \lambda_2$ and consider the weighted error term $\lambda_1 \P(\cE_\tot) + \lambda_2 \P(\cE_\und)$, which equals $\lambda_1 \P(\cE_\era) + (\lambda_1+\lambda_2)\P(\cE_\und)$. Then the following test $\phi$ minimises this weighted error:
\begin{align*}
\phi(x) = 
\begin{cases}
p \quad &\text{ if } \frac{\prior p(x)}{\overline \prior q(x)} \ge \frac{\lambda_2}{\lambda_1}\\
q \quad &\text{ if } \frac{\prior p(x)}{\overline \prior q(x)} \le \frac{\lambda_1}{\lambda_2}\\
\star &\text{ otherwise.}
\end{cases}
\end{align*}
\end{lemma}
\begin{proof}
The proof is essentially identical to \Cref{lemma: np_erasure}.
\end{proof}

We now show one-shot achievability and converse results for the erasure setting. To facilitate the discussion, define the function $T: \cX \to \real$ as
\begin{align*}
T(x) := \log \frac{\prior p(x)}{\overline \prior q(x)}.
\end{align*}
In what follows, we shall always choose $\lambda_1 = 1$ and $\lambda := \lambda_2 > 1$. For this choice, the optimal test in \Cref{lemma: np_erasure_app} gives the regions $R_p = \{x : T(x) \ge \log \lambda \}$ and $R_q = \{x : T(x) \le -\log \lambda\}$.

\begin{lemma}[Achievability]\label{lemma: stein: achievability}
Let $\tau > 0$. There exists a hypothesis test with erasures such that
\begin{align*}
\P(\cE_\era) = \P(T \in (-\tau, \tau)), \quad \text{ and } \quad
\P(\cE_\und) \le e^{-\tau}. 
\end{align*}
\end{lemma}
\begin{proof}
Consider the optimal test in \Cref{lemma: np_erasure} with $\lambda_1 = 1$ and $\lambda_2 = e^\tau$. This test has an erasure probability exactly as stated in \Cref{lemma: stein: achievability}, and its undetected error probability may be upper bounded as:
\begin{align*}
\P(\cE_\und) &= \prior \P_p(T < -\tau) + \overline \prior \P_q(T > \tau)\\
&= \sum_{x:T(x) < -\tau} \prior p(x) + \sum_{x:T(x) > \tau}\overline \prior q(x)\\
&\le \sum_{x:T(x) < -\tau}\overline \prior q(x) e^{-\tau} + \sum_{x:T(x) > \tau} e^{-\tau} \prior p(x)\\
&\le \sum_{x}\overline \prior q(x) e^{-\tau} + \sum_{x} e^{-\tau} \prior p(x)\\
&= e^{-\tau}.
\end{align*}
\end{proof}

\begin{lemma}[Converse]\label{lemma: stein: converse}
Consider any hypothesis test with erasures $\phi$. Set $\lambda_1 = 1$ and $\lambda_2 = \lambda > 1.$ Consider the weighted error term $ \P(\cE_\tot) + \lambda \P(\cE_\und)$, which equals $ \P(\cE_\era) + (1+\lambda)\P(\cE_\und)$. The following lower bound holds for this weighted error:
\begin{align*}
\P(\cE_\era) + (1+\lambda)\P(\cE_\und) \ge \P\left( T \in \left(-\log \lambda, \log \lambda\right)\right).
\end{align*}
\end{lemma}
\begin{proof}
We know that the weighted error of $\phi$ is at least as large as that of the optimal test $\phi^*$ in \Cref{lemma: np_erasure_app}. Let the erasure and undetected probabilities under the optimal test be $\P^{\phi^*}(\cE_\era)$ and $P^{\phi^*}(\cE_\und)$, respectively. Then
\begin{align*}
\P(\cE_\era) + (1+\lambda)\P(\cE_\und) &\ge \P^{\phi^*}(\cE_\era) + (1+\lambda)\P^{\phi^*}(\cE_\und)\\
&\ge \P^{\phi^*}(\cE_\era)\\
&= \P\left( T \in \left(-\log \lambda, \log \lambda\right)\right).
\end{align*}
\end{proof}
We now provide the proof of \Cref{thm: stein_erasure}.
\begin{proof}[Proof of \Cref{thm: stein_erasure}]
Without loss of generality, assume $\kl(p,q) < \kl(q,p)$. 

\paragraph{Case I: $\epsilon < \prior$.} Let $\delta > 0$. We first show that $V_\epsilon \ge \kl(p,q) - \delta$. Consider the test $\phi$ with $\tau = n(\kl(p,q) - \delta)$ in \Cref{lemma: stein: achievability}. Note that $T(X_1, \dots, X_n) = \log \frac{\prior}{\bar \prior} + \sum_{i=1}^n \log\frac{p(X_i)}{q(X_i)}$. Thus, we may apply the weak law of large numbers, we have
\begin{align*}
\P_p(\phi(X_1, \dots, X_n) = \star) &\le \P_p\left( T(X_1, \dots, X_n) < n(\kl(p,q) - \delta)\right) \to 0, \quad \text{ and }\\
\P_q(\phi(X_1, \dots, X_n) = \star) &\le \P_q\left( T(X_1, \dots, X_n) > -n(\kl(p,q) - \delta)\right) \to 0.
\end{align*}
Hence, for all large enough $n$ the erasure probability $\P(\cE_\era^n)$ is at most $\epsilon$. By \Cref{lemma: stein: achievability}, the probability of undetected error is at most 
\begin{align*}
\P(\cE_\und^n) \le e^{-\tau} = e^{-n(\kl(p,q) - \delta)}.
\end{align*}
This shows that $V_\epsilon \ge \kl(p,q) - \delta$. As this holds for any choice of $\delta$, we conclude that $V_\epsilon \ge \kl(p,q).$

To show the converse, suppose $E < V_\epsilon$ is any achievable exponent. This means for all large enough $n$, there exist tests that guarantee $\P(\cE_\era^n) < \epsilon$ and $\P(\cE_\und^n) < e^{-nE}$. For $0 < \delta < \kl(q,p)-\kl(p,q)$, apply the converse result in \Cref{lemma: stein: converse} with $\log \lambda = n(\kl(p,q) + \delta)$ and conclude
\begin{align*}
\epsilon + (1+e^{n(\kl(p,q) + \delta)}) e^{-nE} &\ge \P(\cE_\era) + (1+e^{n(\kl(p,q) + \delta)})\P(\cE_\und)\\
&\ge \P\left( T \in \left(-n(\kl(p,q) + \delta), n(\kl(p,q) + \delta)\right)\right).
\end{align*}
We claim that as $n$ tends to infinity, the right hand side tends to $\prior$. This is because 
\begin{align*}
&\lim_{n \to \infty} \P\left( T \in \left(-n(\kl(p,q) + \delta), n(\kl(p,q) + \delta)\right)\right)\\
&= \lim_{n \to \infty} \prior \P_p\left( T \in \left(-n(\kl(p,q) + \delta), n(\kl(p,q) + \delta)\right)\right)\\
&\qquad+\lim_{n \to \infty} \overline \prior\P_q\left( T \in \left(-n(\kl(p,q) + \delta), n(\kl(p,q) + \delta)\right)\right)\\
&= \prior. \tag{WLLN and $\delta$ being small enough}
\end{align*}
Thus,
\begin{align*}
\epsilon + \lim_{n \to \infty} e^{n(\kl(p,q)+\delta - E)} \ge \prior.
\end{align*}
Since $\epsilon < \prior$, the only way the above inequality can hold is if $E \le \kl(p,q)+\delta$. As this holds for all achievable rates $E$, we conclude $V_\epsilon \le \kl(p,q)+\delta$. Taking $\delta \to 0$, we conclude $V_\epsilon \le \kl(p,q)$. Combining the upper and lower bounds on $V_\epsilon$ we conclude $V_\epsilon = \kl(p,q)$.

\paragraph{Case II: $\epsilon > \prior$.} 
Let $\delta > 0$. We first show that $V_\epsilon \ge \kl(q,p) - \delta$. Consider the test $\phi$ with $\tau = n(\kl(q,p) - \delta)$ in \Cref{lemma: stein: achievability}. By the weak law of large numbers, we have
\begin{align*}
\P_p(\phi(X_1, \dots, X_n) = \star) &\le \P_p\left( \sum_{i=1}^n T(X_i) < n(\kl(q,p) - \delta)\right) \to 1, \quad \text{ and }\\
\P_q(\phi(X_1, \dots, X_n) = \star) &\le \P_q\left( \sum_{i=1}^n T(X_i) > -n(\kl(q,p) - \delta)\right) \to 0.
\end{align*}
Since $\prior < \epsilon$, we conclude that for all large enough $n$ the erasure probability $\P(\cE_\era^n)$ is at most $\epsilon$. By \Cref{lemma: stein: achievability}, the probability of undetected error is at most 
\begin{align*}
\P(\cE_\und^n) \le e^{-\tau} = e^{-n(\kl(q,p) - \delta)}.
\end{align*}
This shows that $V_\epsilon \ge \kl(p,q) - \delta$. As this holds for any choice of $\delta$, we conclude that $V_\epsilon \ge \kl(q,p).$

To show the converse, suppose $E < V_\epsilon$ is any achievable exponent. This means for all large enough $n$, there exist tests that guarantee $\P(\cE_\era^n) < \epsilon$ and $\P(\cE_\und^n) < e^{-nE}$. For $0 < \delta$, apply the converse result in \Cref{lemma: stein: converse} with $\log \lambda = n(\kl(q,p) + \delta)$ and conclude
\begin{align*}
\epsilon + (1+e^{n(\kl(q,p) + \delta)}) e^{-nE} &\ge \P(\cE_\era) + (1+e^{n(\kl(q,p) + \delta)})\P(\cE_\und)\\
&\ge \P\left( T \in \left(-n(\kl(q,p) + \delta), n(\kl(q,p) + \delta)\right)\right).
\end{align*}
We claim that as $n$ tends to infinity, the right hand side tends to $1$. This is because 
\begin{align*}
&\lim_{n \to \infty} \P\left( T \in \left(-n(\kl(q,p) + \delta), n(\kl(q,p) + \delta)\right)\right)\\
&= \lim_{n \to \infty} \prior \P_p\left( T \in \left(-n(\kl(q,p) + \delta), n(\kl(q,p) + \delta)\right)\right)\\
&+\lim_{n \to \infty} \overline \prior\P_q\left( T \in \left(-n(\kl(q,p) + \delta), n(\kl(q,p) + \delta)\right)\right)\\
&= \prior + \overline \prior  \tag{WLLN}\\
&= 1.
\end{align*}
Thus,
\begin{align*}
\epsilon + \lim_{n \to \infty} e^{n(\kl(q,p)+\delta - E)} \ge 1.
\end{align*}
Since $\epsilon < 1$, the only way the above inequality can hold is if $E \le \kl(q,p)+\delta$. As this holds for all achievable rates $E$, we conclude $V_\epsilon \le \kl(q,p)+\delta$. Taking $\delta \to 0$, we conclude $V_\epsilon \le \kl(q,p)$. Combining the upper and lower bounds on $V_\epsilon$ we conclude $V_\epsilon = \kl(q,p)$.
\end{proof}

\subsection{Chernoff and Bayes regimes in the erasure setting}

The results in these two regimes rely on the achievability and converse results in ~\Cref{subsection: stein} and Cramer's theorem from large deviation theory. Throughout this section, we shall assume that $p$ and $q$ are absolutely continuous with respect to each other. Since they are assumed to have discrete support, the moment generating function of $L(X) = \log p(X)/q(X)$ is finite for all $t \in \real$, when $X \sim p$ or $X \sim q$. We state Cramer's theorem first:

\begin{theorem}[Cramer's theorem]\label{thm: cramer}
The logarithmic moment generating function of a random variable is defined as:
\[
\Lambda(t) = \log \mathbb{E} \left[ \exp(t X_1) \right].
\]
Let \( X_1, X_2, \dots \) be a sequence of i.i.d.\ random variables with finite logarithmic moment generating function, i.e., \( \Lambda(t) < \infty \) for all \( t \in \mathbb{R} \). Then the Legendre transform of \( \Lambda \), defined by
\[
\Lambda^*(x) := \sup_{t \in \mathbb{R}} \left( t x - \Lambda(t) \right)
\]
and satisfies:
\[
\lim_{n \to \infty} \frac{1}{n} \log \left( \P\left( \sum_{i=1}^{n} X_i \geq n x \right) \right) = - \Lambda^*(x)
\]
for all \( x > \mathbb{E}[X_1] \).
\end{theorem}

Our goal is to identify the pairs of possible exponential decay rates of $\P(\cE_\era^n)$ and $\P(\cE_\und^n)$, denoted by $A$ and $B$. Note that the decay rate of $\P(\cE_\tot^n)$ is the minimum of these two. Define
\[
B^*(A) \triangleq \sup \left\{ B : \exists n_0, \forall n \geq n_0, \exists \phi \text{ such that } \P(\cE_\era^n) < \exp(-n A), \P(\cE_\und^n) < \exp(-n B) \right\}.
\]

\begin{theorem}[Chernoff regime]\label{thm: chernoff_erasure}
Let $p$ and $q$ be supported on a finite set $\cX$, and assume they are absolutely continuous with respect to each other. For $t \in \real$, define 
\begin{align*}
\psi(t)&= \log \sum_{x \in \cX} p(x)^{1-t}q(x)^{t}.\end{align*}
Let $\psi^*(\cdot)$ denote the Legendre transforms of $\psi$. Then the upper boundary of the set of all achievable pairs $(A,B)$ is characterized by 
\begin{align*}
A(\theta) &= \min\{\psi^*(-\theta), \psi^*(\theta)-\theta\},\\
B(\theta) &= \min\{\psi^*(\theta), \psi^*(-\theta)+\theta\} = A(\theta) + \theta,
\end{align*}
for $\theta \in [0, \min\{\kl(q,p), \kl(p,q)\}]$. 
\end{theorem}

The following corollary is an immediate consequence of \Cref{thm: chernoff_erasure}:
\begin{corollary}\label{cor: Bayes}
Let the optimal Bayes error exponent be \begin{align*}
C_\lambda := \{E~:~\exists n_0, \forall n > n_0, \exists \phi \text{ such that } \P(\cE_\era^n) + \lambda \P(\cE_\und^n) < e^{-nE}\}.
\end{align*}
Then for all $\lambda >0$, $$C_\lambda = \max_{\theta \in [0, \min\{\kl(p,q), \kl(q,p)\}]} \min\{A(\theta), B(\theta)\}.$$
\end{corollary}

We note that the expressions for $A(\theta)$ and $B(\theta)$ have been derived in Sason~\cite{Sas11}, but only in the achievable sense, also using Cramer's theorem. However the converse part, which states that these pairs of $(A(\theta), B(\theta))$ form the upper boundary of the achievable region, was not established in prior works.

\begin{proof}[Proof of \Cref{thm: chernoff_erasure}]
Let $\theta \in [0, \min\{\kl(p,q), \kl(q,p)\}]$. Consider the following test $\phi$:
\begin{align*}
\phi(x_1, \dots, x_n) &= 
\begin{cases}
p \quad \text{ if } \sum_{i=1}^n \log p(x_i)/q(x_i) \ge n\theta,\\
q \quad \text{ if } \sum_{i=1}^n \log p(x_i)/q(x_i) \le -n\theta,\\
\star \quad \text{ if } \sum_{i=1}^n \log p(x_i)/q(x_i) \in (-n\theta,  n\theta).\\
\end{cases}
\end{align*}
The probabilities of errors are given by
\begin{align*}
\P(\cE_\era^n) = \prior \P_p\left(\sum_{i=1}^n \log p(X_i)/q(X_i) \in (-n\theta, n\theta) \right) + \overline \prior \P_q\left(\sum_{i=1}^n \log p(X_i)/q(X_i) \in (-n\theta, n\theta) \right),
\end{align*}
and 
\begin{align*}
\P(\cE_\und^n) = \prior \P_p\left(\sum_{i=1}^n \log p(X_i)/q(X_i) \le -n\theta \right) + \overline \prior \P_q\left(\sum_{i=1}^n \log p(X_i)/q(X_i) \ge -n\theta \right).
\end{align*}
A straightforward application of Cramer's theorem yields 
\begin{align*}
\lim_{n \to \infty} -\frac{1}{n} \log \P(\cE_\era^n) &= A(\theta), \quad \text{ and }\\
\lim_{n \to \infty} -\frac{1}{n} \log \P(\cE_\und^n) &= B(\theta).
\end{align*}
This concludes the achievabilty part of the proof.

To show the converse, suppose $\P(\cE_\era^n) < e^{-nA'}$ and $\P(\cE_\und^n) < e^{-nB'}$ for some $A'$ and $B'$, for all large enough $n$. Applying \Cref{lemma: stein: converse},
\begin{align*}
e^{-nA'} + (1+\lambda)e^{-nB'} &> \P(\cE_\era) + (1+\lambda)\P(\cE_\und)\\
&\ge \P\left( \log \frac{\prior p(X_1, \dots, X_n)}{\overline \prior q(X_1, \dots, X_n)} \in \left(-\log \lambda, \log \lambda\right)\right)\\
&= \P\left( \sum_{i=1}^n \log \frac{p(X_i)}{q(X_i)}\in \left(-\log \lambda - \log (\prior/\overline \prior), \log \lambda - \log (\prior/\overline \prior)\right)\right).
\end{align*} 
Choosing $\lambda = e^{n\theta}$ and using Cramer's theorem again, we derive the lower bound
\begin{align*}
e^{-nA'} + e^{-n(B'-\theta)} &> \prior \P_p \left( \sum_{i=1}^n \log \frac{p(X_i)}{q(X_i)}\in \left(-n\theta - \log \prior/\overline \prior, n\theta - \log \prior/\overline \prior\right)\right)\\
& + \overline \prior \P_q \left( \sum_{i=1}^n \log \frac{p(X_i)}{q(X_i)}\in \left(-n\theta - \log \prior/\overline \prior, n\theta - \log \prior/\overline \prior\right)\right)\\
&= e^{n(A(\theta)+o(1))}.
\end{align*}
The only way this inequality can hold for all large $n$ is when either $A' \le A(\theta)$ or $B' \le A(\theta)+\theta = B(\theta)$, and concedes the proof of the converse.
\end{proof}

\end{document}